\documentclass[a4paper]{amsart}

\usepackage{amsmath,amsthm,amssymb}
\usepackage{graphicx}

\title[Isothermic surfaces]{Isothermic surfaces: conformal geometry,\\
Clifford algebras and integrable systems}
\author{F.E. Burstall}
\address{Department of Mathematical Sciences\\
University of Bath\\Bath BA2 7AY\\UK}
\email{feb@maths.bath.ac.uk}

\renewcommand{\O}[1][n+1,1]{\mathrm{O}(#1)}
\newcommand{\GL}{\mathrm{GL}}
\newcommand{\SO}[1][n+1,1]{\mathrm{SO}(#1)}
\newcommand{\Op}{\mathrm{O}^+(n+1,1)}
\renewcommand{\o}[1][n+1,1]{\mathfrak{o}(#1)}
\newcommand{\SL}[1][\Gamma_n]{\mathrm{SL}(#1)}
\newcommand{\Pin}[1][n+1,1]{\mathrm{Pin}(#1)}
\newcommand{\Spin}[1][n+1,1]{\mathrm{Spin}(#1)}
\renewcommand{\L}{\mathcal{L}}
\renewcommand{\P}{\mathbb{P}} 
\newcommand{\PL}{\mathbb{P}(\L)}
\renewcommand{\lor}{\R^{n+1,1}}
\renewcommand{\d}{\mathrm{d}}
\newcommand{\D}{\mathcal{D}}
\newcommand{\R}{\mathbb{R}}
\newcommand{\T}{\mathcal{T}}
\newcommand{\C}{\mathbb{C}\,}

\newcommand{\Z}{\mathbb{Z}}
\renewcommand{\H}{\mathbb{H}}
\newcommand{\g}{\mathfrak{g}}
\renewcommand{\k}{\mathfrak{k}}
\newcommand{\p}{\mathfrak{p}}
\renewcommand{\a}{\mathfrak{a}}
\newcommand{\del}{\partial}
\newcommand{\delbar}{\overline{\partial}}
\newcommand{\Cl}{C\ell}
\newcommand{\half}{\tfrac{1}{2}}
\newcommand{\mcv}{\mathbf{H}}
\newcommand{\<}{\langle}
\renewcommand{\>}{\rangle}
\newcommand{\Mob}[1][n]{\mathrm{M\ddot{o}b}(#1)}
\newcommand{\gi}{g^{-1}}
\newcommand{\tF}{\widetilde{F}}
\newcommand{\mc}[1]{#1^{-1}\d #1}
\newcommand{\abs}[1]{\lvert#1\rvert}
\newcommand{\cross}{g^{\vphantom{1}}_1g_{12}^{-1}g^{\vphantom{1}}_{21}g_2^{-1}}
\newcommand{\norm}[1]{\lVert#1\rVert}
\newcommand{\sv}{\<v\>}
\renewcommand{\sf}{\<f\>}
\newcommand{\Nf}{\mathcal{N}_{\sf}}
\newcommand{\dt}[1][0]{\left.\frac{\d}{\d t}\right|_{t=#1}}
\newcommand{\dl}[1]{\left.\frac{\del #1}{\del\lambda}\right|_{\lambda=0}}
\newcommand{\Gp}{\mathcal{G}^+}
\newcommand{\Gm}{\mathcal{G}^-}
\newcommand{\Gmb}{\mathcal{G}^-_*}
\newcommand{\G}{\mathcal{G}}
\newcommand{\Ug}[1][-]{\mathcal{U}_{g_{#1}}}
\newcommand{\I}{\sqrt{-1}}
\newcommand{\Um}{\Omega_{\mathrm{alg}}G}
\newcommand{\Up}{\Lambda^+G^\C}
\DeclareMathOperator{\rank}{rank}
\DeclareMathOperator{\sign}{sign}
\DeclareMathOperator{\Ad}{Ad}
\DeclareMathOperator{\tAd}{\widetilde{Ad}}
\DeclareMathOperator{\Hom}{Hom}

\DeclareMathOperator{\trace}{trace}
\DeclareMathOperator{\im}{Im}
\DeclareMathOperator{\Map}{Map}
\DeclareMathOperator{\dom}{dom}
\newcommand{\fh}{\hat{f}}
\newcommand{\hF}{\hat{F}}
\newcommand{\hN}{\hat{N}}
\newcommand{\vh}{\hat{v}}
\newcommand{\wh}{\hat{v}'}
\newcommand{\ha}{\hat{\alpha}}
\newcommand{\hb}{\hat{\beta}}
\newcommand{\hg}{\hat{\gamma}}
\newcommand{\hd}{\hat{\delta}}
\newcommand{\gh}{\hat{g}}
\newcommand{\hh}{\hat{h}}
\newcommand{\hpi}{\hat{\pi}}
\newcommand{\hL}{\hat{L}}
\newcommand{\zbar}{\bar{z}}
\newcommand{\bl}{\bar{\lambda}}

\theoremstyle{plain}
\newtheorem{thm}{Theorem}[section]
\newtheorem{prop}[thm]{Proposition}
\newtheorem{lem}[thm]{Lemma}
\newtheorem{cor}[thm]{Corollary}
\newtheorem*{fact}{Fact}
\newtheorem*{perm}{Bianchi Permutability Theorem}
\newtheorem*{prob}{Problem}
\theoremstyle{definition}
\newtheorem*{notation}{Notation}
\newtheorem*{rem}{Remark}
\newtheorem*{defn}{Definition}
\newtheorem*{eg}{Example}
\newtheorem{ex}{Exercise}[section]
\numberwithin{equation}{section}
\newcommand{\HJ}{Hertrich-Jeromin}
\newcommand{\Bt}{B\"acklund transformation}
\newcommand{\Ch}{Christoffel}
\newcommand{\MC}{Maurer--Cartan}
\newcommand{\HJMN}{\HJ--Musso--Nicolodi}
\newcommand{\HJP}{\HJ--Pedit}
\newcommand{\Bq}{Bianchi quadrilateral}

\begin{document}
\maketitle

\section*{Introduction}

\subsection*{Manifesto}

My aim is to give an account of the theory of isothermic surfaces in
$\R^n$ from the point of view of classical surface geometry and also
from the perspective of the modern theory of integrable systems and
loop groups.

There is some novelty even to the classical theory which arises from
the fact that isothermic surfaces are conformally invariant objects
in contrast, for example, to the more familiar surfaces of constant
Gauss or mean curvature.  Thus we have to do with a second order,
parabolic geometry which has its own flavour quite unlike Euclidean
or Riemannian geometry.  To compute effectively in this setting, we
shall develop an efficient calculus based on Clifford algebras.

The recent renaissance in interest in isothermic surfaces is
principally\footnote{Other motivations are available: see, for
example, the recent work of Kamberov-Pedit-Pinkall \cite{KamPedPin98}
on the Bonnet problem.} due to the fact that they constitute an
integrable system.  I shall attempt to explain in what sense this is
true and how this relates to the classical geometry.  In particular,
I shall show how the loop group formalism provides a context of
considerable generality in which results of Bianchi, Darboux and
others can be understood and generalised.

All of this will take some preparation so let us begin with an
overview of integrable geometry in general and isothermic surfaces in
particular.

\subsection*{Background}

\subsubsection*{What is an integrable system?}
This is a question with many answers of varying degrees of precision,
generality and plausibility!  For our present purposes, I take an
integrable system to be a geometric object or system of PDE with some
(or all) of the following features:
\begin{itemize}
\item an infinite-dimensional symmetry group;
\item the possibility of writing down explicit solutions;
\item a Hamiltonian formulation in which the system is completely
integrable in the sense of Liouville.
\end{itemize}
For Analysts, the prototype example of such a system is the
Korteweg--DeVries equation \cite{DraJoh89} or, perhaps, the non-linear
Schr\"odinger equation \cite{FadTak87} but, for Geometers, the basic
example, already well-known by the end of the 19th Century, is that
of pseudo-spherical surfaces: surfaces in $\R^3$ with constant Gauss
curvature $K=-1$.

Let us recall a little of this theory to fix ideas: let $f: M\to\R^3$
be an isometric immersion with $K=-1$.  According to B\"acklund, (see
\cite[\S120]{Eis60}), one can solve a first order Frobenius integrable
differential equation to obtain a second immersion $\fh:M\to\R^3$,
also with $K=-1$ determined by the geometric conditions that
\begin{enumerate}
\item $\fh-f$ is of constant length and tangent to both $\fh$ and
$f$;
\item normals at corresponding points of $f$ and $\fh$ make constant
angle with each other.
\end{enumerate}
This is the original \Bt\ of pseudo-spherical surfaces.
In this procedure there are two parameters: the angle $\sigma$
between the normals and an initial condition (also an angle) for the
differential equation.

Bianchi \cite[\S121]{Eis60} discovered a beautiful relation between iterated
\Bt s: the permutability theorem.  To describe this, we need a little
notation: for $f$ a pseudo-spherical surface, let $\mathcal{B}_\sigma
f$ denote a \Bt\ of $f$ with angle $\sigma$ between the normals.  Now
start with $f$ and let $f_1=\mathcal{B}_{\sigma_1}f$ and
$f_2=\mathcal{B}_{\sigma_2}f$ be two such \Bt s.  Then Bianchi's
theorem asserts the existence of a fourth pseudo-spherical surface
$\fh$ which is simultaneously a \Bt\ of $f_1$ and $f_2$:
\[
\fh=\mathcal{B}_{\sigma_1}f_2=\mathcal{B}_{\sigma_2}f_1.
\]
Moreover $\fh$ can be computed \emph{algebraically} from $f,f_1,f_2$.
In this way, we begin to see the first two of our desiderata for
integrability: the \Bt s generate an infinite-dimensional symmetry
group acting on the set of pseudo-spherical surfaces and the
permutability theorem shows the possibility of writing down explicit
solutions starting with a simple (possibly degenerate) $f$.

A modern view-point on these classical matters is provided by the
theory of loop groups.  The group generated by the \Bt s can be
identified as the group of rational maps of the Riemann sphere $\P^1$
into the complex orthogonal group $\SO[3,\C]$ satisfying the
conditions:
\begin{enumerate}
\item $g(0)=1$ and $g$ is holomorphic at $\infty$;
\item $g(\lambda)\in\SO[3]$ when $\lambda\in\R$;
\item for all $\lambda\in\dom(g)$, $g(-\lambda)=\tau g(\lambda)$
where $\tau$ is a certain involution of $\SO[3,\C]$.
\end{enumerate}
In this setting, the generators which act by \Bt s are distinguished
by having a pair of simple poles only\footnote{The position of the
poles prescribes the angle $\sigma$ while the residues there amount to the
initial condition of the differential equation.} while the
permutability theorem amounts to an assertion about products of these
generators.  This view-point is expounded in detail in the recent
work of Terng--Uhlenbeck \cite{TerUhl00} and we shall have much to
say about it below.

\subsubsection*{Where does integrability come from?}

A starting point from which all this rich structure can be derived is
a \emph{zero-curvature formulation} of the underlying problem.  That
is, the equations describing the problem should amount to the
flatness of a family of connections depending on an auxiliary
parameter.

Again, we illustrate the basic idea with the example of
pseudo-spherical surfaces: a pseudo-spherical surface $f$ admits
Chebyshev coordinates $\xi,\eta$, that is, asymptotic coordinates
for which the coordinate vector fields $\del/\del\xi,\del/\del\eta$
have unit length.  Now let $\theta:M\to\R$ be the angle between these
vector fields: the Gauss--Codazzi equations for $f$ amount to a
single equation, the \emph{sine-Gordon} equation for $\theta$:
\begin{equation}
\label{eq:1}
\theta_{\xi\eta}=\sin\theta,
\end{equation}
where, here and below, subscripts denote partial differentiation.

Now contemplate the pencil of connections
\[
\nabla^\lambda=\d +
\begin{pmatrix}
0&-\theta_\xi&\lambda\\\theta_\xi&0&0\\-\lambda&0&0
\end{pmatrix}\d\xi+
\begin{pmatrix}
0&0&\lambda^{-1}\cos\theta\\0&0&\lambda^{-1}\sin\theta\\
-\lambda^{-1}\cos\theta&-\lambda^{-1}\sin\theta&0
\end{pmatrix}\d\eta.
\]
By examining the coefficients of $\lambda$ in the curvature of
$\nabla^\lambda$, it is not difficult to show that $\nabla^\lambda$ is
flat for all $\lambda\in\R$ if and only if $\theta$ solves
\eqref{eq:1}.  Thus each pseudo-spherical surface gives rise to a
pencil of flat connections.

In fact, more is true: trivialising each $\nabla^\lambda$ produces, at
least locally, gauge transformations $F_\lambda:M\to\SO[3]$
intertwining $\nabla^\lambda$ and the trivial connection $\d$ and
from these one can construct a $1$-parameter family of
pseudo-spherical surfaces deforming $f$---these turn out to be the
\emph{Lie transforms} of $f$ \cite[\S122]{Eis60}.

These constructions are the starting point of a powerful and rather
general method for establishing integrability.  Indeed, a
zero-curvature formulation of a problem should yield:
\begin{itemize}
\item an action of a loop group on solutions;
\item a spectral deformation of solutions analogous to the Lie
transforms of pseudo-spherical surfaces;
\item explicit solutions via \Bt s or via algebraic geometry.
\end{itemize}

This theory has been fruitfully applied to a number of geometric
problems such as harmonic maps of surfaces into (pseudo-)Riemannian
symmetric spaces \cite{BurFerPed93,BurGue97,Hit90,Uhl89}; isometric
immersions of space forms in space forms \cite{FerPed96A,Ter97}; flat
Egoroff metrics \cite{TerUhl98}---these include (semisimple)
Frobenius manifolds \cite{Dub93,Hit97} and affine spheres
\cite{BobSch99}.

\subsubsection*{Isothermic surfaces in $\R^3$}

I will describe another classical differential geometric theory that
fits into this general picture: this is the theory of
\emph{isothermic surfaces}.  Classically, a surface in $\R^3$ is
isothermic if, away from umbilic points, it admits \emph{conformal
curvature line coordinates}, that is, conformal coordinates that,
additionally, diagonalise the second fundamental form \cite{Cay72}.
Here are some examples:
\begin{itemize}
\item surfaces of revolution;
\item quadrics;
\item minimal surfaces and, more generally, surfaces of constant mean
curvature.
\end{itemize}
There is a second characterisation of isothermic surfaces due to \Ch\
\cite{Chr67} : a surface $f:M\to\R^3$ is isothermic if and only if,
locally, there is a second surface, a \emph{dual surface},
$f^c:M\to\R^3$ with parallel tangent planes to those of $f$ which
induces the same conformal structure but opposite orientation on
$M$.  It is this view-point we shall emphasise below.

Isothermic surfaces were studied intensively at the turn of the 20th
century and a rich transformation theory of these surfaces was
developed that is strikingly reminiscent of that of pseudo-spherical
surfaces.  Darboux \cite{Dar99} discovered a transformation of
isothermic surfaces very like the \Bt\ of pseudo-spherical surfaces:
again the transform is effected by solving a Frobenius integrable
system of differential equations and again there is a geometric
construction only now the surface and its transform are enveloping
surfaces of a sphere congruence rather than focal surfaces of a line
congruence.  Moreover, Bianchi \cite{Bia05} proved the analogue of
his permutability theorem for these Darboux transformations.  Again,
Bianchi \cite{Bia05A}, and independently, Calapso \cite{Cal03} found
a spectral deformation of isothermic surfaces, the $T$-transform,
strictly analogous to the Lie transform of pseudo-spherical surfaces.

However, there is one important difference between the two theories:
that of pseudo-spherical surfaces is a Euclidean theory while that of
isothermic surfaces is a conformal one---the image of an isothermic
surface by a conformal diffeomorphism of $\R^3\cup\{\infty\}$ is also
isothermic.

Such an intricate transformation theory strongly suggests the
presence of an underlying integrable system.  That this is indeed the
case was established by Cie\'sli\'nski--Goldstein--Sym
\cite{CieGolSym95} who wrote down a zero-curvature formulation of the
Gauss--Codazzi equations of an isothermic surface.  This work was
taken up in \cite{BurHerPed97} where the conformal invariance of the
situation was emphasised and the underlying integrable system was
identified as an example of the \emph{curved flat} system of
Ferus--Pedit \cite{FerPed96}.  A new view-point on these matters was
provided by the Berlin school and their collaborators who developed a
beautiful quaternionic formalism for treating surfaces in
$4$-dimensional conformal geometry
\cite{BurFerLes00,Her97,HerMusNic,HerPed97,KamPedPin98,PedPin98}. In
particular, Hertrich-Jeromin--Pedit \cite{HerPed97} discovered a
description of Darboux transformations via solutions of a Riccati
equation which gives an extraordinarily efficient route into the
heart of the theory.

\subsection*{Overview}

My purpose in this paper is two-fold: firstly, I want to describe how
the entire theory of isothermic surfaces of $\R^3$ can be carried
through for isothermic surfaces in $\R^n$ \emph{with no loss of
integrable structure}.  Secondly, I shall show how this theory can be
profitably described using the loop group formalism and, in
particular, how to identify Darboux transformations with the dressing
action of \emph{simple factors} in the spirit of Terng--Uhlenbeck
\cite{TerUhl98,TerUhl00}.  In this way, I hope to exhibit the common
mechanism underlying the classical geometry of both pseudo-spherical
and isothermic surfaces.  Along the way, I shall describe a very
efficient method for doing conformal geometry which was inspired by
the quaternionic formalism of \HJP\ and, in fact, simultaneously
generalises and (at least when $n=4$) simplifies their approach.

Having declared our aims, let us turn to a more detailed description
of the topics we treat.  These can be grouped under three headings:
isothermic surfaces, loop groups and conformal geometry.

\subsubsection*{Isothermic surfaces in $\R^n$}

An isothermic surface in $\R^n$ can be defined just as in the
classical situation: either as a surface admitting conformal
curvature line coordinates (although we must now demand that the
surface have flat normal bundle in order for curvature lines to be
defined) or as a surface that admits a dual surface, that is, a
second surface with parallel tangent planes to the first, the same
conformal structure and opposite orientation.  That these two
characterisations locally coincide is due to Palmer \cite{Pal88}.

The starting point of our study is the observation that two
immersions $f,f^c:M\to\R^n$ are dual isothermic surfaces if and only
if
\begin{equation}
\label{eq:2}
\d f\wedge \d f^c =0,
\end{equation}
where we multiply the coefficients of these $\R^n$-valued $1$-forms
using the product of the Clifford algebra $\Cl_n$ of $\R^n$.
Equation \eqref{eq:2} is the integrability condition for a Riccati
equation involving an auxiliary parameter $r\in\R^\times$: 
\[
\d g =r g\d f^c g-\d f
\]
where again all multiplications take place in $\Cl_n$.  We now
construct a new isothermic surface $\fh$ by setting $\fh=f+g$: this
is the Darboux transform of $f$.  We show that, just as in the
classical case, $f$ and $\fh$ are characterised by the conditions
that they have the same conformal structure and curvature lines and
are the enveloping surfaces of a $2$-sphere congruence.  This is
perhaps, a surprising result: a generic congruence of $2$-spheres in
$\R^n$ has no enveloping surfaces at all!

This approach to the Darboux transform is a direct extension of that
of \HJP\ for the case $n=3,4$ and we follow their methods to prove the
Bianchi permutability theorem for Darboux transforms.  This proceeds
by establishing an explicit algebraic formula for the fourth
isothermic surface which comes from the ansatz that corresponding
points on the four surfaces in the Bianchi configuration should have
constant (Clifford algebra) cross-ratio.  We shall find some a priori
justification for this ansatz.  Further analysis of the Clifford
algebra cross-ratio allows us to extend other results of Bianchi in
this area to $n$ dimensions: in particular, we prove that the Darboux
transform of a Bianchi quadrilateral is another such.

Again, isothermic surfaces in $\R^n$ are conformally invariant and
admit a spectral deformation, the $T$-transform.  We explain the
intricate relationships between the T-transforms and Darboux
transforms of an isothermic surface and its dual.

Examples of isothermic surfaces in $\R^3$ are provided by surfaces of
constant mean curvature.  In fact, non-minimal CMC surfaces can be
characterised as those isothermic surfaces whose dual surface is also
a Darboux transform \cite{Bia05} and such surfaces are preserved by a
co-dimension $1$ family of Darboux transforms.  In $\R^n$, we find an
exactly analogous theory for \emph{generalised $H$-surfaces}, that
is, surfaces which admit a parallel isoperimetric section in the
sense of Chen \cite{Che73A}.  In fact, these methods have a wider
applicability: applying the same formal arguments in a different
algebraic setting establishes the existence of a family of B\"acklund
transformations of Willmore surfaces in $S^4$ \cite{BurFerLes00}.

We complete our extension of the classical theory to $n$ dimensions
by considering the approach of Calapso \cite{Cal03} who showed that
an isothermic surface together with its $T$-transforms amounts to a
solution $\kappa:M\to\R$ of the Calapso equation
\begin{equation}
\label{eq:3}
\Delta\biggl(\frac{\kappa_{xy}}{\kappa}\biggr)+2(\kappa^2)_{xy}=0.
\end{equation}
A straightforward generalisation of the analysis in
\cite{BurHerPed97} shows that the same is true in $\R^n$ if
\eqref{eq:3} is replaced by a \emph{vector Calapso equation}:
\begin{gather*}
\kappa_{xy}=\psi\kappa\\
\Delta\psi+2\biggl(\sum_{i=1}^{n-2}\kappa_i\biggr)_{xy}=0
\end{gather*}
for $\kappa:M\to\R^{n-2}$ and $\psi:M\to\R$.  Moreover, we identify
$\kappa$ as (the components of) the conformal Hopf differential of
the isothermic surface.

Several of these results have been proved independently by Schief
\cite{Sch} who, in particular, established the existence and
Bianchi permutability of Darboux transforms in this context as well
as the description of isothermic surfaces via the vector Calapso
equation.

Curved flats are submanifolds of a symmetric space on whose tangent
spaces the curvature operator vanishes \cite{FerPed96}.  Curved flats
admit a zero-curvature representation and so the methods of
integrable systems theory apply.  A main result of \cite{BurHerPed97}
is that an isothermic surface in $\R^3$ together with a Darboux
transform $\fh$ constitute a curved flat $(f,\fh):M\to S^3\times
S^3\setminus\Delta$ in the space of pairs of distinct points in
$S^3=\R^3\cup\{\infty\}$.  This last is a pseudo-Riemannian symmetric
space for the diagonal action of the M\"obius group of conformal
diffeomorphisms of $S^3$.  This result goes through unchanged in the
$n$-dimensional setting where we identify the spectral deformation of
curved flats with the $T$-transforms of the factors.  In fact, more
is true: a curved flat and its spectral deformations give rise via a
limiting procedure (Sym's formula) to a certain map of $M$ into a
tangent space $\p$ of the symmetric space.  We call these
\emph{$\p$-flat maps} and show that the converse holds: a $\p$-flat
map gives rise to a family of curved flats.  In the case of
isothermic surfaces, a $\p$-flat map is the same as an isothermic
surface together with a dual surface and our procedure produces the
family of the $T$-transforms of this dual pair.  By passing to frames
of this family, one obtains an extended object which can be viewed as
a map from $M$ into an infinite dimensional group of holomorphic maps
$\C\to\O[n+2,\C]$.  This is the key to the application of loop group
methods to isothermic surfaces to which we now turn.

\subsubsection*{Loop groups}

There is a very general mechanism, pervasive in the theory of
integrable systems, for constructing a group action on a space of
solutions.  Here is the basic idea: let $\G$ be a group with
subgroups $\G_1,\G_2$ such that $\G_1\G_2=\G$ and
$\G_1\cap\G_2=\{1\}$.  Then $\G_2\cong \G/\G_1$ so that we get an
action of $\G$ and, in particular, $\G_1$ on $\G_2$.  In concrete
terms, for $g_i\in\G_i$, the product $g_1g_2$ can be written in a
unique way
\begin{equation}
\label{eq:4}
g_1g_2=\gh_2\gh_1
\end{equation}
with $\gh_i\in\G_i$ and then the action is given by
\[
g_1\# g_2=\gh_2.
\]
More generally, when $\G_1\G_2$ is only open in $\G$, one gets a local
action.

The case of importance to us is when the $\G_i$ are groups of
holomorphic maps from subsets of the Riemann sphere $\P^1$ to a
complex Lie group $G^\C$ distinguished by the location of their
singularities.  For example, in our applications to isothermic
surfaces, we take $\G_2$ to be a group of holomorphic maps
$\C\to\O[n+2,\C]$ and $\G_1$ a group of rational maps from $\P^1$ to
$\O[n+2,\C]$ which are holomorphic near $0$ and $\infty$.  The whole
point is that, as we have indicated above, isothermic surfaces give
rise to certain maps, \emph{extended flat frames}, $M\to\G_2$ of a
type that is preserved by the point-wise action of $\G_1$.  In this
way, we find a local action of $\G_1$ on the set of (dual pairs of)
isothermic surfaces and, more generally, on the set of $\p$-flat maps.

This is a phenomenon that is not peculiar to isothermic surfaces: the
key ingredient is that the extended frame is characterised completely
by the singularities of its derivative and this ingredient is shared
by many integrable systems with a zero-curvature representation (see
\cite{TerUhl00} for many examples).

It remains to compute this action which amounts to performing the
factorisation \eqref{eq:4}.  In general, this is a Riemann--Hilbert
problem for which explicit solutions are not available.  However, it
is philosophy developed by Terng and Uhlenbeck
\cite{TerUhl98,TerUhl00,Uhl89,Uhl92} that there should be certain
basic elements of $\G_1$, the \emph{simple factors}, for which the
factorisation \eqref{eq:4} can be computed explicitly and, moreover,
the action of these simple factors should amount to B\"acklund-type
transformations of the underlying geometric problem.  A difficulty
with this approach is that simple factors for a given situation are
constructed on an ad hoc basis.  We shall propose a concrete
characterisation of simple factors which has the status of a theorem
when the underlying geometry is that of a compact Riemannian
symmetric space and that of an ansatz in non-compact situations (the
case of relevance to isothermic surfaces).  

For isothermic surfaces, we show that the action of the simple
factors we find in this way amount to the Darboux transformations.
As a consequence, we find simple complex-analytic arguments that
provide a second proof of the circle of results around Bianchi
permutability which apply in a variety of contexts.

\subsubsection*{Conformal geometry and Clifford algebras}

The underlying setting for our theory of isothermic surfaces is that
of conformal geometry: the basic objects of study are conformally
invariant as are many of our ingredients and constructions: sphere
congruences, Darboux and $T$-transforms.  This explains the
appearance of the indefinite orthogonal group $\Op$ and its
complexification $\O[n+2,\C]$ as this is precisely the group of
conformal diffeomorphisms of $S^n=\R^n\cup\{\infty\}$.  Indeed, $S^n$
can be identified with the projective light-cone of $\lor$ and then
the projective action of $\O$ is by conformal diffeomorphisms giving
an isomorphism of an open subgroup $\Op$ with the M\"obius group.

The presence of Clifford algebras in this context, while not new, is
not so well known and deserves further comment.  Everyone knows how
the conformal diffeomorphisms of the Riemann sphere are realised on
$\C$ by the action of $\SL[2,\C]$ through linear fractional
transformations.  There is a completely analogous theory in higher
dimensions due to Vahlen \cite{Vah02} that replaces $\SL[2,\C]$ with
a group of $2\times 2$ matrices with entries in a Clifford algebra.
Here is the basic idea: instead of working with $\Op$, we pass to a
double cover and work with an open subgroup of $\Pin$ which is itself
a multiplicative subgroup of the Clifford algebra $\Cl_{n+1,1}$ of
$\lor$.  The point now is that $\Cl_{n+1,1}$ is isomorphic to the
algebra of $2\times2$ matrices with entries in the Clifford algebra
$\Cl_n$ of $\R^n$.  Moreover, a theorem of Vahlen identifies the
matrices that comprise the double cover of $\Op$ and once again these
act on $\R^n$ by linear fractional transformations.

This beautiful formalism is well suited to the study of isothermic
surfaces: it makes the action of the M\"obius group on $\R^n$
particularly easy to understand and leads to extremely compact
formulae in moving frame calculations and elsewhere.  While Vahlen's
ideas have been used in hyperbolic geometry (see, for example,
\cite{ElsGruMen87,ElsGruMen90,Wad90}) and harmonic analysis
\cite{GilMur91}, I believe that this is the first time\footnote{See,
however, Cie\'sli\'nski's direct use of $\Cl_{4,1}$ in his study of
isothermic surfaces in $\R^3$ \cite{Cie97A}.}  that these methods have
been used in a thorough-going way to do conformal differential
geometry.

\subsection*{Road Map}
 
To orient the Reader, we briefly outline the contents of each section
of the paper.

Section 1 is preparatory in nature: we describe the light-cone model
of the conformal $n$-sphere and introduce submanifold geometry in
this context.  We set up the approach via Clifford algebras and use
it to prove some preliminary results.

Section 2 contains our account of the classical geometry of
isothermic surfaces in $\R^n$.  We define \Ch, Darboux and
$T$-transformations of these surfaces and investigate the
permutability relations between them.  We consider the special case
of generalised $H$-surfaces and digress to contemplate the vector
Calapso equation.

Section 3 is devoted to curved flats.  We describe the relation
between curved flats and $\p$-flat maps and how it specialises to
give the relation between a dual pair of isothermic surfaces and
their $T$-transforms.  We shall see that much of our preceding theory
of isothermic surfaces is unified by this curved flats
interpretation.

Section 4 deals with loop groups and how they may be applied to study
curved flats in general and isothermic surfaces in particular.  We
give a general discussion of simple factors and then specialise to
give a detailed account of the case of isothermic surfaces.

Section 5 rounds things off with brief descriptions of recent
developments and some open problems.

\subsection*{A note on the text}

This work had its genesis in lecture notes for a short course on
``Integrable systems in conformal geometry'' given at Tsing Hua
University in January 1999 but has evolved into a statement of
Everything I Know About Isothermic Surfaces.  However, this final
version has retained something of its origins in that I have given a
somewhat leisurely account of background material and also in that I
have set a large number of exercises.  These exercises are an
integral part of the exposition and, among other things, contain most
of the computations where no New Idea is needed.  Solutions may
become available at
\begin{quotation}
\url{http://www.maths.bath.ac.uk/~feb/taiwan-solutions.html}
\end{quotation}
and Readers are warmly invited to contribute their own!

\subsection*{Acknowledgements}

This paper is the product of a lengthy period of research during
which I have incurred many debts of gratitude.  I have enjoyed the
hospitality of the Dipartimento di Matematica\footnote{Visits to Rome
were supported by MUNCH and the Short-Term Mobility Program of the
CNR.} ``G. Castelnuovo'' at the Universit\`a di Roma ``La Sapienza'',
SFB288 \emph{Differential Geometry and Quantum Physics} at the
Technische Universit\"at Berlin and the National Centre for
Theoretical Sciences, Tsing Hua University, Taiwan.  Moreover I have
benefited greatly from conversations with A.~Bobenko, M.~Br\"uck,
D.~Calderbank, B.-Y.~Chen, J.~Cie\'sli\'nski, J.-H.~Eschenburg, C.
McCune, E.~Musso, L.~Nicolodi, F.~Pedit, U.~Pinkall, the participants
of the $H$-seminar at SFB288 in the autumn of 1998, C.L.~Terng and my
audience at Tsing Hua University in January 1999.

Special thanks are due to C.-L. Terng for the invitation to lecture
in Taiwan, W.~Schief for informing me of his work and, above all, to
Udo~Hertrich-Jeromin who has had a decisive influence on my thinking
about isothermic surfaces and conformal geometry.

\subsubsection*{Note added in October 2000}
\label{sec:note-added-proof}

Some time after this paper was written, I have had the opportunity to
read a wonderful book by Tzitz\'eica \cite{Tzi24} which contains a
completely different approach to isothermic surfaces: Tzitz\'eica
studies surfaces in an $n$-dimensional projective quadric that
support a conjugate net with equal Laplace invariants.  He observes
that, when $n=3$, these are exactly the isothermic surfaces (the
conjugate net is that formed by curvature lines while the quadric is
the projective light-cone of our exposition) and develops a theory of
Darboux transformations of such surfaces for arbitrary $n$.  It is
not hard to see that, for any $n$, Tzitz\'eica's surfaces amount
(locally) to precisely the isothermic surfaces in the conformal
compactification of some $\R^{p,q}$ with $p+q=n$.  Thus, in this way,
isothermic surfaces and their Darboux transformations have been known
in this generality since 1924!

\section{Conformal geometry and Clifford algebras}
\label{sec:preliminaries}

\subsection{Conformal geometry of $S^n$}
\label{sec:conf-geom-sn}

Recall that a map $\phi:(M,g)\to(M,g)$ of a Riemannian manifold is
\emph{conformal} if $\d \phi$ preserves angles.  Analytically this
means
\[
\phi^*g=e^{2u} g
\]
for some $u:M\to\R$.

Here are some conformal maps of (open sets of) $\R^n$:
\begin{enumerate}
\item Euclidean motions: $\phi^*g=g$;
\item Dilations: $x\mapsto rx$, $r\in\R^+$;
\item Inversions in hyperspheres: for fixed $p\in\R^n$, $r\in\R^+$,
these are $\phi:\R^n\setminus\{p\}\to\R^n$ given by
\[
\phi(x)=p+r^2\frac{x-p}{\norm{x-p}^2}.
\]
\begin{figure}[htbp]
\begin{center}
\includegraphics{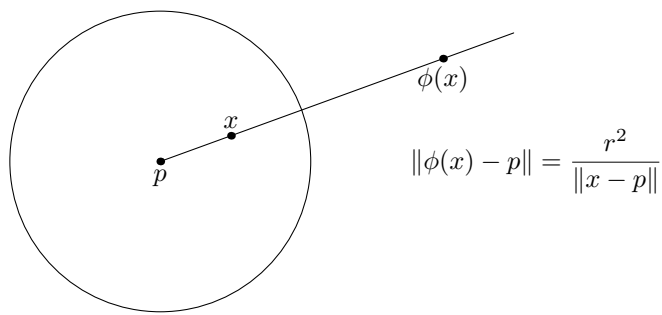}
\caption{Inversion in the sphere of radius $r$ about $p$}
\label{fig:1}
\end{center}
\end{figure}
\end{enumerate}

\begin{ex}
Show that such inversions are conformal.
\end{ex}

A theorem of Liouville states:
\begin{thm}
\label{th:1}
For $n\geq 3$, any conformal map $\Omega\subset\R^n\to\R^n$ is the
restriction to $\Omega$ of a composition of Euclidean motions,
dilations and inversions.
\end{thm}
For a proof, see do~Carmo \cite{doC76}.

\begin{figure}[htbp]
\begin{center}
\includegraphics{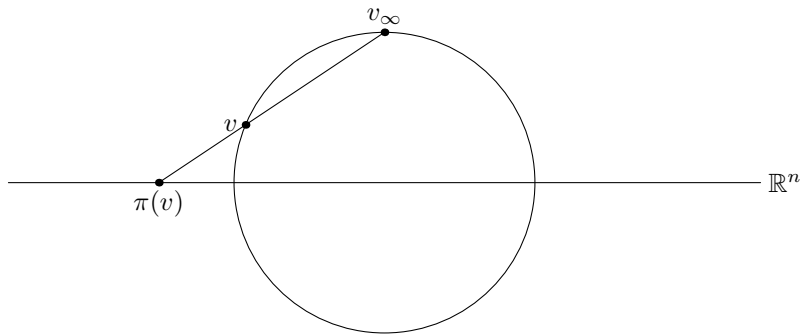}
\end{center}
\caption{Stereographic projection}
\label{fig:2}
\end{figure}

It is natural to extend the definition of the inversions
$\phi:\R^n\setminus\{p\}\to\R^n$ by setting $\phi(p)=\infty$ and
$\phi(\infty)=p$ and so viewing $\phi$ as a conformal diffeomorphism
of the $n$-sphere $\R^n\cup\{\infty\}=S^n$.  To make sense of this,
recall that the conformal geometries of $\R^n$ and
$S^n\subset\R^{n+1}$ are linked by stereographic projection: choosing
a ``point at infinity'' $v_\infty\in S^n$, we have a conformal
diffeomorphism
$\pi:S^n\setminus\{v_\infty\}\to\<v_\infty\>^\perp\cong\R^n$ as in
Figure~\ref{fig:2}.
\begin{ex}
Prove:
\begin{enumerate}
\item $\pi$ is a conformal diffeomorphism;
\item $S\subset\R^n$ is a $k$-sphere if and only if
$\pi^{-1}(S)\subset S^n$ is a $k$-sphere;
\item $V\subset\R^n$ is an affine $k$-plane if and only if
$\pi^{-1}(V)\cup\{v_\infty\}\subset S^n$ is a $k$-sphere containing
$v_\infty$.  Thus ``planes are spheres through infinity''.
\end{enumerate}
\end{ex}

Under stereographic projection, inversions in hyperspheres extend to
conformal diffeomorphisms of $S^n$ as do Euclidean motions and
dilations (these fix $v_\infty$) and, in this way, we are led to
consider the \emph{M\"obius group} $\Mob$ of conformal
diffeomorphisms of $S^n$.

To go further, it is very convenient to introduce another model of
the $n$-sphere discovered by Darboux\footnote{For a modern account,
see Bryant \cite{Bry84}.} \cite{Dar87}.  For this, we contemplate the
Lorentzian space $\R^{n+1,1}$: a real $(n+2)$-dimensional vector space
equipped with an inner product $(\,,\,)$ of signature $(n+1,1)$ so
that there is an orthonormal basis $e_1,\dots,e_{n+2}$ with
\[
(e_i,e_i)=
\begin{cases}
1&i<n+2\\-1&i=n+2.
\end{cases}
\]

Inside $\R^{n+1,1}$, we distinguish the \emph{light-cone} $\L$:
\[
\L=\{v\in\R^{n+1,1}\setminus\{0\}:(v,v)=0\}.
\]
\begin{ex}
$\L$ is a submanifold of $\R^{n+1,1}$.
\end{ex}
Clearly, if $v\in\L$ and $r\in\R^\times$ then $rv\in\L$ so that
$\R^\times$ acts freely on $\L$ and we may take the quotient
$\PL\subset\mathbb{P}(\R^{n+1,1})$:
\[
\PL=\L/\R^\times=\{\ell\subset\lor:\text{$\ell$ is a $1$-dimensional
isotropic subspace}\}.
\]
The point of this is that $\PL$ has a conformal structure with
respect to which it is conformally diffeomorphic to $S^n$ with its
round metric.  Indeed, let us fix a unit time-like vector
$t_0\in\lor$ (thus $(t_0,t_0)=-1$) and set
\[
S_{t_0}=\{v\in\L:(v,t_0)=-1\}.
\]
For $v\in S_{t_0}$, write $v=v^\perp+t_0$ so that $v^\perp \perp t_0$
and note:
\[
0=(v,v)=(v^\perp,v^\perp)+(t_0,t_0)=(v^\perp,v^\perp)-1.
\]
Thus the projection $v\mapsto v^\perp$ is a diffeomorphism
$S_{t_0}\to S^n$ onto the unit sphere in $\<t_0\>^\perp\cong\R^{n+1}$
which is easily checked to be an isometry.
\begin{ex}
\begin{enumerate}
\item For $v\in\L$, $(t_0,v)\neq 0$.
\item Deduce that each line $\ell\in\PL$ intersects $S_{t_0}$ in
exactly one point.
\end{enumerate}
\end{ex}
Thus we have a diffeomorphism $\ell\mapsto S_{t_0}\cap\ell:\PL\to
S_{t_0}$ whose inverse is the canonical projection
$\pi:v\mapsto\<v\>:\L\to\PL$ restricted to $S_{t_0}$.
\begin{ex}
Suppose that $t'_0$ is another unit time-like vector.  Show that the
composition $S_{t'_0}\stackrel{\pi}{\to}\PL\cong S_{t_0}$ is a
conformal diffeomorphism which is \emph{not} an isometry unless
$t'_0=\pm t_0$.
\end{ex}

To summarise the situation: \emph{for each unit time-like $t_0$, $\pi$
restricts to a diffeomorphism $S_{t_0}\to\PL$ and each such
diffeomorphism induces a conformally equivalent metric on $\PL$.}

\begin{ex}
Let $g$ be any Riemannian metric on $S^n$ in the conformal class of
the round metric.  Show that there is an isometric embedding
$(S^n,g)\to\L$.
\end{ex}

Having identified $\PL$ with the conformal $n$-sphere, we can use a
similar argument to describe stereographic projection in this model
by replacing time-like $t_0$ with $v_\infty\in\L$: fix
$v_0,v_\infty\in\L$ with $(v_0,v_\infty)=-\half$ (so that, in
particular, $\<v_0\>\neq\<v_\infty\>$) and set
\[
E_{v_\infty}=\{v\in\L:(v,v_\infty)=-\half\}.
\]
\begin{ex}
For $v\in\L\setminus\<v_\infty\>$, show that $(v,v_\infty)\neq0$.
\end{ex}
Thus we have a diffeomorphism $\ell\mapsto
E_{v_\infty}\cap\ell:\PL\setminus\{\<v_\infty\>\}\to E_{v_\infty}$.
Moreover, $E_{v_\infty}$ is isometric to a Euclidean space: indeed,
set $\R^n=\<v_0,v_\infty\>^\perp$, a subspace of $\lor$ on which the
inner product is definite.
\begin{ex}
\begin{enumerate}
\item There is an isometry $E_{v_\infty}\to\R^n$ given by
\begin{equation}
\label{eq:5}
v\mapsto v-v_0+2(v,v_0)v_\infty
\end{equation}
with inverse
\begin{equation}
\label{eq:6}
x\mapsto x+v_0+(x,x)v_\infty.
\end{equation}
\item Verify that the composition $\PL\setminus\<v_\infty\>\cong
E_{v_\infty}\to\R^n$ really is stereographic projection.

More precisely, set $t_0=v_0+v_\infty$, $x_0=v_0-v_\infty$ so that
$(t_0,t_0)=-1=-(x_0,x_0)$.  Let $S^n$ be the unit sphere in
$\<t_0\>^\perp$.  Then the composition
\[
S^n\setminus\{x_0\}\to
S_{t_0}\setminus\{2v_\infty\}\stackrel{\pi}{\to}
\PL\setminus\<v_\infty\>\to\R^n=\<v_0,v_\infty\>^\perp=\<t_0,x_0\>^\perp
\]
is stereographic projection.
\end{enumerate}
\end{ex}

The beauty of this model is that it \emph{linearises} conformal
geometry.  For example, observe that the set of hyperspheres in $S^n$
is parametrised by the set $\P^+(\lor)$ of space-like lines, that is,
$1$-dimensional subspaces on which the inner product is positive
definite.  Indeed, if $L\subset\lor$ is such a line then
$L^\perp\cong\R^{n,1}$ so that $\P(\L\cap L^\perp)\cong S^{n-1}$.
Choosing $v_0,v_\infty$ and so a choice of stereographic projection,
this correspondence becomes quite explicit:
\begin{ex}\label{ex:1}
Let $L\in\P^+(\lor)$ and fix $s\in L$ of unit length.  Let $s^\perp$
be the orthoprojection of $s$ onto $\R^n=\<v_0,v_\infty\>^\perp$ and
set $S_L=\P(\L\cap L^\perp)$.
\begin{enumerate}
\item If $\<v_\infty\>\in S_L$, that is, $(v_\infty,s)=0$, then the
stereo-projection of $S_L\setminus\{\<v_\infty\>\}$ is the hyperplane
\[
\{x\in\R^n:(x,s^\perp)=-(s,v_0)\}.
\]
\item If $\<v_\infty\>\not\in S_L$, then the stereo-projection of
$S_L$ is the sphere centred at $-s^\perp/2(s,v_\infty)$ of radius
$1/2\abs{(s,v_\infty)}$.
\item Stereo-projection intertwines reflection in the hyperplane
$L^\perp\subset\lor$ with reflection or inversion in the plane or
sphere determined by $S_L$.
\end{enumerate}
\end{ex}

Now contemplate the orthogonal group $\O$ of $\lor$, that is,
\[
\O=\{T\in\mathrm{GL}(n+2,\R):\text{$(Tu,Tv)=(u,v)$ for all
$u,v\in\lor$}\}.
\]
The linear action of $\O$ on $\lor$ preserves $\L$ and the set of
lines in $\L$ and so descends to an action on $\PL$.  Moreover, for
$t_0$ a unit time-like vector and $T\in\O$, $T$ restricts to give an
isometry $S_{t_0}\to S_{Tt_0}$ so that the induced map on $\PL$ is a
\emph{conformal} diffeomorphism.  In this way, we have found a
homomorphism $\O\to\Mob$ which is, in fact, a double cover:
\begin{thm}
The sequence $0\to\Z_2\to\O\to\Mob\to0$ is exact.
\end{thm}
\begin{proof}
Any $T$ in the kernel of our homomorphism must preserve each
light-line and so has each light-line as an eigenspace.  This forces
$T$ to be a multiple of the identity matrix $I$ and then $T\in\O$
gives $T=\pm I$.

Thus the main issue is to see that our homomorphism is onto.
However, by Liouville's Theorem, $\Mob$ is generated by reflections
in hyperplanes and inversions in hyperspheres: indeed, any Euclidean
motion is a composition of reflections while a dilation is a
composition of two inversions in concentric spheres.  On the other
hand, we have seen in Exercise~\ref{ex:1} that all these reflections
and inversions are induced by reflections in $L^\perp$ for
$L\in\P^+(\lor)$ a space-like line.  Such reflections are certainly
in $\O$ and we are done.
\end{proof}

In fact, we can do better: the light cone $\L$ has two
components\footnote{Non-collinear elements $v_0,v_\infty\in\L$ are in
the same component if and only if $(v_0,v_\infty)<0$.} $\L^+$ and
$\L^-$ which are transposed by $v\mapsto-v$.  Correspondingly, $\O$
has four components distinguished by the sign of the determinant and
whether or not the components of $\L$ are preserved.  Denote by $\Op$
the subgroup of $\O$ that preserves $\L^\pm$.  Then $-I\not\in\Op$
and we deduce:
\begin{thm}\label{th:2}
$\Op\cong\Mob$.
\end{thm}
The two components of $\Op$ are the orientation preserving and
orientation reversing conformal diffeomorphisms of $\PL$.

\subsection{Submanifold geometry in $\PL$}
\label{sec:subm-geom-pl}

\subsubsection{Submanifolds and normal bundles}
\label{sec:subm-norm-bundl}

Contemplate the projection $\pi:\L\to\PL$, $\pi(v)=\<v\>$.  Clearly,
$T_v\L=\<v\>^\perp$ while $\ker\d \pi_v=\<v\>$ so we have an
isomorphism
\[
\d\pi_v:\sv^\perp/\sv\cong T_{\sv}\PL.
\]
Scaling $v$ leaves $\sv^\perp/\sv$ unchanged but scales the
isomorphism:
\begin{ex}\label{ex:2}
For $r\in\R^\times$ and $X\in\sv^\perp/\sv$,
$\d\pi_{rv}(rX)=\d\pi_v(X)$.
\end{ex}
More invariantly, we have an isomorphism
$\Hom(\sv,\sv^\perp/\sv)\cong T_{\sv}\PL$ given by
\[
B\mapsto\d\pi_v(Bv)
\]
which is well-defined by Exercise~\ref{ex:2}.

For $\ell\in\PL$, the inner product on $\lor$ induces a positive
definite inner product on $\ell^\perp/\ell$ and if $v\in\ell^\times$
lies in some round sphere $S_{t_0}$ then projection along $\ell$ is
an isometry $T_vS_{t_0}\to\ell^\perp/\ell$.  We therefore conclude
that the isomorphism $\d\pi_v: \ell^\perp/\ell\to T_\ell\PL$ is
conformal.

Now let $M$ be a manifold and $\phi:M\to\PL$ an immersion.  We study
$\phi$ by studying its \emph{lifts}, that is, maps $f:M\to\L$ with
$\pi\circ f=\phi$.  Since the principal $\R^\times$-bundle $\pi:\L\to\PL$
is trivial (each $S_{t_0}$ is the image of a section!) there are many
such lifts.  Moreover, in view of Theorem~\ref{th:2}, it suffices to
consider lifts $f:M\to\L^+$.  If $f$ is one such, then any other is
of the form $e^uf$ for some $u:M\to\R$.

So let $f:M\to\L^+$ be a lift of $\phi=\sf$.  Then $\d\pi_f$ gives an
isomorphism
\[
\sf^{-1}T\PL\cong\sf^\perp/\sf
\]
under which the derivative of $\sf$ is given by $\d f\mod\sf$.  In
particular, $\sf$ is an immersion if and only if, for each $X\in TM$,
\[
\d_X f\wedge f\neq 0.
\]
Scaling the lift scales the isomorphism:
\[
\d(e^u f)=e^u(\d uf+\d f)\equiv e^u\d f\mod\sf
\]
from which we see that the image of $\d f$ in $\sf^\perp/\sf$ is
independent of the choice of lift.  Thus, when $\sf$ is an immersion,
orthogonal decomposition gives a well-defined \emph{weightless normal
bundle}\footnote{Strictly speaking, the normal bundle to $\sf$ is
$\Hom(\sf,\Nf)=\sf^*\otimes\Nf\subset\sf^{-1}T\PL$.  Since we will
mostly deal with lifts $f$ of $\sf$, we shall ignore this distinction
which, in any case, amounts only to tensoring with a trivial line
bundle.}  $\Nf$:
\begin{equation}
\label{eq:7}
\sf^\perp/\sf=\im\d f(TM)\oplus\Nf
\end{equation}
which is M\"obius invariant: for $T\in\O$,
\[
\mathcal{N}_{T\sf}=T\Nf.
\]
\begin{notation}
For $s$ a section of $\sf^\perp$, write
\[
s+\sf= [s]^T+[s]^\perp
\]
according to the decomposition \eqref{eq:7}.
\end{notation}

\subsubsection{Conformal invariants}
\label{sec:conformal-invariants}

We construct conformal invariants of submanifolds by finding
$\Op$-invariant properties of lifts that do not depend on the choice
of said lift.

Firstly, we have the conformal class of the metric induced by $f$ on
$M$:
\[
\bigl(\d(e^uf),\d(e^uf)\bigr)=e^{2u}\bigl(\d f,\d f\bigr)
\]
since $(\d f,f)=0$ ($\d f$ is $\sf^\perp$-valued).

Now let $N$ be a section of $\sf^\perp$ such that $[N]=N+\sf$ is a
section $\Nf$.  Thus
\[
(N,f)=(N,\d f)=0
\]
whence
\[
(\d N,f)=-(N,\d f)=0
\]
so that $\d N$ is $\sf^\perp$-valued also.  Moreover, if $[N]=[N']$
so that $N'=N+\mu f$, for some function $\mu:M\to\R$, we have
\[
\d N'=\d N+\mu\d f+\d \mu f
\]
whence
\[
\d N'\equiv \d N+\mu\d f\mod \sf.
\]
In particular,
\[
[\d N]^\perp=[\d N']^\perp
\]
so that we can define a conformally invariant connection
$\nabla^\perp$ on $\Nf$ by
\[
\nabla^\perp[N]=[\d N]^\perp.
\]
Further,
\[
\d f^{-1}([\d N']^T)=\d f^{-1}([\d N]^T)+\mu
\]
so that shape operators are well-defined up to addition of multiples
of the identity and scaling (as the lift varies).  In particular, the
eigenspaces of shape operators, \emph{the principal curvature
directions}, are well-defined.

To summarise: given an immersion $\sf:M\to\PL$, we obtain in a
M\"obius invariant way:
\begin{enumerate}
\item A conformal class of metrics on $M$;
\item A weightless normal bundle $\Nf$ with normal connection $\nabla^\perp$;
\item The conformal class of \emph{trace-free} shape operators.
\end{enumerate}
\begin{rem}
A more precise formulation of these invariants can be obtained by
viewing $\sf$ as the sub-bundle of the trivial bundle $M\times\lor$
whose fibre at $p\in M$ is $\<f(p)\>\subset\lor$.  Following Calderbank
\cite{Cal98}, our conformal class of metrics on $M$ can be viewed as
an honest metric on $TM\otimes\sf$ via
\[
(X\otimes f,Y\otimes f)=(\d_Xf,\d_Yf),
\]
$X,Y\in TM$.
In the same way, the conformal class of trace-free shape operators
can be viewed as a single trace-free quadratic form taking values in
$\Nf\otimes\sf$.  We shall return to this viewpoint on conformal
submanifold geometry elsewhere.
\end{rem}

When $M$ is an orientable surface, our analysis can be refined
somewhat.  In this case, $M$ becomes a Riemann surface so let
$z=x+iy$ be a holomorphic coordinate and take a lift $f:M\to\L^+$.
We define a local section $K_{\sf}$ of $\Nf^*$ by
\[
K_{\sf}(N+\sf)=\sqrt{2}\frac{(f_{zz},N)}{\sqrt{(f_z,f_{\zbar})}}
\]
where, here and below, we use subscripts to denote partial
differentiation.
\begin{ex}
$K_{\sf}$ is well-defined and independent of the choice of lift $f$
in $\L^+$.
\end{ex}
It is clear that $K_{\sf}$ is equivariant under the action of $\Op$:
for $T\in\Op$,
\[
K_{T\sf}\circ T=K_{\sf}
\]
and so is conformally invariant.  As for the dependence on the
holomorphic coordinate $z$, we see that $K_{\sf}$ should be viewed a
density with values in $\Nf^*$, that is, as a section of $(\bigwedge
^{1,0}M)^{3/2}\otimes (\bigwedge^{0,1}M)^{-1/2}\otimes\Nf^*$.

We call $K_{\sf}$ the \emph{conformal Hopf differential} of $\sf$ and
will return to this topic in
Section~\ref{sec:isoth-surf-via}.\label{page:conf-Hopf}

\subsubsection{Spheres and sphere congruences}
\label{sec:spher-sphere-congr}

We have already seen that the hyperspheres in $S^n$ are parametrised
by the space $\P^+(\lor)$ of space-like lines.  In the same way, the
Grassmannian $G^+_k(\lor)$ of space-like $k$-planes in $\lor$
parametrises co-dimension $k$ spheres in $S^n$ \cite{Roz48}.  Indeed, any such
sphere is of the form $S_\Pi=\P(\Pi^\perp\cap\L)$ for a unique
$\Pi\in G^+_k(\lor)$.

Now let $\sf:M\to\PL$ be an immersion and $\Pi\in G^+_k(\lor)$.  For
$p\in M$, we see that $\sf(p)\in S_\Pi$ if and only if $\Pi\perp
f(p)$ while $\sf$ is tangent to $S_\Pi$ at $p$ if, in addition,
$\Pi\perp\im\d f_p$.

\begin{defn}
A \emph{congruence of $k$-spheres} is a map $\Pi:M\to
G^+_{n-k}(\lor)$ of a $k$-dimensional manifold into the space of
$k$-spheres.
\end{defn}

An immersion $\sf:M\to\PL$ \emph{envelopes} the congruence $\Pi$ if,
for each $p\in M$, the sphere $S_{\Pi(p)}$ has first order contact
with $\sf$ at $p$.  This amounts to demanding that
\begin{subequations}
\label{eq:8}
\begin{align}
(\Pi,f)&=0\label{eq:9}\\(\Pi,\d f)&=0.\label{eq:10}
\end{align}
\end{subequations}
It can be shown that under mild (open) conditions, a congruence of
hyperspheres has two enveloping hypersurfaces.  In higher
co-dimension, there need not be \emph{any} enveloping submanifolds.

In view of \eqref{eq:8}, there is a close relationship between the
normal bundle $\Nf$ of an immersion $\sf$ and the sphere congruences
that envelope $\sf$.  Indeed, if $\Pi:M\to G_{n-k}^+$ is such a
congruence, then \eqref{eq:9} says that $\Pi\subset\sf^\perp$ and
then \eqref{eq:10} shows that projection along $\sf$ is a isomorphism
$\Pi\cong\Nf$.  Moreover, this isomorphism is parallel with respect
to the connection on $\Pi$ induced by flat differentiation in $\lor$
and $\nabla^\perp$ on $\Nf$.  In particular, the honest normal bundle
of a lift $f$ lying in some Riemannian model of $\PL$ gives an
enveloping sphere congruence.

\begin{eg}
Let $f:M\to E_{v_\infty}\subset\L^+$ be a lift lying in a copy of
Euclidean space and let $\Pi=\d f(TM)^\perp\subset
f^{-1}TE_{v_\infty}=\<f,v_\infty\>^\perp$.  Then $\Pi$ is an
enveloping sphere congruence and since $(\Pi,v_\infty)=0$ we see that
each sphere $S_{\Pi(p)}$ meets the point at infinity $\<v_\infty\>$.
Thus, after stereo-projection, $S_{\Pi(p)}$ is a plane and $\Pi$ is
the congruence $p\mapsto \d f(T_pM)$ of tangent planes to $f$.
\end{eg}

For a more substantial example, let $f:M\to\L^+$ be a lift of an
immersion of a $k$-dimensional manifold and let $\mcv_f$ be the mean
curvature vector of $f$:
\[
\mcv_f=\frac{1}{k}\trace\nabla\d f
\]
where $\nabla$ is the connection on $TM\otimes f^{-1}T\lor$ induced
by flat differentiation on $\lor$ and the Levi--Civita connection for
the metric $(\d f,\d f)$ on $M$ (the trace is, of course, computed
with respect to this metric also).
\begin{ex}
\begin{enumerate}
\item The sub-bundle $Z_{\sf}=\<f,\d f,\mcv_f\>^\perp\subset\lor$
depends only on $\sf$ and not on the choice of lift.
\item For $T\in\O$, $TZ_{\sf}=Z_{T\sf}$.
\end{enumerate}
\end{ex}
Moreover, for $e_1,\dots,e_k$ orthonormal with respect to $(\d f,\d
f)$, we have
\[
(f,\mcv_f)=-\frac{1}{k}(\d_{e_i} f,\d_{e_i} f)\neq 0
\]
whence, at each point $\<f,\d f,\mcv_f\>$ spans a $(k+2)$-plane on
which the inner product has signature $(k+1,1)$ so that each
$Z_{\sf}(p)$ is a space-like $(n-k)$-plane and $Z_{\sf}:M\to
G^+_{n-k}(\lor)$ is a M\"obius invariant enveloping sphere
congruence.  Geometrically, for a Euclidean lift $f:M\to
E_{v_\infty}$, this is the sphere congruence for which the sphere
tangent to $f$ at $p$ has the same mean curvature vector at $p$ as
$f$.

$Z_{\sf}$ is the \emph{central sphere congruence} \cite{Tho23} or
\emph{conformal Gauss map} \cite{Bry84} of $\sf$.\label{page:conf-gauss}

This construction comes alive when $k=2$ where it becomes a
fundamental tool in the theory of Willmore surfaces \cite{Bry84,Eji88}.
We shall meet this congruence again when we discuss Calapso's
approach \cite{Cal03,Cal15} to isothermic surfaces in
Section~\ref{sec:isoth-surf-via}.

\subsection{Clifford algebras in conformal geometry}
\label{sec:cliff-algebr-conf}

We are going to develop an extraordinarily efficient calculus for
conformal geometry using Clifford algebras that is especially well
adapted to working in the familiar Euclidean setting.  This will take
a little preparation so we begin by summarising the main idea.

We already know that the orthogonal group $\Op$ is isomorphic to the
M\"obius group $\Mob$.  We shall take a double cover and work instead
with an open subgroup of $\Pin$ which lies in the Clifford algebra
$\Cl_{n+1,1}$ of $\lor$.  A priori, it is not so clear why this is a
useful strategy.  However, there is a simple isomorphism of algebras
between $\Cl_{n+1,1}$ and the algebra $\Cl_n(2)$ of $2\times2$
matrices in the Clifford algebra $\Cl_n$ of $\R^n$.  The image under
this isomorphism of $\Pin$ is identified by a theorem of Vahlen
\cite{Vah02}.  Using this model, conformal diffeomorphisms of $\R^n$
become linear fractional transformations and the method of the moving
frame simplifies massively as one only has to do with $2\times 2$
matrices rather than the $(n+2)\times (n+2)$ matrices of the $\O$
formulation.

A good general reference for Clifford algebras is the text of
Michelsohn--Lawson \cite[Chapter 1]{LawMic89} while a clear account
of the relation between Clifford algebras and M\"obius
transformations can be found in the monograph of Porteous\footnote{In
fact, our approach differs slightly in the details from that in
\cite{Por95} since our conformal diffeomorphisms act on vectors (the
$\R^n$ that generates $\Cl_n$) rather than \emph{hypervectors}
(spanned by $1$ and some $\R^{n-1}\subset\R^n$).  In this we have
followed \cite{Wad90}.}  \cite[Chapters 18 and 23]{Por95}.

\subsubsection{Clifford algebras}
\label{sec:clifford-algebras}

Let $\R^{p,q}$ denote a $(p+q)$-dimensional vector space equipped
with an inner product of signature $(p,q)$ (that is, $p$ positive
directions and $q$ negative ones) and let $\Cl_{p,q}$ denote its
Clifford algebra.  Thus $\Cl_{p,q}$ is an associative algebra with
unit $1$ of dimension $2^{p+q}$ which contains $\R^{p,q}$ and is
generated by $\R^{p,q}$ subject only to the relations
\[
vw+wv=-2(v,w)1.
\]

$\Cl_{p,q}$ has a universal property which ensures the existence of
the following (anti-) involutions uniquely determined by their action
on the generators $\R^{p,q}$:
\begin{enumerate}
\item $a\mapsto\tilde{a}$: the \emph{order involution} with
$\tilde{v}=-v$ for $v\in\R^{p,q}$.
\item $a\mapsto a^t$: the \emph{transpose anti-involution} with
$v^t=v$ for $v\in\R^{p,q}$.
\item $a\mapsto\bar{a}$: the \emph{conjugate anti-involution} with
$\bar{v}=-v$ for $v\in\R^{p,q}$.
\end{enumerate}

\begin{ex}
For $a\in\Cl_{p,q}$, $\tilde{a}=\bar{a}^t$.
\end{ex}

The invertible elements $\Cl^\times_{p,q}$ form a multiplicative
group which acts on $\Cl_{p,q}$ via the \emph{twisted
adjoint action}:
\[
\tAd(g)a=ga\tilde{g}^{-1}
\]
\begin{ex}
$\tAd:\Cl_{p,q}^\times\to\GL(\Cl_{p,q})$ is a representation.
\end{ex}
Inside $\Cl^\times_{p,q}$ we distinguish the \emph{Clifford group}
$\Gamma_{p,q}$ given by
\[
\Gamma_{p,q}=\{g\in\Cl^{\times}_{p,q}:
\tAd(g)\R^{p,q}\subset\R^{p,q}\}
\]
The twisted adjoint action therefore restricts to give a
representation of $\Gamma_{p,q}$ on $\R^{p,q}$.
\begin{fact}
$\Gamma_{p,q}$ is generated by
$\R_\times^{p,q}=\{v\in\R^{p,q}:(v,v)\neq0\}=\R^{p,q}\cap\Cl^\times_{p,q}$.
\end{fact}
\begin{ex}
For $v\in\R_\times^{p,q}$, $\tAd(v):w\mapsto
vw\tilde{v}^{-1}=-vwv^{-1}$ is reflection in the hyperplane
orthogonal to $v$.
\end{ex}
As a consequence, each $\tAd(g)\in\O[p,q]$, the orthogonal group of
$\R^{p,q}$, and $\tAd:\Gamma_{p,q}\to\O[p,q]$ is a homomorphism which
has all reflections in its image and so is surjective by the
Cartan--Dieudonn\'e theorem.  Moreover, $\ker\tAd=\R^\times$ so that
we have an exact sequence:
\[
0\to\R^\times\to\Gamma_{p,q}\stackrel{\tAd}{\to}\O[p,q]\to 0.
\]

For $g\in\Gamma_{p,q}$, set $N(g)=g\bar{g}$, the \emph{norm} of $g$.
Writing $g=v_1\dots v_n$, with each $v_i\in\R_\times^{p,q}$, we see
that
\begin{align*}
N(g)=g\bar{g}=(v_1\dots v_n)\overline{(v_1\dots v_n)}&=
(v_1\dots v_n)(\bar{v}_n\dots\bar{v}_1)\\
&=\prod_{i=1}^n(v_i,v_i)\in\R^\times
\end{align*}
since $v_i\bar{v}_i=-v_i^2=(v_i,v_i)$.  From this we learn:
\begin{ex}
\begin{enumerate}
\item $N:\Gamma_{p,q}\to\R^\times$ is a homomorphism.
\item For $g\in\Gamma_{p,q}$, $N(g)=N(\bar{g})$.
\end{enumerate}
\end{ex}

Now let $\Cl^0_{p,q}$, $\Cl^1_{p,q}$ denote the $+1$ and $-1$
eigenspaces respectively of the order involution so that
$\Cl_{p,q}=\Cl^0_{p,q}\oplus\Cl^1_{p,q}$ is a $\Z_2$-graded algebra.
We define subgroups $\Pin[p,q]$ and $\Spin[p,q]$ of $\Gamma_{p,q}$ by
\begin{align*}
\Pin[p,q]&=\{g\in\Gamma_{p,q}:N(g)=\pm1\}\\
\Spin[p,q]&=\Pin[p,q]\cap\Cl^0_{p,q}.
\end{align*}
Then we have exact sequences:
\begin{gather*}
0\to\Z_2\to\Pin[p,q]\to\O[p,q]\to 0\\
0\to\Z_2\to\Spin[p,q]\to\SO[p,q]\to 0
\end{gather*}
where $\SO[p,q]=\O[p,q]\cap\SL[p+q,\R]$.

The Lie algebra $\o[p,q]$ of $\Pin[p,q]$ is the commutator
$[\R^{p,q},\R^{p,q}]\subset\Cl^0_{p,q}$ which acts on $\R^{p,q}$ by
the derivative of $\tAd$:
\[
\xi\cdot v=\xi v-v\tilde{\xi}=[\xi,v]
\]
since $\tilde{\xi}=\xi$.

Before leaving these generalities, we record some simple facts that
will be useful later on:
\begin{ex}
\label{ex:3}
For $g\in\Gamma_{p,q}$, $g^t,\bar{g},g^{-1}$ are all collinear.
Deduce:
\begin{enumerate}
\item For $v\in\R_\times^{p,q}$, $w\mapsto vwv$ is a symmetric
endomorphism of $\R^{p,q}$;
\item For $d\in\Gamma_{p,q}$, $d^td\in\R^\times$ and $w\mapsto d^twd$
is a conformal automorphism of $\R^{p,q}$: for $(v,w)\in\R^{p,q}$,
\[
(d^twd,d^twd)=(d^td)^2(v,w).
\]
\end{enumerate}
\end{ex}

\subsubsection{Vahlen matrices}
\label{sec:vahlen-matrices}

We now specialise to the case $(p,q)=(n+1,1)$ and arrive at the whole
point of our application of Clifford algebras: write $\Cl_n$ for
$\Cl_{n,0}$, the Clifford algebra of Euclidean $\R^n$ and contemplate
the algebra $\Cl_n(2)$ of $2\times 2$ matrices with entries in
$\Cl_n$.  I claim that
\[
Cl_n(2)\cong\Cl_{n+1,1}.
\]
Since both algebras have dimension $2^{n+2}$, this amounts to finding
a $(n+2)$-dimensional subspace $V$ of $\Cl_n(2)$ such that:
\begin{enumerate}
\item $v^2=-Q(v)I$ for all $v\in V$ where $I$ is the unit (identity
matrix) in $\Cl_n(2)$ and $Q$ is a quadratic form of signature
$(n+1,1)$;
\item $V$ generates $\Cl_n(2)$.
\end{enumerate}
For this, we take
\[
V=\left\{
\begin{pmatrix}
x&\lambda\\\mu&-x
\end{pmatrix}:\text{$x\in\R^n$, $\lambda,\mu\in\R$}\right\}
\]
and observe that
\[
\begin{pmatrix}
x&\lambda\\\mu&-x
\end{pmatrix}^2=
\begin{pmatrix}
-x^2+\lambda\mu&0\\0&-x^2+\lambda\mu
\end{pmatrix}=(-x^2+\lambda\mu)I.
\]
Thus we have light-like vectors $v_0,v_\infty\in V$ with
$(v_0,v_\infty)=-\half$ given by
\[
v_0=
\begin{pmatrix}
0&0\\1&0
\end{pmatrix},\qquad v_\infty=
\begin{pmatrix}
0&1\\0&0
\end{pmatrix}
\]
and $V$ therefore has an inner product of signature $(n+1,1)$.
\begin{ex}
$V$ generates $\Cl_n(2)$.
\end{ex}
This establishes the claim and henceforth we shall write $\lor$ for
$V\subset\Cl_n(2)$.

In fact, we have a little more: the decomposition of $\Cl_n(2)$ into
diagonal and off-diagonal matrices gives us a decomposition
\[
\lor=\R^n\oplus\R^{1,1}
\]
and fixed light-vectors $v_0,v_\infty\in\R^{1,1}$ lying in a
component $\L^+$ of $\L$.  Conversely, each such decomposition of
$\lor$ with chosen light-vectors in $\R^{1,1}$ gives us an
isomorphism $\Cl_n(2)\cong\Cl_{n+1,1}$.

The distinguished light-vectors $v_0,v_\infty$ give us a ready-made
stereographic projection
$\PL\setminus\<v_\infty\>\to\R^n=\<v_0,v_\infty\>^\perp$.  Indeed,
\[
E_{v_\infty}=\{v\in\L:(v,v_{\infty})=-\half\}= \left\{
\begin{pmatrix}
x&-x^2\\1&-x
\end{pmatrix}:x\in\R^n\right\}
\]
and the stereo-projection of \eqref{eq:5} reads
\[
\begin{pmatrix}
x&-x^2\\1&-x
\end{pmatrix}\mapsto
\begin{pmatrix}
x&0\\0&-x
\end{pmatrix}=x\in\R^n
\]
with inverse
\[
x\mapsto\begin{pmatrix} x&-x^2\\1&-x
\end{pmatrix}.
\]

The various (anti-)involutions on $\Cl_n(2)$ are readily identified:
\begin{ex}
For $\displaystyle
\begin{pmatrix}
a&b\\c&d
\end{pmatrix}\in\Cl_n(2)\cong\Cl_{n+1,1}$,
\begin{align*}
\begin{pmatrix}
a&b\\c&d
\end{pmatrix}^-&=
\begin{pmatrix}
d^t&-b^t\\-c^t&a^t
\end{pmatrix}\\
\begin{pmatrix}
a&b\\c&d
\end{pmatrix}^t&=
\begin{pmatrix}
\bar{d}&\bar{b}\\\bar{c}&\bar{a}
\end{pmatrix}\\
\begin{pmatrix}
a&b\\c&d
\end{pmatrix}^\sim&=
\begin{pmatrix}
\tilde{a}&-\tilde{b}\\-\tilde{c}&\tilde{d}
\end{pmatrix}
\end{align*}
\end{ex}

Now let $g=\displaystyle\begin{pmatrix}
a&b\\c&d\end{pmatrix}\in\Gamma_{n+1,1}$.  Since
$N(g)=g\bar{g}\in\R^\times$ we deduce from
\[
N(g)=\begin{pmatrix} a&b\\c&d
\end{pmatrix}
\begin{pmatrix}
d^t&-b^t\\-c^t&a^t
\end{pmatrix}=
\begin{pmatrix}
ad^t-bc^t&ba^t-ab^t\\cd^t-dc^t&da^t-cb^t
\end{pmatrix}\in\R^\times I
\]
that
\begin{gather*}
ad^t-bc^t=da^t-cb^t\in\R^\times\\
cd^t=dc^t\qquad ab^t=ba^t.
\end{gather*}
Moreover, $N(g)=N(\bar{g})$ gives
\begin{gather*}
ad^t-bc^t=d^ta-b^tc\\
c^ta=a^tc\qquad d^tb=b^td.
\end{gather*}
These are all necessary conditions for the matrix $g$ to lie in
$\Gamma_{n+1,1}$.  The full story is the content of Vahlen's theorem:
\begin{thm}\label{th:3}
$\displaystyle\begin{pmatrix} a&b\\c&d\end{pmatrix}\in\Gamma_{n+1,1}$
if and only if $a,b,c,d\in\Gamma_n\cup\{0\}$ with
\begin{enumerate}
\item $ad^t-bc^t\in\R^\times$;
\item $ac^t,bd^t,a^tb,c^td\in\R^n$.
\end{enumerate}
\end{thm}
\begin{ex}
For $a,c\in\Gamma_n\cup\{0\}$, $ac^t\in\R^n$ if and only if
$a^tc\in\R^n$.  Then take transposes to get $ca^t,c^ta\in\R^n$ also.
\end{ex}

We now restrict attention to the open subgroup $\SL$ of $\Pin$ given
by
\[
\SL=N^{-1}\{1\}.
\]
This has two components $\SL\cap\Cl_{n+1,1}^0$ and
$\SL\cap\Cl_{n+1,1}^1$ and double covers $\Op\cong\Mob$.
\begin{rem}
$\displaystyle\begin{pmatrix}a&b\\c&d\end{pmatrix}\in\Cl^0_{n+1,1}$
if and only if $a,d\in\Cl^0_n$ and $b,c\in\Cl^1_n$.
\end{rem}

Our formalism gives a beautiful description of the action of $\SL$ on
$\R^n$ by linear fractional transformations:
$g=\displaystyle\begin{pmatrix}a&b\\c&d\end{pmatrix}\in\SL$ induces a
conformal diffeomorphism of $\R^n\cup\{\infty\}$ which we denote by
$x\mapsto g\cdot x$.  To compute this, note that $g\bar{g}=1$ whence
$\tilde{g}^{-1}=g^t$ so that
\[
\tAd(g)v=gvg^t.
\]
Embedding $\R^n$ as usual into $\L^+$ by inverse stereo-projection,
\[
x\mapsto
\begin{pmatrix}
x&-x^2\\1&-x
\end{pmatrix}
\]
we have
\begin{align*}
\tAd(g)\begin{pmatrix} x&-x^2\\1&-x
\end{pmatrix}&=
\begin{pmatrix}
a&b\\c&d
\end{pmatrix}
\begin{pmatrix}
x&-x^2\\1&-x
\end{pmatrix}
\begin{pmatrix}
\bar{d}&\bar{b}\\\bar{c}&\bar{a}
\end{pmatrix}\\&=
\begin{pmatrix}
(ax+b)\overline{(cx+d)}&(ax+b)\overline{(ax+b)}\\(cx+d)\overline{(cx+d)}&(cx+d)\overline{(ax+b)}
\end{pmatrix}.
\end{align*}
\begin{ex}
For $x\in\R^n$ and $c,d\in\Gamma_n\cup\{0\}$,
$cx+d\in\Gamma_n\cup\{0\}$ and, in particular,
$(cx+d)\overline{(cx+d)}\in\R$.
\end{ex}
In the case at hand, either $cx+d=0$ in which case
\[
\tAd(g)
\begin{pmatrix}
x&-x^2\\1&-x
\end{pmatrix}=(ax+b)\overline{(ax+b)}
\begin{pmatrix}
0&1\\0&0
\end{pmatrix}\in\<v_\infty\>
\]
so that $g\cdot x=\infty$ or else
\[
\tAd(g)
\begin{pmatrix}
x&-x^2\\1&-x
\end{pmatrix}=(cx+d)\overline{(cx+d)}
\begin{pmatrix}
(ax+b)(cx+d)^{-1}&*\\1&*
\end{pmatrix}
\]
with stereo-projection $(ax+b)(cx+d)^{-1}\in\R^n$.
\begin{ex}
Show that
\[
\tAd(g)v_\infty=
\begin{cases}
\displaystyle c\bar{c}
\begin{pmatrix}
ac^{-1}&*\\1&*
\end{pmatrix}&\text{if $c\neq 0$;}\\
a\bar{a}v_\infty&\text{if $c=0$.}
\end{cases}
\]
Otherwise said, $g\cdot\infty=ac^{-1}\in\R^n\cup\{\infty\}$.
\end{ex}

To summarise: the action of $g$ as a conformal diffeomorphism of
$\R^n\cup\{\infty\}$ is given by
\[
g\cdot x=(ax+b)(cx+d)^{-1}.
\]

\begin{eg}
$\displaystyle
\begin{pmatrix}
0&-1\\1&0
\end{pmatrix}\in\SL$ acts by $x\mapsto -x^{-1}=x/\norm{x}^2$: this is inversion
in the unit sphere.
\end{eg}

Having understood the groups involved, let us briefly consider the
Lie algebra $\o=[\lor,\lor]$.
\begin{ex}
Show that
\[
\o=\left\{
\begin{pmatrix}
\xi&x\\y&-\xi^t
\end{pmatrix}:x,y\in\R^n,\xi\in[\R^n,\R^n]\oplus\R\right\}
\]
\end{ex}
Note that the decomposition of $\o$ into diagonal and off-diagonal
pieces,
\[
\o=\k\oplus\p,
\]
$\k=[\R^n,\R^n]\oplus\R$, $\p=\R^n\oplus\R^n$, is a symmetric
decomposition.  Indeed, $\k$, $\p$ are, respectively the $+1$ and
$-1$ eigenspaces of the involution in $\Op$ which is $+1$ on
$\R^{1,1}$ and $-1$ on $\R^n$.  The corresponding symmetric space will
play a starring role in Section~\ref{sec:curved-flats}.

\begin{rem}
We have confined our exposition to the case $\lor$ of direct
relevance to the theory we wish to develop.  However, the analogous
theory holds for any $\R^{p+1,q+1}$.  Again
$\Cl_{p+1,q+1}\cong\Cl_{p,q}(2)$ and the analog of Vahlen's theorem
identifies $\Gamma_{p+1,q+1}$ (with the refinement that
$\Gamma_n\cup\{0\}$ is replaced by the monoid generated by all
elements of $\R^{p,q}$ whether invertible or not).  Again, the
projective light cone in $\R^{p+1,q+1}$ is the conformal
compactification of $\R^{p,q}$ and we arrive at a description of the
conformal group of $\R^{p,q}$ in terms of $2\times2$ matrices with
entries in $\Cl_{p,q}$ and linear fractional transformations
\cite{Por95}.  It is hard not to hope that the methods elaborated
here may have applications in this more general setting.  A good test
case for this would be to take $(p,q)=(3,1)$ where these ideas
describe the symmetry group $\O[4,2]$ of Lie sphere geometry \cite{Cec92}.
\end{rem}

\subsubsection{Moving frames}
\label{sec:moving-frames}

We have now arrived at the model of conformal geometry with which we
shall work for the rest of this paper.  Let us summarise this
picture: we work with the ``Euclidean'' model $\R^n\cup\{\infty\}$ of
the conformal $n$-sphere using stereo-projection \eqref{eq:5} to
identify $\R^n\cup\{\infty\}$ with
$E_{v_\infty}\cup\{v_\infty\}\subset\L^+$ and so, via $\pi$, with
$\PL$.  The projective action of $\SL$ on $\PL$ induces an
action\footnote{Thus the action on $E_{v_\infty}\cup\{v_\infty\}$ is
the linear action $\tAd$ on $\L$ followed by rescaling to ensure that
the end result lies in $E_{v_\infty}\cup\{v_\infty\}$.} on
$E_{v_\infty}\cup\{v_\infty\}$ and so on $\R^n\cup\{\infty\}$ by
conformal diffeomorphisms.  We have seen that this action on
$\R^n\cup\{\infty\}$ is given by
\begin{subequations}
\label{eq:11}
\begin{align}
g\cdot x&=(ax+b)(cx+d)^{-1}\\
g\cdot\infty&=ac^{-1},
\end{align}
\end{subequations}
for $g=\displaystyle
\begin{pmatrix}
a&b\\c&d
\end{pmatrix}\in\SL$.

In what follows, we shall study maps $f:M\to\R^n$ and also pairs of
maps $f,\fh:M\to\R^n$.  A useful technique for this is the method of
the moving frame: a \emph{frame} for $f$ is a map $F:M\to\SL$ such
that
\[
f=F\cdot0.
\]
\begin{eg}
$F=\displaystyle
\begin{pmatrix}
1&f\\0&1
\end{pmatrix}$ frames $f$.
\end{eg}

Similarly, a frame for the (ordered) pair $(f,\fh)$ is a map
$F:M\to\SL$ such that
\[
f=F\cdot0\qquad\fh=F\cdot\infty.
\]
In this case, with $F=\displaystyle
\begin{pmatrix}
a&b\\c&d
\end{pmatrix}$, we have
\[
f=bd^{-1}\qquad \fh=ac^{-1}
\]
so that
\begin{equation}
\label{eq:12}
F=
\begin{pmatrix}
\fh c&fd\\c&d
\end{pmatrix}
\end{equation}
and the determinant condition of Theorem~\ref{th:3} reads
\begin{equation}
\label{eq:13}
(\fh-f)cd^t=1.
\end{equation}
In fact, once this condition is satisfied, $F$ defined by
\eqref{eq:12} automatically satisfies the remaining conditions of
Vahlen's Theorem and so lies in $\SL$:
\begin{ex}
If $f$ and $\fh$ never coincide\footnote{This is a necessary
condition for the pair to be framed since $g\cdot0\neq g\cdot\infty$
for any $g\in\SL$.} and $c,d\in\Gamma_n$ satisfy \eqref{eq:13}, then
$F$ defined by \eqref{eq:12} lies in $\SL$.
\end{ex}
\begin{eg}\label{eg:1}
The pair $(f,\fh)$ is framed by $\displaystyle
\begin{pmatrix}
\fh(\fh-f)^{-1}&f\\(\fh-f)^{-1}&1
\end{pmatrix}$.
\end{eg}

The point of using frames is that maps into a group are essentially
determined by their derivative.  For $F:M\to\SL$, consider the
\emph{\MC} form of $F$ given by $B=\mc{F}$: this is a $1$-form with
values in $\o$.  Differentiating the $\Cl_n(2)$-valued equation
\[
\d F=FB
\]
gives the \MC\ equations
\begin{equation}
\label{eq:14}
\d B+B\wedge B=0
\end{equation}
where multiplication in $\Cl_n(2)$ is used to multiply the
coefficients of $B$ in $B\wedge B$.

Conversely, given such a $1$-form $B$ satisfying \eqref{eq:14}, we can
locally\footnote{That is, on simply connected subdomains of $M$.}
integrate \cite{Ste83} to find a map $F:M\to\SL$ with $\mc{F}=B$ which
is unique up to left multiplication by constants in $\SL$.  The \MC\ 
equations \eqref{eq:14} amount to ``structure equations'' for the
immersions framed by $F$.
\begin{ex}\label{ex:4}
Set $g=\fh-f$ and put $F=\displaystyle
\begin{pmatrix}
\fh g^{-1}&f\\g^{-1}&1\end{pmatrix}$: we have seen that $F$ frames
$(f,\fh)$.  Show that
\[
\mc{F}=
\begin{pmatrix}
(\d f)g^{-1}&\d f\\-\gi(\d\fh)\gi&-\gi\d f
\end{pmatrix}.
\]
\end{ex}

\subsection{Exterior calculus on $\Omega\otimes\Cl_n$ and applications}
\label{sec:exter-calc-omeg}

\subsubsection{Clifford algebra valued differential forms}
\label{sec:cliff-algebra-valu}

Let $M$ be a manifold and $\Omega$ the exterior algebra of
differential forms on $M$.  Consider the space $\Omega\otimes\Cl_n$
of $\Cl_n$-valued forms on $M$.  Since $\Cl_n$ is an associative
algebra, we may extend exterior multiplication to
$\Omega\otimes\Cl_n$ by using the product in $\Cl_n$ to multiply
coefficients.  Thus for monomials
$a\omega_1,b\omega_2\in\Omega\otimes\Cl_n$ with $a,b:M\to\Cl_n$,
$\omega_i\in\Omega$:
\[
a\omega_1\wedge b\omega_2=(ab)\omega_1\wedge\omega_2.
\]
In particular, for $f,g\in\Omega^0\otimes\Cl_n$, that is,
$f,g:M\to\Cl_n$, the exterior product is just pointwise
multiplication.

Similarly, we extend the exterior derivative by
\[
\d(a\omega)=\d a\wedge\omega+a\d\omega.
\]
Since $\Cl_n$ is not, in general, commutative, exterior
multiplication on $\Omega\otimes\Cl_n$ is no longer super-commutative:
\[
\alpha\wedge\beta\neq\pm\beta\wedge\alpha.
\]
However, it is not difficult to establish:
\begin{prop}
For $\alpha\in\Omega^p\otimes\Cl_n$, $\beta\in\Omega^q\otimes\Cl_n$
and $f\in\Omega^0\otimes\Cl_n$:
\begin{enumerate}
\item $\alpha f\wedge\beta=\alpha\wedge f\beta$ (this is a special
case of the associativity of $\wedge$).
\item $\tilde{\alpha}\wedge\tilde{\beta}=(\alpha\wedge\beta)^\sim$.
\item $\alpha^t\wedge\beta^t=(-1)^{pq}(\beta\wedge\alpha)^t$.
\item $\bar{\alpha}\wedge\bar{\beta}=(-1)^{pq}\overline{(\beta\wedge\alpha)}$.
\item
$\d(\alpha\wedge\beta)=\d\alpha\wedge\beta+(-1)^p\alpha\wedge\d\beta$.
\item $\d^2=0$.
\end{enumerate}
\end{prop}

\begin{ex}
If $g:M\to\Cl^\times_n$, differentiate $g\gi=1$ to conclude:
\[
\d\gi=-\gi(\d g)\gi.
\]
\end{ex}

This exterior calculus will be our main computational tool for much
of these lectures.

\subsubsection{A lemma on commuting forms in $\Omega^1\otimes\R^n$}
\label{sec:lemma-comm-forms}

With an eye to a basic application to isothermic surfaces, we prove:
\begin{lem}
\label{th:4}
Let $V$ be a real vector space with $\dim V\geq2$ and
$\alpha,\beta:V\to\R^n$ non-zero linear maps with $\alpha$ injective.
Consider $\alpha\wedge\beta:\bigwedge^2V\to\Cl_n$.  Then
\begin{equation}
\label{eq:15}
\alpha\wedge\beta=0
\end{equation}
if and only if the following conditions are satisfied:
\begin{enumerate}
\item $\dim V=2$;
\item There is $\lambda\in\R^+$ such that
$(\beta,\beta)=\lambda(\alpha,\alpha)$;
\item $\im\alpha=\im\beta$;
\item $\det(\alpha^{-1}\circ\beta)<0$.
\end{enumerate}
Thus $\alpha$ and $\beta$ have the same image, induce conformally
equivalent inner products on $V$ but opposite orientations.
\end{lem}
\begin{proof}
Suppose first that $\alpha\wedge\beta=0$.  Choose an orthonormal
basis $e_1,\dots,e_n$ of $\R^n$ so that $\im\alpha=\<e_1,\dots,e_m\>$
and a basis $\omega_1,\dots,\omega_m$ of $V^*$ so that
\[
\alpha=\sum_{i\leq m}e_i\otimes\omega_i.
\]
Write
\[
\beta=\sum_{j\leq m}e_j\otimes\eta_j+\sum_{l>m}e_l\otimes\eta_l
\]
for some $\eta_1,\dots,\eta_n\in V^*$.  Then \eqref{eq:15} reads
\begin{equation}
\label{eq:16}
\sum_{i,j\leq m}e_ie_j\omega_i\wedge\eta_j+
\sum_{\substack{i\leq m\\l>m}}e_ie_l\omega_i\wedge\eta_l=0.
\end{equation}
The elements $1, e_ie_j$ ($i<j$) are linearly dependent in $\Cl_n$
while $e_ie_j=-e_je_i$, for $i\neq j$, whence taking coefficients in
\eqref{eq:16} gives:
\begin{subequations}
\begin{align}
\sum_{i\leq m}\omega_i\wedge\eta_i&=0&\quad&\label{eq:17}\\
\omega_i\wedge\eta_j&=\omega_j\wedge\eta_i&&\text{for $1\leq i<j\leq
m$.}
\label{eq:18}\\
\omega_i\wedge\eta_l&=0&&\text{for $1\leq i\leq m$,
$l>m$.}\label{eq:19}
\end{align}
\end{subequations}
{}From \eqref{eq:19}, we see that each $\eta_l=0$, for $l>m$, so that
$\im\beta\subset\im\alpha$.

Applying the Cartan Lemma to \eqref{eq:17}, we get, for each $j\leq
m$,
\[
\eta_j=\sum_{i\leq m}a_{ji}\omega_i
\]
with $a_{ij}=a_{ji}$.  Thus, \eqref{eq:18} becomes, for fixed
$i,j\leq m$,
\[
\sum_k a_{jk}\omega_i\wedge\omega_k=\sum_k
a_{ik}\omega_j\wedge\omega_k
\]
from which we conclude that $a_{ik}=0$ whenever $k\neq i,j$ and
$a_{ii}=-a_{jj}$.  If we can choose $i,j,k$ all distinct we quickly
conclude that all $a_{ij}=0$ so that $\beta=0$: a contradiction.  We
must therefore have $\dim V=2$ and
\begin{align*}
\eta_1&=a_{11}\omega_1+a_{12}\omega_2\\
\eta_2&=a_{12}\omega_1-a_{11}\omega_2
\end{align*}
with $a_{11}$ and $a_{12}$ not both zero.  We now have
\begin{align*}
(\beta,\beta)&=\eta_1^2+\eta_2^2\\
&=(a_{11}^2+a_{12}^2)(\omega_1^2+\omega_2^2)=(a_{11}^2+a_{12}^2)(\alpha,\alpha)
\end{align*}
and
\[
\det( \alpha^{-1}\circ\beta)=\det(a_{ij})=-a_{11}^2-a_{12}^2<0.
\]

The converse is more direct: let $\dim V=2$ and choose $v_1,v_2$ an
orthonormal basis of $V$ with respect to $(\alpha,\alpha)$ and set
$Z=v_1+iv_2\in V^{\C}$.  Thus $(\alpha(Z),\alpha(Z))=0$. Now
$\im\alpha=\im\beta$ gives $\beta(Z)\in\<\alpha(Z),\alpha(\bar{Z})\>$
while $(\beta,\beta)=\lambda(\alpha,\alpha)$ forces
$(\beta(Z),\beta(Z))=0$ so that $\beta(Z)$ is parallel to either
$\alpha(Z)$ or $\alpha(\bar{Z})$.  Finally,
$\det(\alpha^{-1}\circ\beta)<0$ forces the second possibility to hold
so that there is $\mu\in\C$ such that
\[
\beta(Z)=\mu\alpha(\bar{Z}),\quad\beta(\bar{Z})=\bar{\mu}\alpha(Z).
\]
Then
\begin{align*}
\alpha\wedge\beta(Z,\bar{Z})&=\alpha(Z)\beta(\bar{Z})-\alpha(\bar{Z})\beta(Z)\\
&=\bar{\mu}\alpha(Z)^2-\mu\alpha(\bar{Z})^2=0
\end{align*}
since $\alpha(Z)^2=-(\alpha(Z),\alpha(Z))=0$ and similarly for
$\alpha(\bar{Z})^2$.
\end{proof}

\subsubsection{More on sphere congruences}
\label{sec:more-sphere-congr}

Let $f,\fh:M\to\R^n$ be immersions of a $k$-dimensional manifold.  We
give a simple analytic condition for $f$ and $\fh$ to envelope the
same sphere congruence:
\begin{prop}
\label{th:5}
Let $g=\fh-f$.  Then $f$ and $\fh$ envelope the same sphere
congruence if and only if
\begin{equation}
\label{eq:20}
\im\d\fh=\im g\d f\gi.
\end{equation}
\end{prop}
\begin{proof}
The hypothesis \eqref{eq:20} means that $\im\d\fh=\im\rho_g\d f$
where $\rho_g=\tAd(g)$ is reflection in the hyperplane orthogonal to
$g$.

Now fix $p\in M$ and restrict attention to a $(k+1)$-dimensional
affine space\footnote{This space is uniquely determined unless
$g(p)\in\d f(T_pM)$.} containing $\fh(p)$ and $\d f(T_pM)+f(p)$.
Certainly, any $k$-sphere (or $k$-plane) tangent to $f$ and $\fh$ at
$p$ lie in this space.  Now any $k$-sphere containing $f(p)$ and
$\fh(p)$ must have centre on the hyperplane orthogonal to $g(p)$
through $\half(f(p)+\fh(p))$ and so is stable under reflection in this
hyperplane (which interchanges $f(p)$ and $\fh(p)$).  Moreover, there
is a unique $k$-sphere (or possibly $k$-plane) of this kind whose
tangent space at $f(p)$ is $\d f(T_pM)$.  The tangent space to this
sphere at $\fh(p)$ is therefore $\rho_g\d f(T_pM)$ which is tangent
to $\fh$ at $p$ if and only if $\rho_g\d f(T_pM)=\d\fh(T_pM)$. 
\end{proof}
\begin{ex}
In the situation of Proposition~\ref{th:5}, if $N$ is the unit
normal of $f$ pointing towards the centres of the sphere congruence,
then the radii $r$ of the spheres are given by
\[
1/r=-2(\gi,N).
\]
\end{ex}

\begin{rem}
$\rho_g$ restricts in this case to an isomorphism between the normal
bundles of $f$ and $\fh$ which, since $\rho_g$ is an isometry, is
parallel for the normal connections on those bundles.
\end{rem}
For future use, we record:
\begin{ex}\label{ex:5}
With $\im\d\fh=\im g\d f\gi$ and $N$ normal to $f$,
\[
\d(gN\gi)=g\bigl(\d N-2(\gi, N)\gi\d\fh g-2(\gi,N)\d f\bigr)\gi.
\]
\end{ex}

We conclude our present discussion of sphere congruences by
considering a degenerate case that we wish to exclude from further
discussions:
\begin{defn}
A sphere congruence $S:M\to G_{n-k}^+(\lor)$ is \emph{full} if there
is no fixed hyperplane $\Pi\subset\lor$ containing all the
$(n-k)$-planes $S(p), p\in M$.
\end{defn}
Let us contemplate non-full congruences.  The geometry of this
condition depends on the signature of the inner product when
restricted to $\Pi$:
\begin{enumerate}
\item If $\Pi$ has signature $(n,1)$, all spheres in the congruence
cut the hypersphere determined by the space-like line $\Pi^\perp$
orthogonally.
\item If $\Pi^\perp\in\PL$ (that is, $\Pi$ has signature $(n,0)$)
then all spheres in the congruence contain the point $\Pi^\perp$.
\item If $\Pi$ is space-like then all spheres in the congruence lie
totally geodesically in the round $n$-sphere determined by a unit
time-like vector in $\Pi^\perp$.
\end{enumerate}
Now restrict attention to the case where $\Pi$ has non-degenerate
inner product.  Then our non-full sphere congruence is stable under
the M\"obius transformation $R$ induced by reflection in $\Pi$.  As a
consequence, if $f$ envelopes the congruence, so does $R\circ f$.  We
give an analytic condition, refining that of Proposition~\ref{th:5},
for two enveloping surfaces to arise this way.
\begin{prop}
\label{th:6}
$(f,\fh)$ envelope a non-full sphere congruence with $\fh=R\circ f$
if and only if there is a function $\mu:M\to\R^\times$ such that
\begin{equation}
\label{eq:21}
\gi\d\fh\gi=\mu\d f.
\end{equation}
\end{prop}
\begin{proof}
Suppose first that $\fh=R\circ f$ where $R$ is the M\"obius
transformation induced by reflection in a non-degenerate hyperplane
$\Pi\subset\lor$.  Let $F:M\to\SL$ be the frame of $(f,\fh)$ given by
\[
F=
\begin{pmatrix}
\fh\gi&f\\\gi&1
\end{pmatrix}
\]
(cf. Exercise~\ref{ex:4}) and fix $v\in\Pi^\perp$ with $v^2=\pm 1$.
Up to a scaling, $\tAd(F)v_\infty$ is the reflection in $\Pi$ of
$\tAd(F)v_0$ whence $v\in\<\tAd(F)v_0,\tAd(F)v_\infty\>$ so that
\begin{equation}
\label{eq:22}
v=\Ad(F)
\begin{pmatrix}
0&\pm e^{-u}\\e^u&0
\end{pmatrix},
\end{equation}
for some $u:M\to\R$.  Since $v$ is constant, we have
\begin{align*}
0&=\tAd(F^{-1})\d v\\
&=\d
\begin{pmatrix}
0&\pm e^{-u}\\e^u&0
\end{pmatrix}+
\left[\mc{F},\begin{pmatrix}
0&\pm e^{-u}\\e^u&0
\end{pmatrix}\right]\\
&=
\begin{pmatrix}
0&\mp e^{-u}\d u\\e^u\d u&0
\end{pmatrix}+
\left[
\begin{pmatrix}
(\d f)\gi&\d f\\-\gi(\d\fh)\gi&-\gi\d f
\end{pmatrix},
\begin{pmatrix}
0&\mp e^{-u}\d u\\e^u\d u&0
\end{pmatrix}\right]\\
&=
\begin{pmatrix}
e^u\d f\pm e^{-u}\gi(\d\fh)\gi&
\mp e^{-u}(\d u -\d f\gi-\gi\d f)\\
e^u(\d u -\d f\gi-\gi\d f)&\mp e^{-u}\gi(\d\fh)\gi-e^u \d f
\end{pmatrix}
\end{align*}
where we have used Exercise~\ref{ex:4} to compute $\mc{F}$.

Thus we have two equations
\begin{subequations}
\begin{align}
e^{2u}\d f&=\mp \gi\d\fh\gi\label{eq:23}\\
\d u&=\{\gi,\d f\},\label{eq:24}
\end{align}
\end{subequations}
where $\{\,,\,\}$ is the anti-commutator in $\Cl_n$.  The first of
these is our desired equation~\eqref{eq:21} with
\begin{equation}
\label{eq:25}
\mu=\pm e^{2u}.
\end{equation}

Conversely, if \eqref{eq:21} holds, define $u$ by \eqref{eq:25} so
that \eqref{eq:23} holds and define $v$ by \eqref{eq:22}.  At each
point $p\in M$, $v(p)\in\<\tAd(F(p))v_0,\tAd(F(p))v_\infty\>$ so that
the enveloping sphere at $p$ is defined by a $(n-k)$-plane lying in
the hyperplane $\<v(p)\>^\perp$.  Moreover, reflection in this
hyperplane permutes the light-lines spanned by $\tAd(F(p))v_0$ and
$\tAd(F(p))v_\infty$ so that $\fh(p)=R_p(f(p))$ where $R_p$ is the
corresponding M\"obius transformation.  We will therefore be done if
we can show that $v$ is constant which, since \eqref{eq:23} holds by
construction, amounts to establishing \eqref{eq:24}.  However,
differentiating \eqref{eq:23} gives
\begin{align*}
2e^{2u}\d u\wedge\d f&=\pm\gi\d g\gi\wedge\d\fh\gi
\mp\gi\d\fh\wedge\gi\d g\gi\\
&=\mp\gi\d f\wedge\gi\d\fh\gi\pm\gi\d\fh\gi\wedge\d f\gi\\
&=e^{2u}\gi\d f\wedge\d f-e^{2u}\d f\wedge\d f\gi\\
&=e^{2u}(\gi\d f\wedge\d f+\d f\gi\wedge\d f
-\d f\wedge\gi\d f-\d f\wedge\d f\gi)\\
&=2e^{2u}\{\gi,\d f\}\wedge\d f.
\end{align*}
Equation~\eqref{eq:24} follows immediately and the proof is complete.
\end{proof}

\section{Isothermic surfaces: classical theory}
\label{sec:isoth-surf-class}

\subsection{Isothermic surfaces and their duals}
\label{sec:isoth-surf-their}

Let $f:M\to\R^{n}$ be an immersion of a surface $M$.  We begin with a
problem studied by Christoffel \cite{Chr67} for $n=3$ and Palmer
\cite{Pal88} for $n$ arbitrary: under what conditions is there a
second immersion $f^{c}:M\to\R^{n}$, a \emph{dual surface} of $f$,
such that:
\begin{enumerate}
\item $f$ and $f^{c}$ have parallel tangent planes: $\d f(T_{x}M)=\d
f^{c}(T_{x}M)$, for all $x\in M$;
\item $f$ and $f^{c}$ induce conformally equivalent metrics on $M$:
\[
(\d f,\d f)=\lambda(\d f^{c},\d f^{c}),
\]
for some $\lambda:M\to\R^{+}$.
\item $\d f^{-1}\circ\d f^{c}:TM\to TM$ is orientation-reversing:
$\det(\d f^{-1}\circ\d f^{c})<0$.
\end{enumerate}

In view of Lemma~\ref{th:4}, these conditions have a compact formulation in
our Clifford algebra formalism: viewing $\d f$ and $\d f^{c}$ as
$\Cl_{n}$-valued $1$-forms, they amount to
\[
\d f\wedge\d f^{c}=0.
\]
This motivates our main definition:
\begin{defn}
An immersion $f:M\to\R^{n}$ is \emph{isothermic} if there is a
non-constant map $f^{c}:M\to\R^{n}$ such that
\begin{equation}
\label{eq:26}
\d f\wedge\d f^{c}=0.
\end{equation}
\end{defn}

Note that, away from the zeros of $\d f^{c}$ (about which more
below), $f^{c}$ is a dual surface of $f$ and is itself isothermic
with dual surface $f$ since
\[
0=(\d f\wedge\d f^{c})^{t}=-\d f^{c}\wedge\d f.
\]

\begin{eg}
For $n=4$, $\Cl_{4}=\H(2)$ with $\R^{4}=\H$ embedded in $\H(2)$ via
\[
q\mapsto
\begin{pmatrix}
0&q\\-\bar{q}&0
\end{pmatrix}.
\]
Then, viewing $\d f$ and $\d f^{c}$ as $\H$-valued $1$-forms,
equation~\eqref{eq:26} reads
\[
\d f\wedge \d\bar{f^{c}}=0=\d\bar{f}\wedge\d f^{c}
\]
which is the characterisation of isothermic surfaces in $\R^{4}$
given by Hertrich-Jeromin--Pedit \cite{HerPed97}.
\end{eg}

Let $f:M\to\R^{n}$ be isothermic with dual $f^{c}$ and equip $M$ with
the conformal structure induced by $f$.  Define a quadratic
differential $Q_{f}:\otimes^{2}T^{1,0}M\to\C$ by $Q_{f}=(\d f,\d
f^{c})^{2,0}$.
\begin{lem}
$Q_{f}$ is a holomorphic quadratic differential.
\end{lem}
\begin{proof}
Choose a holomorphic coordinate $z$ on $M$.  We must show that
\[
(f_{z},f^{c}_{z})^{\vphantom{c}}_{\zbar}=0.
\]
As in Lemma~\ref{th:4}, there is a function $\mu$ with
\[
f^{c}_{z}=\mu f^{\vphantom{c}}_{\zbar},
\qquad f^{c}_{\zbar}=\bar{\mu}f^{\vphantom{c}}_{z}
\]
so that
\begin{align*}
(f^{\vphantom{c}}_{z},f^{c}_{z})^{\vphantom{c}}_{\zbar}&=
(f^{\vphantom{c}}_{z\zbar},f^{c}_{z})+(f^{\vphantom{c}}_{z},f^{c}_{z\zbar})\\
&=\mu(f_{z\zbar},f_{\zbar})+\bigl(f_{z},(\bar{\mu} f_{z})_{z}\bigr)\\
&=\half\mu(f_{\zbar},f_{\zbar})_{z}+\bar{\mu}_{z}(f_{z},f_{z})+
\half\bar{\mu}(f_{z},f_{z})_{z}=0
\end{align*}
since both $(f_{\zbar},f_{\zbar})$ and $(f_{z},f_{z})$ vanish by the
conformality of $f$.
\end{proof}

\begin{cor}
$Q_{f}$ and so $\d f^{c}$ vanish on at most a discrete set.
\end{cor}
Thus $f^{c}$ is at worst a branched conformal immersion.

We have now seen that an isothermic immersion $f:M\to\R^{n}$ equips
$M$ with a conformal structure and a non-zero holomorphic quadratic
differential $Q=Q_{f}\in\Gamma(\otimes^2 T^{1,0}M)$.  Otherwise said,
$(M,Q)$ is a \emph{polarised Riemann surface} in the sense of
\cite{HerMusNic}.

Moreover, we can recover $\d f^c$ from $f$ and this data: for any
holomorphic coordinate $z$ on $M$, write $Q=q\d z^2$, $f^c_z=\mu
f^{\vphantom{c}}_{\zbar}$ so that
\[
q=(f^{\vphantom{c}}_z,f^c_z)=\mu(f_z,f_{\zbar})
\]
whence
\begin{equation}
  \label{eq:27}
 f_z^c=qf_{\zbar}/(f_z,f_{\zbar}).
\end{equation}
Equation~(\ref{eq:27}) can be given an invariant formulation as
follows: for any map $g:M\to\R^n$ of a Riemann surface, write
\[
\d g=\del g+\delbar g
\]
where $\del g\in C^\infty(T^{1,0}M\otimes \C^n)$ and $\delbar
g=\overline{\del g}$ (thus, locally, $\del g=g_z \d z$).  Then
(\ref{eq:27}) reads:
\begin{equation}
  \label{eq:28}
  \del f^c=\frac{Q\delbar f}{(\d f,\d f)},
\end{equation}
where we have used tensor product to multiply powers of $T^{1,0}M$
and $T^{0,1}M$ and contraction to divide them.

To summarise: a conformal immersion $f$ of a polarised Riemann
surface $(M,Q)$ is isothermic with $Q_f=Q$ if and only if the
$1$-form
\[
\eta=\frac{1}{(\d f,\d f)}(Q\delbar f+\overline{Q}\del f)
\]
is exact\footnote{This latter condition is what Kamberov \cite{Kam97}
calls \emph{globally isothermic} when $n=3$.}.  Then $\d f^c=\eta$.

To make contact with the classical notion of an isothermic surface,
we compute the condition for the $1$-form $\eta$ to be closed: this
is
\[
\left(\frac{qf_{\zbar}}{(f_z,f_{\zbar})}\right)_{\zbar}=
\left(\frac{qf_{z}}{(f_z,f_{\zbar})}\right)_z.
\]
A short calculation using the holomorphicity of $q$ and the
conformality of $f$ reduces this to
\begin{equation}
  \label{eq:29}
  q(f_{\zbar\zbar})^\perp=\bar{q}(f_{zz})^\perp
\end{equation}
where ${}^\perp$ denotes the component in the normal bundle of $f$.
Away from the (isolated) zeros of $Q$, we may locally choose $z=x+iy$
so that $q=1$ and then (\ref{eq:29}) amounts to
\[
(f_{xy})^\perp=0
\]
so that $\del/\del x$ and $\del/\del y$ diagonalise the shape
operator $A^{N}$ of any normal $N$ to $f$.  We therefore conclude
that:
\begin{enumerate}
\item All shape operators of $f$ commute so that $f$ has flat normal
bundle;
\item $x,y$ are conformal curvature line (CCL) coordinates on $M$
(that is, conformal coordinates with respect to which each second
fundamental form is diagonal).
\end{enumerate}

These last two conditions constitute the classical definition of an
isothermic surface \cite{Cay72,Dar99} and in particular, we have
proved the following result of Christoffel ($n=3$) and Palmer:
\begin{thm}[\cite{Chr67,Pal88}]
Let $f$ have flat normal bundle and CCL coordinate $z=x+iy$ with $(\d
f,\d f)=e^{2u}\d z\d\zbar$.  Then the $\R^{n}$-valued $1$-form
defined by
\[
\eta=e^{-2u}\bigl(f_{\zbar}\d z+f_{z}\d\zbar\bigr)
\]
is closed and so locally is $\d f^{c}$ whence $f$ is isothermic with
dual $f^{c}$.
\end{thm}

How unique is the dual of an isothermic surface?  Certainly, if
$f^{c}$ is dual to $f$ then so is any $rf^c+k$ for constants
$r\in\R^{\times}$ and $k\in\R^{n}$ and then $Q_f$ becomes $rQ_f$.
With one interesting exception, these are the only possibilities: if
$f^{c}$ and $\tilde{f}^{c}$ are both duals of $f$ then, for any
holomorphic coordinate $z$, we have a function $\mu$ for which
\[
f^{c}_{z}=\mu\tilde{f}^{c}_{z},\qquad
f^{c}_{\zbar}=\bar{\mu}\tilde{f}^{c}_{\zbar}.
\]
Taking normal and tangential components of mixed derivatives of
$\tilde{f}^{c}$ gives
\begin{align*}
\mu f^{c}_{z\zbar}&=\bar{\mu}f^{c}_{z\zbar}\\
\mu_{\zbar}f_{z}^{c}&=\bar{\mu}_{z}f^{c}_{\zbar}
\end{align*}
so that $\mu$ is holomorphic.  Moreover $\mu$ is real (and so
constant) unless $f^{c}_{z\zbar}=0$, that is, unless $f^{c}$ is
minimal.  In this latter case, for any normal vector field $N$ to $f$
(and so $f^{c}$ also), we have
$(f^{c}_{\zbar},N^{\vphantom{c}}_{z})=0$ whence $(f_{z},N_{z})$
vanishes also and $f$ is totally umbilic.  Thus $f$ takes values in a
plane or $2$-sphere and $f^{c}$ is a minimal surface in $\R^{3}$.

We have therefore proved:
\begin{prop}\label{th:7}
Let $f:M\to\R^{n}$ be a full\footnote{that is, the image of $f$ is
not contained in any affine hyperplane.} isothermic surface.  Then the
dual $f^{c}$ of $f$ is unique up to scale and translations unless
$n=3$ and $f$ has image in a $2$-sphere.
\end{prop}

\subsubsection*{An example}

Let $f:M\to\R^{n}$ be an isometric immersion with mean curvature
vector $\mcv$, that is,
\[
\mcv=\half\trace\nabla\d f
\]
where $\nabla$ is the connection on $T^{*}M\otimes f^{-1}T\R^{n}$
induced by the Levi--Civita connections of $M$ and $\R^{n}$.

Following Chen \cite{Che73A}, a unit normal vector field $N$ of $f$
is said to be an \emph{isoperimetric section} if $(\mcv,N)$ is
constant and a \emph{minimal section} if $(\mcv,N)=0$.  Of course,
when $n=3$, $N$ is isoperimetric, respectively minimal, if and only
if $f$ has constant mean curvature, respectively is minimal, and this
motivates the following terminology:
\begin{defn}
A surface is said to be a \emph{generalised $H$-surface} if it admits
a parallel isoperimetric section.
\end{defn}

Generalised $H$-surfaces provide a class of examples of isothermic
surfaces in view of:
\begin{prop}
\label{th:8}
Let $f:M\to\R^{n}$ be an immersion and $N:M\to\R^{n}$ a unit normal
vector field not equal\footnote{This is to exclude the case where $f$
has image in a hyper-sphere and $N$ is the normal to that sphere.} to
any $rf+k$ for constants $r\in\R$, $k\in\R^{n}$.  Let $\phi=M\to\R$.
Then
\begin{enumerate}
\item $\phi N$ is dual to $f$ if and only if $\phi$ is constant and
$N$ is a parallel minimal section.
\item $f+\phi N$ is dual to $f$ if and only if $\phi$ is constant and
$N$ is a parallel isoperimetric section with
\[
(\mcv,N)=1/\phi.
\]
\end{enumerate}
\end{prop}
\begin{proof}
For $z$ a holomorphic coordinate on $M$,
\[
(\mcv,N)=-\frac{(N_{z},f_{\zbar})}{(f_{z},f_{\zbar})}.
\]
Now $f+\phi N$ is dual to $f$ if and only if $(f+\phi N)_{z}$ is
parallel to $f_{\zbar}$ or, equivalently,
\begin{subequations}
\begin{align}
(\phi_{z}N,N_{1})+\phi(N_{z},N_{1})&=0\label{eq:30}\\
(f_{z},f_{\zbar})+\phi (N_{z},f_{\zbar})&=0,\label{eq:31}
\end{align}
\end{subequations}
for any normal $N_{1}$ to $f$.  Taking $N_{1}=N$ in \eqref{eq:30}
gives $\phi_{z}=0$ and then \eqref{eq:30} asserts that $N$ is parallel
while \eqref{eq:31} asserts that $\phi=1/(\mcv,N)$ so that $N$ is
isoperimetric.

The case of minimal $N$ is similar.
\end{proof}

Let us collect some special cases of classical interest:
\begin{enumerate}
\item Let $f:M\to\R^{3}$ have constant mean curvature $H\neq 0$ and
Gauss map $N:M\to S^{2}$ so that $\mcv=HN$.  We see that $f$ is
isothermic with parallel dual surface
\[
f^{c}=f+\frac{1}{H}N
\]
which also has constant mean curvature $H$.  Moreover, the Hopf
differential of $f$ is $-HQ_f$.
\item If $f:M\to\R^{3}$ is minimal then $f$ is isothermic with its
Gauss map as dual surface and $Q_f$ is the negative of its Hopf
differential.
\item Let $N:M\to S^{n-1}\subset\R^{n}$ be isothermic then
Christoffel's formula~\eqref{eq:28} provides a dual surface $f$ for
which $N$ is a parallel minimal section.  For $n=3$, this is
particularly interesting since \emph{any} conformal map $M\to S^{2}$
is locally isothermic.  Indeed, if
$g:\Omega\subset\C\to\C\cup\{\infty\}$ is a meromorphic function,
inverse stereo-projection gives us a conformal map $N:M\to S^{2}$.
Moreover, since $S^{2}$ is totally umbilic, all directions are
curvature directions so that any holomorphic coordinate $z$ on
$\Omega$ is CCL.  Now let $f$ be holomorphic on $\Omega$ and set $Q=f
\d z^2$.  The formula~\eqref{eq:28} for the dual of $N$ now reads
\[
N^{c}_{z}=\frac{f}{g'}\bigl(\half(1-g^{2}),\tfrac{i}{2}(1+g^{2}),g\bigr)
\]
which we recognise as the Weierstrass--Enneper formula for the
minimal surface $N^{c}$ with Hopf differential $-f\d z^{2}$.

This explains the lack of uniqueness discussed in
Proposition~\ref{th:7}: the Weierstrass--Enneper formula requires a
choice of Hopf differential to prescribe a minimal surface with Gauss
map $N$.
\item Let $f:M\to S^{3}\subset\R^{4}$ be minimal with polar surface
$f^{\perp}:M\to S^{3}$ (thus $f^{\perp}\perp\< f,\d f(TM)\>$).  Then
$f^{\perp}$, viewed as a section of the normal bundle of $f$ in
$\R^{4}$, is certainly a parallel minimal section and so is dual to
$f$ in $\R^{4}$.  The symmetry of the situation ensures that
$f^{\perp}$ is minimal in $S^{3}$ also.
\item Similarly, if $f:M\to S^{3}$ has constant mean curvature $H\neq
0$ in $S^{3}$ then $f^{\perp}$ is a parallel isoperimetric section
with $(\mcv,f^{\perp})=H$ so that $f+\frac{1}{H}f^{\perp}$ is dual to
$f$.  Further, this parallel surface has constant mean curvature in a
sphere of radius $\sqrt{1+1/H^{2}}$.
\end{enumerate}

\subsection{Transformations of isothermic surfaces}
\label{sec:transf-isoth-surf}

A triumph in the classical study of isothermic surfaces was the
discovery by Darboux, Bianchi and Calapso
\cite{Bia05,Bia05A,Cal15,Dar99} of large families of transformations
of isothermic surfaces and permutability theorems relating these.
In this section, we shall show that this classical transformation
theory goes through unchanged in arbitrary co-dimension.

In all that follows, it will be convenient to fix a choice of dual
surface up to translation which amounts to fixing $Q_f$.  We
therefore fix a polarised Riemann surface $(M,Q)$ and refine our
basic definition:
\begin{defn}
A conformal immersion $f:(M,Q)\to\R^n$ is \emph{isothermic} if there
is a map $f^c:M\to\R^n$ with
\[
\d f\wedge \d f^c=0,\qquad (\d f,\d f^c)^{2,0}=Q.
\]
We say that $f^c$, which is unique up to translation, is the
\emph{Christoffel transform} of $f$.
\end{defn}

Thus $f^c$ is a particular choice of scaling for the dual surface of
$f$ and, away from its branch points, is isothermic on $(M,Q)$ with
Christoffel transform $f$.

\subsubsection{Conformal invariance}
\label{sec:conformal-invariance}

At first sight, our theory requires the notion of parallel tangent
planes and so is purely Euclidean.  However, at least locally, our
constructions are conformally invariant:
\begin{prop}
Let $f:(M,Q)\to\R^n$ be isothermic and $T\in\Mob$.  Then, locally,
$T\circ f$ is isothermic on $(M,Q)$.
\end{prop}
\begin{proof}
It suffices to check that the inversion
$f'=-f^{-1}:(M,Q)\to\R^n$ is isothermic.  Now
$\d f'=f^{-1}\d f f^{-1}$ and we put $\eta=f\d f^c f$.  Then
\begin{align*}
\d f'\wedge\eta&=f^{-1}\d f\wedge\d f^c f=0\\
(\d f',\eta)^{2,0}&=(\d f,\d f^c)^{2,0}=Q\\
\intertext{while} \d\eta&=\d f\wedge\d f^c f-f\d f^c\wedge\d f=0.
\end{align*}
Thus, locally, $\eta=\d(f')^c$ with $(f')^c$ the \Ch\ transform of
$f'$.
\end{proof}

\subsubsection{The Darboux transform}
\label{sec:darboux-transform}

Inspired by \HJ--Pedit \cite{HerPed97}, we begin with a (temporarily)
unmotivated definition:
\begin{defn}
Let $f,\fh:M\to\R^n$ be non-constant with $g=\fh-f$.  Say that $\fh$
is a \emph{Darboux transform} of $f$, or that $(f,\fh)$ are a
\emph{Darboux pair}, if
\begin{equation}
\label{eq:32}
\d f\wedge g^{-1}\d\fh g^{-1}=0.
\end{equation}
\end{defn}
Note that, in this case, $f$ is a Darboux transform of $\fh$ also.

Observe that if $(f,\fh)$ is a Darboux pair then
\begin{align*}
\d(g^{-1}\d\fh g^{-1})&=
\d g^{-1}\wedge\d\fh g^{-1}-g^{-1}\d\fh\wedge\d g^{-1}\\
&=-g^{-1}\d g\gi\wedge\d\fh\gi+\gi\d\fh\wedge\gi\d g\gi\\
&=\gi\d f\gi\wedge\d\fh\gi-\gi\d\fh\gi\wedge\d f\gi=0
\end{align*}
in view of~\eqref{eq:32} and its transpose.  Thus, locally,
$g^{-1}\d\fh g^{-1}=\d f^c$ for $f^c$ a dual surface to $f$.  Thus
$f$ is isothermic and so is $\fh$ with dual given by $\d\fh^c=\gi\d
f\gi$.

Moreover,
\[
(\d f^c,\d f^c)=( g^{-1}\d\fh g^{-1},g^{-1}\d\fh
g^{-1})=g^{-4}(\d\fh,\d\fh)
\]
so that $f$ and $f^c$ induce the same conformal structure on $M$.
Further, using Exercise~\ref{ex:3},
\[
(\d f,\d f^c)=(\d f,\gi\d\fh\gi)=(\gi\d f\gi,\d\fh)=(\d\fh^c,\d\fh)
\]
so that $Q_f=Q_{\fh}$ and $f$ and $\fh$ induce the same polarisation
on $M$ also.

To summarise:
\begin{thm}
If $(f,\fh)$ are a Darboux pair, then $f$ and $\fh$ are isothermic on
the same polarised Riemann surface.
\end{thm}

Now let $(M,Q)$ be a polarised Riemann surface and let
$f:(M,Q)\to\R^n$ be isothermic with \Ch\ transform $f^c$.  We seek
Darboux transforms of $f$.  If $\fh=f+g$ is a Darboux transform then,
from~\eqref{eq:32}, we have that $\gi\d\fh\gi$ is the derivative of a
dual surface to $f$ so that, for some $r\in\R^\times$,
\[
\d\fh=\d f+\d g=rg\d f^c g.
\]
We rearrange this into a Riccati equation for $g$:
\begin{equation}
\label{eq:33}
\d g=g\d f^c g-\d f.
\end{equation}
The integrability condition for \eqref{eq:33} is easily
checked\footnote{We shall see an illuminating proof below on
page~\pageref{page:darboux-system}.} to be the isothermic condition
\[
\d f\wedge\d f^c=0
\]
so that, for any initial condition, we may locally solve \eqref{eq:33}
for $g$ and then defining $\fh$ by $\fh=f+g$, we have
\[
\gi\d\fh\gi=r\d f^c
\]
so that $\fh$ is a Darboux transform of $f$.

\begin{notation}
For future use, fix a base-point $o\in M$ and let $f:(M,Q)\to\R^n$ be
isothermic with \Ch\ transform $f^c$.  We denote by $\D_r^v f$ the
Darboux transform $\fh=f+g$ where $g$ solves \eqref{eq:33} with
$\fh(o)=v$.

If we do not wish to emphasise the initial condition, we shall simply
write $\D_r f$.
\end{notation}

So let $\D_r f=\fh=f+g$ be a Darboux transform of $f:(M,Q)\to\R^n$.
The demand that $Q_{\fh}=Q$ fixes the \Ch\ transform of $\fh$ so that
\begin{equation}
\label{eq:34}
d\fh^c=r^{-1}\gi\d f\gi.
\end{equation}
On the other hand, a well-known symmetry of Riccati equations tells
us that $\gi$ must solve a Riccati equation also: indeed,
\begin{align*}
\d(rg)^{-1}&=-r^{-1}\gi\d g\gi\\
&=r^{-1}\gi\d f\gi-\d f^c\\
&=r(rg)^{-1}\d f(rg)^{-1}-\d f^c.
\end{align*}
Thus $(rg)^{-1}$ solves the $r$-Riccati equation for $f^c$ so that
$\widehat{f^c}=f^c+(rg)^{-1}=\D_r f^c$ with
\begin{equation}
\label{eq:35}
\d\widehat{f^c}=r^{-1}\gi\d f\gi.
\end{equation}
Comparing equations \eqref{eq:34} and \eqref{eq:35}, we see that
$\widehat{f^c}=\fh^c$ up to a translation and we have proved a
theorem due to Bianchi \cite{Bia05} when $n=3$ and \HJ--Pedit
\cite{HerPed97} when $n=4$:
\begin{thm}\label{th:9}
\Ch\ and Darboux transforms commute.
\end{thm}
Thus, once we have fixed the \Ch\ transform of $f$, we have a
\emph{unique} \Ch\ transform of any $\D_r f$ with all ambiguity of
scaling and translation removed.  Otherwise said, we may think of the
Darboux transform as a transform of \emph{\Ch\ pairs}
$(f,f^c)\mapsto(f+g,f^c+(rg)^{-1})$.

It is humbling to discover that the geometrical description of this
construction of the Darboux transform of $f^c$ was already known to
Bianchi even though he did not have the Riccati equation: let
$P,P_1,\bar{P},\bar{P}_1$ denote corresponding points on
$f,\fh,f^c,\widehat{f^c}$, he writes \cite[p.~105]{Bia05}:
\begin{quotation}
I segmenti $PP_1$, $\bar{P}\bar{P}_1$ sono paralleli ed il prodotto
delle loro lunghezze \`e constante $=2/m$.
\end{quotation}
(Our $r$ is Bianchi's $m/2$.)

\begin{ex}\label{ex:6}
Show that if $\fh=\D_r f$ then $f=\D^{f(o)}_r\fh$.  Thus $f=\D_r\D_r
f$.
\end{ex}

To justify our terminology and make contact with the classical
literature, we turn to the geometry of our constructions.  So let
$(f,\fh)$ be a Darboux pair.  In view of \eqref{eq:32} and Lemma~\ref{th:4},
we see that
\[
\im\d f=\im g\d\fh\gi
\]
so that Proposition~\ref{th:5} tells us that $f$ and $\fh$ are enveloping
surfaces of a $2$-sphere congruence $S$.  We have already seen that
$f$ and $\fh$ induce the same conformal structure on $M$: in
classical terms, this means that $S$ is a \emph{conformal sphere
congruence}.  Again, we know that $Q_f=Q_{\fh}$ so that $f$ and $\fh$
have the same curvature lines.  This condition was also well known in
the classical literature: recall that a $2$-sphere congruence $S$
induces a parallel isomorphism of the normal bundles of its
enveloping surfaces $f$ and $\fh$ via $N\mapsto gN\gi$.  The
congruence is \emph{Ribaucour} if corresponding normals have the
same principal directions, that is, if the shape operators $A^N$ and
$\hat{A}^{gN\gi}$ of $f$ and $\fh$ commute for each normal $N$ to $f$.

We therefore conclude: \emph{a Darboux pair consists of the
enveloping surfaces of a conformal Ribaucour congruence of
$2$-spheres.}

If we exclude degenerate cases, the converse is also true: recall
that a sphere congruence $S:M\to G^+_{n-2}(\R^{n+1,1})$ is
\emph{full} if its image is not contained in some fixed hyperplane.
\begin{ex}
Suppose that $S:M\to G^+_{n-2}(\R^{n+1,1})$ is not full and has an
enveloping surface $f:M\to\L$.  Reflect $f$ in the fixed
hyperplane\footnote{If this hyperplane has degenerate metric, take
the second surface to be constant.} to get a second enveloping
surface $\fh$ and so conclude that $S$ is conformal and Ribaucour.
\end{ex}

We now have:
\begin{prop}
Let $f,\fh$ envelope a full conformal Ribaucour $2$-sphere
congruence $S$ and suppose $f$ and $\fh$ have no umbilic points in
common.  Then $(f,\fh)$ is a Darboux pair.
\end{prop}
\begin{proof}\newcommand{\bh}{\hat{b}}
Let $N$ be normal along $f$ so that $gN\gi$ is normal along $\fh$.
The second fundamental form $\bh^{gN\gi}$ of $\fh$ along $gN\gi$ is
given by
\[
\bh^{gN\gi}_{U,V}=-\bigl(\d_U(gN\gi),\d_V\fh\bigr)
\]
and we know from Exercise~\ref{ex:5} that
\[
\d (gN\gi)=g\bigl( \d N-2(\gi,N)\gi\d\fh g-2(\gi,N)\d f \bigr)\gi
\]
whence
\[
\bh^{gN\gi}_{U,V}=
\Bigl(g\bigl(\d_U N-2(\gi,N)\d_U f\bigr)\gi,\d\fh_V \Bigr)-
2(\gi,N)\bigl( \d\fh_U,\d\fh_V \bigr).
\]
It is not difficult to check that if $(\d N,\d f)=2(\gi, N)(\d f,\d
f)$ at some point then the same identity is also true for $\fh$ at
that point:
\[
(\d(gN\gi),\d\fh)=-2(\gi,gN\gi)(\d\fh,\d\fh).
\]
Thus our exclusion of common umbilics prevents this possibility
occurring for all $N$.  So choose $N$ and principal vectors $X,Y$,
orthogonal with respect to $f$, such that the tangential component of
$\d_X N-2(\gi, N)\d_X f$ is a non-zero multiple of $d_X f$.  Since
$S$ is conformal and Ribaucour, $X,Y$ are orthogonal for $\fh$ and
principal for $gN\gi$ so that
\[
0=\bigl( g(\d_X N-2(\gi,N)\d_X f)\gi,\d_Y\fh \bigr)
\]
whence
\[
0=(d_X f,\gi\d_Y\fh g).
\]
We therefore conclude that there are are functions $\mu_1,\mu_2$ such
that
\begin{align*}
\d_X f&=\mu_1\gi\d_X\fh g\\
\d_Y f&=\mu_2\gi\d_Y\fh g.
\end{align*}
Since $(\d f,\d f)$ and $(\d\fh,\d\fh)$ are conformally equivalent,
we get $\mu_1^2=\mu_2^2$ and there are only two possibilities: either
$\mu_1=-\mu_2$ which quickly gives
\[
\d f\wedge \gi\d\fh g=0
\]
so that $(f,\fh)$ are a Darboux pair, or,
\[
\d f=\mu_1\gi\d\fh g
\]
which, by Proposition~\ref{th:6}, forces $S$ to be non-full.
\end{proof}
Thus, modulo umbilics, a Darboux pair is exactly a pair of enveloping
surfaces of a full conformal Ribaucour $2$-sphere congruence and,
for $n=3$, it is this latter formulation that Darboux gave
\cite{Dar99}.

Darboux's own construction of the Darboux transforms of a given
isothermic surface in $\R^3$ proceeded by solving a linear
differential system in $\R^{4,1}$ with the algebraic constraint that
the solution lie in the light cone $\L$.  It is instructive to
compare this approach with ours: for $(f,f^c)$ a \Ch\ pair,
contemplate the Lie algebra valued $1$-form $B\in\Omega^1\otimes\o$
given by\label{page:darboux-system}
\[
B=
\begin{pmatrix}
0&\d f\\r\d f^c&0
\end{pmatrix}.
\]
The \MC\ equations $\d B+\half[B\wedge B]=0$ reduce in this case to
the isothermic condition $\d f\wedge\d f^c=0$ so that the linear
differential system
\begin{equation}
\label{eq:36}
\d\omega+B\omega=0
\end{equation}
for $\omega:M\to\R^{n+1,1}$ is integrable (indeed, one integrates the
\MC\ equations to find $F:M\to\O$ with $F^{-1}\d F=B$ and then
solutions of \eqref{eq:36} are given by $\omega=F^{-1}\omega_0$ for
constant $\omega_0$).  Clearly, $(\omega,\omega)$ is an integral of
\eqref{eq:36} so that, in particular, $\L$ is preserved by the
integral flows.

The linear system \eqref{eq:36} with the constraint
$(\omega,\omega)=0$ is, up to gauge, the system considered by
Darboux\footnote{For a recent account of Darboux's approach see
\cite{BurHerPed97}.}.

Now let $\omega:M\to\L\subset\Cl_n(2)$ be a solution of \eqref{eq:36}
given by
\[
\omega=
\begin{pmatrix}
v&s\\t&-v
\end{pmatrix}
\]
and let $g:M\to\R^n\cup\{\infty\}$ be its stereo-projection.  Thus
$g=v/t$.  I claim that $g$ solves our Riccati equation \eqref{eq:33}.
Indeed, the action of $\o$ on $\R^{n+1,1}\subset\Cl_n(2)$ is by
commutator of Clifford matrices so that \eqref{eq:36} reads
\[
\begin{pmatrix}
\d v&\d s\\\d t&-\d v
\end{pmatrix}+
\left[
\begin{pmatrix}
0&\d f\\r\d f^c&0
\end{pmatrix},
\begin{pmatrix}
v&s\\t&-v
\end{pmatrix}
\right]=0
\]
from which we get
\begin{align*}
\d v&=sr\d f^c-t\d f\\
\d t&=-r\d f^c v-rv\d f^c
\end{align*}
whence
\begin{align*}
\d g=\frac{1}{t}\d v-\frac{\d t}{t^2}v&=\frac{s}{t}r\d f^c-\d f+r\d
f^c\frac{v^2}{t^2}+r\frac{v}{t}\d f^c\frac{v}{t}\\
&=rg\d f^c g-\d f
\end{align*}
where we have used $v^2+st=0$ (since $\omega$ is $\L$-valued).

We end our present discussion of the Darboux transform by
characterising the frames of Darboux pairs.  Recall that a frame of a
pair of maps $(f,\fh)$ is a map $F:M\to\SL$ such that
\[
f=F\cdot 0\qquad\fh=F\cdot\infty.
\]
We prove:
\begin{thm}
\label{th:10}
Let $F$ frame $(f,\fh)$ with
\[
\mc{F}=
\begin{pmatrix}
\alpha&\beta\\\gamma&\delta
\end{pmatrix}\in\Omega^1\otimes\o.
\]
Then $(f,\fh)$ is a Darboux pair if and only if
$\beta\wedge\gamma=0$.

In this case, if $f$ is isothermic with respect to a polarisation
$Q$, $\fh=\D_r f$ where $r$ is given by
\[
(\beta,\gamma)^{2,0}=-rQ.
\]
\end{thm}
\begin{proof}
The first thing to note is that the conditions on $\beta,\gamma$ are
independent of the choice of frame: if $\hF$ is another frame of
$(f,\fh)$ then $\hF=Fk$ for $k:M\to\SL$ with $k\cdot0=0$ and
$k\cdot\infty=\infty$.  Thus
\[
k=
\begin{pmatrix}
a&0\\0&a^{-t}
\end{pmatrix}
\]
for $a:M\to\Gamma_n$, and setting
\[
\mc{\hF}=
\begin{pmatrix}
\ha&\hb\\\hg&\hd
\end{pmatrix}
\]
we readily compute that
\[
\hb=a^{-1}\beta a^{-t},\qquad\hg=a^t\gamma a.
\]
Thus
\begin{align*}
\hb\wedge\hg&=a^{-1}(\beta\wedge\gamma)a\\
(\hb,\hg)&=-\half(a^{-1}\beta\gamma a+a^t\gamma\beta a^{-t})=(\beta,\gamma)
\end{align*}
where the last equality follows since $a^t$ and $a^{-1}$ are
collinear and so have the same adjoint action.

Thus we are free to choose a convenient frame to establish the
theorem.  With $g=\fh-f$, we take
\[
F=
\begin{pmatrix}
\fh\gi&f\\\gi&1
\end{pmatrix}
\]
so that, by Exercise~\ref{ex:4},
\[
\mc{F}=
\begin{pmatrix}
(\d f)g^{-1}&\d f\\-\gi(\d\fh)\gi&-\gi\d f
\end{pmatrix}.
\]
Thus $\beta=\d f$, $\gamma=-\gi(\d\fh)\gi$ and the vanishing of
$\beta\wedge\gamma$ is precisely the condition that $(f,\fh)$ is a
Darboux pair of isothermic surfaces.

Moreover, if $f:(M,Q)\to\R^n$ is isothermic with \Ch\ transform $f^c$
and $\fh=\D_rf$, we have
\[
\gi(\d\fh)\gi=r\d f^c
\]
so that
\[
(\beta,\gamma)^{2,0}=-r(\d f,\d f^c)^{2,0}=-rQ.
\]
\end{proof}
\begin{rem}
The geometric content of this result is that a Darboux pair is the
same as a \emph{curved flat} in the symmetric space $S^n\times
S^n\setminus\Delta$ of pairs of distinct points in $S^n$.  We shall
explain this in Section~\ref{sec:curved-flats}.
\end{rem}
\subsubsection{The $T$-transform}
\label{sec:t-transforms}

Our final family of transformations, discovered in the classical
setting by Calapso \cite{Cal03} and Bianchi \cite{Bia05A}, have a
slightly different flavour: the construction proceeds by solving a
\MC\ equation to build a frame of the new surface.  In particular,
these new surfaces are only defined up to the action of M\"obius
group.

We begin with an isothermic surface $f:(M,Q)\to\R^n$, its \Ch\ 
transform $f^c$ and a parameter $r\in\R$.  We have already seen that
the $\o$-valued $1$-form $B_r$ given by
\[
B_r=
\begin{pmatrix}
0&\d f\\r\d f^c&0
\end{pmatrix}
\]
solves the \MC\ equations so that, locally, we may integrate to find
$F_r:M\to\SL$ with $F_r^{-1}\d F_r=B_r$.  Of course, $F_r$ is only
determined up to left translation by a constant in $\SL$.  Now $F_r$
frames the pair $f_r,\fh_r:M\to\R^n\cup\{\infty\}$ given by
\[
f_r=F_r\cdot0,\qquad \fh_r=F_r\cdot\infty
\]
and, since
\[
\d f\wedge(r\d f^c)=0,\qquad (\d f,r\d f^c)^{2,0}=rQ,
\]
we immediately deduce from Theorem~\ref{th:10}:
\begin{thm}\label{th:11}
For $r\neq 0$, $(f_r,\fh_r)$ are a Darboux pair of isothermic
surfaces.  Moreover $f_r$ (and so $\fh_r$) are isothermic with
respect to $(M,Q)$ and
\[
\fh_r=\D_{-r}f_r.
\]
\end{thm}
We denote $f_r$ by $\T_r f$ and, following Bianchi, call it a
$T$-transform of $f$.  Note that $\T_r f$ is only determined up to
the action of $\Mob$.

When $r=0$ we may take
\begin{equation*}
F_0=
\begin{pmatrix}
1&f\\0&1
\end{pmatrix}
\end{equation*}
so that $f_0=f$ and $\fh_0\equiv\infty$.  Thus we take $\T_0 f=f$
modulo $\Mob$.

Our construction seems to depend in an essential way on the frames we
obtained by integrating $B_r$.  However, one can use any frame of $f$
as a starting point\footnote{I am grateful to Udo Hertrich-Jeromin
for explaining this point to me.}: indeed, any frame of $f$ is of the
form $\tF_0=F_0P$ where $P:M\to\SL$ has $P\cdot0=0$ and so is of the
form
\begin{equation*}
P=
\begin{pmatrix}
p_1&0\\p_2&p_3
\end{pmatrix}
\end{equation*}
with $p_1p_3^t=1$.

Then $\tF_r=F_r P$ frames $f_r$, that is, $F_rP\cdot0=F_r\cdot 0=f_r$.
Moreover,
\begin{align*}
\tF_r^{-1}\d\tF_r&=P^{-1} B_r P+P^{-1}\d P\\
&=\tF_0^{-1}\d\tF_0+rP^{-1}
\begin{pmatrix}
0&0\\\d f^c&0
\end{pmatrix}P
\end{align*}
and a short computation using $p_1p_3^t=1$ gives:
\begin{align*}
\tF_0^{-1}\d\tF_0&=
\begin{pmatrix}
*&p_3^t\d fp_3\\ *&*
\end{pmatrix}\\
P^{-1}\begin{pmatrix} 0&0\\\d f^c&0
\end{pmatrix}P&=
\begin{pmatrix}
0&0\\ p^{-1}_3\d f^c p^{-t}_3&0
\end{pmatrix}.
\end{align*}
The key point now is that $p^{-1}_3\d f^c p^{-t}_3$ is constructed
from $p_3^t\d fp_3$ in exactly the same way as $\d f^c$ is
constructed from $\d f$, that is, via \Ch's formula \eqref{eq:28}.
Indeed,
\begin{align*}
p^{-1}_3\del f^c p^{-t}_3&=\frac{1}{(\d f,\d f)}p_3^{-1}(Q\delbar f)p_3^{-t}\\
&=\frac{(p_3p_3^t)^2}{(p_3^t \d f p_3,p_3^t \d f p_3)}p_3^{-1}(Q\delbar
f)p_3^{-t}\\
&=\frac{1}{(p_3^t \d f p_3,p_3^t \d f p_3)}Q p_3^t \delbar f p_3.
\end{align*}
For $\alpha\in\Omega^1\otimes\R^n$ a conformal $1$-form, that is
$(\alpha,\alpha)^{2,0}=0$, write
\[
\alpha=\alpha'+\alpha''
\]
with $\alpha'\in\Omega^{1,0}\otimes\C^n$ and
$\overline{\alpha'}=\alpha''$ and define
$\alpha^c\in\Omega^1\otimes\R^n$ by
\[
\alpha^c=\frac{1}{(\alpha,\alpha)}(Q\alpha''+\overline{Q}\alpha').
\]
Our last calculation now reads
\[
(p_3^t \d fp_3)^c=p_3^{-1}\d f^c p_3^{-t}
\]
and we have proved
\begin{thm}\label{th:12}
Let $\tF$ frame an isothermic surface $f:(M,Q)\to\R^n$ with
\[
\tF^{-1}\d\tF=
\begin{pmatrix}
\alpha&\beta\\\gamma&\delta
\end{pmatrix}.
\]
Then
\[
\widetilde{B}_r=\begin{pmatrix} \alpha&\beta\\\gamma&\delta
\end{pmatrix}+
r\begin{pmatrix} 0&0\\\beta^c&0
\end{pmatrix}.
\]
solves the \MC\ equations and if $\tF_r^{-1}\d\tF_r=\widetilde{B}_r$
then $\tF_r$ frames $\T_r f$.
\end{thm}

As a first application, let us show that, in analogy with the Lie
transform of $K$-surfaces, $\T_r$ gives an action of $\R$ on
isothermic surfaces modulo $\Mob$.  Indeed, $F_r$ frames $f_r$ with
\[
F_r^{-1}\d F_r=
\begin{pmatrix}
0&\d f\\r\d f^c &0
\end{pmatrix}
\]
while
\[
\mc{F_{s+r}}=
\begin{pmatrix}
0&\d f\\(s+r)\d f^c& 0
\end{pmatrix}=
\mc{F_r}+s
\begin{pmatrix}
0&0\\\d f^c&0
\end{pmatrix}
\]
so that $f_{s+r}=\T_s f_r$ and we have a theorem proved by \HJMN\ 
\cite{HerMusNic} for the case $n=3$.
\begin{thm}
$\T_{s+r}=\T_s\circ\T_r$ modulo $\Mob$.
\end{thm}

Again, we can compare the $T$-transforms of $f$ and $f^c$: for
$r\in\R^\times$, define $R_r$ by
\[
R_r=
\begin{pmatrix}
0&\sign(r)/\sqrt{\abs{r}}\\\sqrt{\abs{r}}&0
\end{pmatrix}
\]
so that $R_r\cdot0=\infty$ and $R_r\cdot\infty=0$.
\begin{ex}
\label{ex:7}
\[
\Ad R_r^{-1}
\begin{pmatrix}
\alpha&\beta\\\gamma&\delta
\end{pmatrix}=\begin{pmatrix}
\delta&\gamma/r\\r\beta&\alpha
\end{pmatrix}.
\]
\end{ex}
Then $\hF_r=F_rR_r$ frames $\fh_r=\D_{-r}f$ but, using
Exercise~\ref{ex:7}, we have
\[
\mc{\hF_r}=\Ad R_r^{-1}(\mc{F_r})=
\begin{pmatrix}
0&\d f^c\\r\d f&0
\end{pmatrix}
\]
so that $\hF_r$ also frames $\T_r f^c$.  We have therefore proved a
result due to Bianchi \cite{Bia05A} when $n=3$:
\begin{thm}\label{th:13}
$\T_r f^c=\D_{-r}\T_r f$ modulo $\Mob$.
\end{thm}

Similarly, we can compute the interaction of Darboux transforms and
$T$-transforms: let $f:(M,Q)\to\R^n$ be isothermic with \Ch\
transform $f^c$ and let $\fh=\D_r f$.  As usual, frame $(f,\fh)$ with
\[
F_0=
\begin{pmatrix}
\fh\gi&f\\\gi&1
\end{pmatrix}
\]
so that
\[
\mc{F_0}=
\begin{pmatrix}
(\d f)g^{-1}&\d f\\-\gi(\d\fh)\gi&-\gi\d f
\end{pmatrix}=
\begin{pmatrix}
(\d f)g^{-1}&\d f\\-r\d f^c&-\gi\d f
\end{pmatrix}.
\]
Now let $F_s$ frame $(f_s,\fh_s)$ where $F_s$ solves
\[
\mc{F_s}=\mc{F_0}+s
\begin{pmatrix}
0&0\\\d f^c&0
\end{pmatrix}=\begin{pmatrix}
(\d f)g^{-1}&\d f\\(s-r)\d f^c&-\gi\d f
\end{pmatrix}.
\]
Then Theorem~\ref{th:12} tells us that $f_s=\T_s f$ while, from
Theorem~\ref{th:10}, we have $\fh_s=\D_{r-s}f_s$.

On the other hand, set $\hF_s=F_sR_{s-r}$ so that $\hF_s$ frames
$(\fh_s,f_s)$ and $\hF_0$ frames $(\fh,f)$.  Using
Exercise~\ref{ex:7}, we have
\begin{align*}
\mc{\hF_s}&=\Ad R_{s-r}^{-1}(\mc{F_s})=
\begin{pmatrix}
-\gi\d f&\d f^c\\(s-r)\d f&(\d f)\gi
\end{pmatrix}\\
&=\mc{\hF_0}+s
\begin{pmatrix}
0&0\\\d f&0
\end{pmatrix}.
\end{align*}
Thus, by Theorem~\ref{th:12}, $\fh_s=\T_s\fh$ and we conclude, as
have \HJMN\ when $n=3$:
\begin{thm}
$\T_s\D_r f=\D_{r-s}\T_s f$ modulo $\Mob$.
\end{thm}

\subsection{Darboux transforms of generalised $H$-surfaces}
\label{sec:darb-transf-gener}

Recall that a special class of isothermic surfaces is furnished by
the generalised $H$-surfaces.  In view of Proposition~\ref{th:8}, we
may characterise these as surfaces $f$ with a unit normal section $N$
such that, for some constant $H\in\R$, $Hf+N$ is dual to $f$:
\[
\d f\wedge(H\d f+\d N)=0.
\]

Fix such an $f$ and seek Darboux transforms of the same kind.  For
simplicity we choose the polarisation $Q$ so that $f^c=Hf+N$ (when
$n=3$, this amounts to taking $-Q$ to be the Hopf differential of
$f$).  In this case, our Riccati equation has a conserved quantity.
Indeed, if $g$ solves
\[
\d g=rg\d f^cg-\d f,
\]
define $I:M\to\R$ by
\[
I=rHg^2-r\{g,N\}-1
\]
where $\{\,,\,\}$ is the anti-commutator in $\Cl_n$:
$\{g,N\}=-2(g,N)$.  We compute:
\begin{align*}
\d I&= rH\{g,\d g\}-r\{\d g,N\}-r\{g,\d N\}\\
&=rH\{g,rg\d f^cg-\d f\}-r\{rg\d f^cg-\d f,N\}-r\{g,\d f^c-H\d f\}\\
&=rH\{g,rg\d f^cg\}-r\{rg\d f^cg,N\}-r\{g,\d f^c\}
\end{align*}
where we have used $\{\d f,N\}=0$ ($N$ is normal to $f$) and $\d N=\d
f^c-H\d f$.  Rearranging this last equation and exploiting $\{\d
f^c,N\}=0$ yields
\begin{align*}
\d I&=rHg^2\{rg,\d f^c\}-r\{g,N\}\{rg,\d f^c\}-r\{g,\d f^c\}\\
&=I\{rg,\d f^c\}.
\end{align*}
This is a linear differential equation for $I$ and so, in particular,
$I$ vanishes identically if it vanishes at a single point.  We
therefore conclude:
\begin{lem}
If $rHg(o)^2-r\{g(o),N(o)\}=1$ then
\begin{equation}
\label{eq:37}
rHg^2-r\{g,N\}\equiv1.
\end{equation}
\end{lem}

\begin{ex}
For \emph{any} Darboux transform $f+g$ of \emph{any} isothermic
surface $f$, show that $\{g,\d f^c\}$ is a closed $1$-form.
\end{ex}

Now let $g$ satisfy \eqref{eq:37} and contemplate $\hN=-gNg^{-1}$: a
unit normal to $\fh=f+g$.  We know that the \Ch\ transform $\fh^c$ of
$\fh$ is given by
\[
\fh^c=f^c+(rg)^{-1}=Hf+N+r^{-1}g^{-1}.
\]

On the other hand, \eqref{eq:37} tells us that $r^{-1}=Hg^2-\{g,N\}$
and a simple computation gives:
\begin{equation}
\label{eq:38}
\fh^c=H(f+g)-gNg^{-1}=H\fh+\hN.
\end{equation}
Thus $\hN$ is a parallel isoperimetric section for $\fh$ with
$(\hat{\mcv},\hN)=H$ and we have proved yet another theorem which is
due to Bianchi \cite{Bia05} in the classical setting:
\begin{thm}
Let $f$ be a generalised $H$-surface with $(\mcv,N)=H$ and choose
initial condition $v\in\R^n$ and parameter $r\in\R^\times$ so that
$g(o)=v-f(o)$ satisfies
\[
rHg(o)^2-r\{g(o),N(o)\}=1.
\]
Then the Darboux transform $\D^v_r$ is also a generalised $H$-surface
with the same $H$.
\end{thm}
Thus of the $(n+1)$-dimensional family of Darboux transforms of a
generalised $H$-surface, an $n$-dimensional family also produce
generalised $H$-surfaces.

When $H\neq 0$, the conserved quantity \eqref{eq:37} has a simple
geometric interpretation: multiplying by $H$ and completing the
square gives
\[
\frac{H}{r}-1\equiv(Hg-N)^2=\bigl(H\fh-(Hf+N)\bigr)^2
\]
or, equivalently,
\[
\bigl(\fh-(f+\tfrac{1}{H}N)\bigr)^2\equiv\frac{1}{Hr}-\frac{1}{H^2}.
\]
Recall that $f+\tfrac{1}{H}N$ is the parallel generalised $H$-surface
dual to $f$ and conclude that $\fh$ lies on the tube of radius
$\sqrt{1/H^2-1/Hr}$ about this parallel surface.  In particular, we
must have
\[
\frac{1}{Hr}\leq\frac{1}{H^2}.
\]

The extreme case $H=r$ is not without interest: here
$\fh=f+\tfrac{1}{H}N$ so that $\fh$ is simultaneously dual to $f$ and
a Darboux transform of $f$.  In fact, this property characterises
generalised $H$-surfaces\footnote{For CMC surfaces in $\R^3$, this
was known to Bianchi \cite[footnote p.~132]{Bia05}, see also
\cite{HerPed97}.} with $H\neq 0$:
\begin{ex}\label{ex:8}
Let $f:(M,Q)\to\R^n$ be isothermic and $H\in\R^\times$.  Show that
the following are equivalent:
\begin{enumerate}
\item $f$ admits a parallel isoperimetric section $N$ with
$(\mcv,N)=H$.
\item $f$ has a Darboux transform which is also dual to $f$: $\D_r
f=rH^{-2}f^c$.
\item $f$ has a unit normal $N$ such that $N/H$ solves a Riccati
equation of $f$.
\end{enumerate}
\end{ex}

\subsection{Bianchi permutability and the Clifford algebra cross-ratio}
\label{sec:bianchi-perm-cliff}

We begin by stating a permutability theorem for Darboux transforms
that was proved by Bianchi \cite{Bia05} when $n=3$, \HJP\ 
\cite{HerPed97} when $n=4$ and, independently of this writer, Schief
\cite{Sch} in full generality:
\begin{perm}
Let $f:(M,Q)\to\R^n$ be isothermic, $r_1,r_2\in\R^\times$ and
$f_1=\D_{r_1}f$, $f_2=\D_{r_2}f$ distinct Darboux transforms of $f$.
Then there is a fourth isothermic surface $\fh$ such that
\[
\fh=\D_{r_2}f_1=\D_{r_1}f_2.
\]
\end{perm}

In these notes, we shall give two proofs of this result using rather
different ideas.  The first relies on the Clifford algebra
cross-ratio to which we now turn:
\begin{defn}
Let $v_0,v_1,v_2,v_3$ be distinct points in $\R^n$.  The
\emph{Clifford algebra cross-ratio} of these points is given by
\begin{align*}
C(v_0,v_1,v_2,v_3)&=(v_1-v_0)(v_2-v_1)^{-1}(v_2-v_3)(v_3-v_0)^{-1}\\
&=(v_0-v_1)(v_1-v_2)^{-1}(v_2-v_3)(v_3-v_0)^{-1}\in\Cl_n.
\end{align*}
\end{defn}

This cross-ratio is almost invariant under the action of the M\"obius
group:
\begin{ex}
Let $v_0,v_1,v_2,v_3$ be distinct points in $\R^n$ and $T\in\SL$ with
\[
T=
\begin{pmatrix}
a&b\\c&d
\end{pmatrix}.
\]
\begin{enumerate}
\item Show that $T\cdot v_1-T\cdot
v_0=(cv_0+d)^{-t}(v_1-v_0)(cv_0+d)^{-1}$.

\textbf{Hint}: recall that $a^td-c^tb=1$ and that $a^tc\in\R^n$ so
that $a^tc=c^ta$.
\item Write
$C(v_0,v_1,v_2,v_3)=(v_1-v_0)(v_1-v_2)^{-1}(v_3-v_2)(v_3-v_0)^{-1}$
and deduce that
\[
C(T\cdot v_0,T\cdot v_1,T\cdot v_2,T\cdot v_3)=
(cv_0+d)^{-t}C(v_0,v_1,v_2,v_3)(cv_0+d)^{t}.
\]
\end{enumerate}
\end{ex}

In particular, the condition that four points have real cross-ratio
is conformally invariant.  In fact, we can say more: 
\begin{prop}[\cite{Cie97}]\label{th:14}
$C(v_0,v_1,v_2,v_3)=r\in\R$ if and only if $v_0,v_1,v_2,v_3$ lie on a
circle and have real cross-ratio $r$.
\end{prop}
\begin{proof}
Possibly after a M\"obius transformation, we may assume that
$v_0,v_1,v_2,v_3$ lie on a $\R^2\subset\R^n$ so that their
cross-ratio lies in $\Cl_2=\H$.  Write $\H=\C\oplus j\C$.  Then
$\R^2=j\C$ and, writing $v_i=jz_i$, we see that
\[
C(v_0,v_1,v_2,v_3)=jC_\C(v_0,v_1,v_2,v_3)j^{-1}=
\overline{C_\C(v_0,v_1,v_2,v_3)}
\]
where $C_\C$ is the usual complex cross-ratio which is real if and
only if the $z_i$ are concircular\footnote{Indeed, possibly after a
second M\"obius transformation, we may assume $z_0,z_1,z_2$ are real
and then solve for $z_3$:
$z_3=(z_2(z-z_0)+rz_0(z_2-z_1))/(r(z_2-z_1)+(z_1-z_0))\in\R$.}.
\end{proof}

The relevance of the cross-ratio to Bianchi permutability comes from
the following considerations: with $f$ isothermic and
$f_i=\D_{r_i}f=f+g_i$, $i=1,2$, distinct Darboux transforms, suppose
that the theorem is true so that we have
$\fh=\D_{r_2}f_1=\D_{r_1}f_2$ and write
\[
\fh=f_1+g_{12}=f_2+g_{21}.
\]
Now
\[
\d\fh=r_2g_{12}\d f_1^cg_{12}=
\frac{r_2}{r_1}g^{\vphantom{1}}_{12}g_1^{-1}\d f g^{\vphantom{1}}_{12}g_1^{-1}
\]
and, in the same way, we also have
\[
\d\fh=
\frac{r_1}{r_2}g^{\vphantom{1}}_{21}g_2^{-1}\d f g^{\vphantom{1}}_{21}g_2^{-1}.
\]
Equating these, we arrive at
\begin{equation}
\label{eq:39}
\frac{r_2^2}{r_1^2}\d fg_1^{-1}g^{\vphantom{1}}_{12}g_{21}^{-1}g^{\vphantom{1}}_2=
\cross\d f.
\end{equation}
Taking Clifford algebra norms of both sides gives
\begin{equation}
\label{eq:40}
\frac{r_2^2}{r_1^2}=g_1^2g_{12}^{-2}g_{21}^{2}g_2^{-2}
\end{equation}
and \eqref{eq:39} becomes
\begin{equation}
\label{eq:41}
[\cross,\d f]=0.
\end{equation}
Now recall that if $f, f+g$ envelope a $2$-sphere congruence,
$N\mapsto -gNg^{-1}$ is a parallel isomorphism of normal bundles.  In
the present setting, we therefore arrive at two such isomorphisms
between the normal bundles of $f$ and $\fh$ and we make the
\emph{ansatz} that these coincide\footnote{When $n=3$, this amounts to
choosing a sign.}: that is, for $N$ normal to $f$, we assume,
\[
g^{\vphantom{1}}_{21}g^{\vphantom{1}}_2Ng_2^{-1}g_{21}^{-1}=
g^{\vphantom{1}}_{12}g^{\vphantom{1}}_1Ng_1^{-1}g_{12}^{-1}.
\]
Rearranging this and multiplying by $g_1^2g_2^{-2}$ gives us
\[
[\cross,N]=0
\]
which taken together with \eqref{eq:41} tells us that $\cross$
commutes with all of $\R^n$ and so is central in $\Cl_n$.  Moreover,
using
\[
g_1+g_{12}=\fh-f=g_2+g_{21}
\]
one checks that $\cross\in\R^n\cdot\R^n\subset\Cl_n$.  However, when
$n>2$, $\R^n\cdot\R^n$ intersects the centre of $\Cl_n$ in $\R$ alone
so we conclude that $\cross\in\R$ and, in view of \eqref{eq:40}, we
must have
\[
C(f,f_1,\fh,f_2)=\cross=\pm\frac{r_2}{r_1}.
\]
To fix the sign, we consider the degenerate case where $r_1=r_2$
where, according to Exercise~\ref{ex:6}, we may take $\fh=f$ and then
the cross-ratio is $1=r_2/r_1$.  We therefore conclude that we should
have
\begin{equation}
\label{eq:42}
C(f,f_1,\fh,f_2)=\frac{r_2}{r_1}
\end{equation}
\begin{rem}
In the case $n=4$, \HJP\ arrive at the quaternionic version of the
same \emph{ansatz} by considerations coming from the theory of
discrete isothermic surfaces \cite{BobPin96,HerHofPin99}. 
\end{rem}

\begin{ex}
\begin{enumerate}
\item If \eqref{eq:42} holds, show that
\begin{equation}
\label{eq:43}
g^{\vphantom{1}}_{12}=
(g^{\vphantom{1}}_1-g^{\vphantom{1}}_2)
r_1g_2^{-1}(r_2g_1^{-1}-r_1g_2^{-1})^{-1}.
\end{equation}
\item Deduce from~\eqref{eq:43} that
\begin{equation}
\label{eq:44}
\fh=(r_2f^{\vphantom{1}}_1g_1^{-1}-r_1f^{\vphantom{1}}_2g_2^{-1})
(r_2g_1^{-1}-r_1g_2^{-1})^{-1}.
\end{equation}
\end{enumerate}
\end{ex}

To prove our theorem, it remains to show that if $\fh$ is defined by
\eqref{eq:44} then we really do have $\fh=\D_{r_2}f_1=\D_{r_1}f_2$.
To show the first of these amounts to proving that
\[
\d g_{12}=r_2g_{12}\d f_1^c g_{12}-\d f_1,
\]
that is,
\[
\d g_{12}=\frac{r_2}{r_1}g_{12}g_1^{-1}\d fg_1^{-1} g_{12}-r_1g_1\d
f^cg_1.
\]
This is a tedious but straightforward verification using
\eqref{eq:43}.
\begin{ex}
Check the grisly details!
\end{ex}

This completes the proof of the permutability theorem and gives us
more.  In fact, we have shown (as has Schief \cite{Sch}):
\begin{thm}\label{th:15}
Let $f$ be isothermic with distinct Darboux transforms
$f_1=\D_{r_1}f$ and $f_2=\D_{r_2}f$.  Then there is a fourth surface
$\fh=\D_{r_2}f_1=\D_{r_1}f_2$ such that corresponding points on
$f,f_1,\fh,f_2$ are concircular with real cross-ratio $r_2/r_1$.
\end{thm}

We call four surfaces in the configuration of Theorem~\ref{th:15} a
\emph{Bianchi quadrilateral}.

Our explicit formula for the fourth surface of a \Bq\ allows us to
give \emph{algebraic} proofs\footnote{All the material in the
remainder of this section resulted from conversations with Udo
Hertrich-Jeromin.} of several results of Bianchi \cite{Bia05}
concerning the geometry of such configurations which immediately
extend to our $n$-dimensional setting.

First, let us consider the \Ch\ transform of a \Bq: let
$(f,f_1,\fh,f_2)$ be such a quadrilateral and contemplate the \Ch\ 
transforms $f^c$, $f_1^c=\D_{r_1}f^c=f^c+(r_1g_1)^{-1}$,
$f_2^c=\D_{r_2}f^c=f^c+(r_2g_2)^{-1}$.  We now have three rival \Ch\ 
transforms of $\fh$: $f_1^c+(r_2g_{12})^{-1}$,
$f_2^c+(r_1g_{21})^{-1}$ and $\fh^c$ given by the permutability
theorem so as to make $(f^c,f_1^c,\fh^c,f_2^c)$ a \Bq\footnote{That
$\fh^c$ is also a \Ch\ transform follows from Theorem~\ref{th:9}.}.
Of course, these three possibilities can only differ by constants
but, in fact, they coincide exactly:
\begin{ex}
If $g_1,g_2,g_{12}\in\R^n$ are given by \eqref{eq:43} then
$g_1^c=(r_1g_1)^{-1}$, $g_2^c=(r_2g_2)^{-1}$,
$g_{12}^c=(r_2g_{12})^{-1}$ also satisfy \eqref{eq:43}:
\[
g^c_{12}=(g^c_1-g^c_2)r_1(g^c_2)^{-1}(r_2(g^c_1)^{-1}-r_1(g^c_2)^{-1})^{-1}.
\]
\end{ex}
Thus $\fh^c=f_1^c+(r_2g_{12})^{-1}$ and, by symmetry,
$\fh^c=f_2^c+(r_1g_{21})^{-1}$.  To summarise:
\begin{thm}
\label{th:16}
The \Ch\ transform of a \Bq\ is also a \Bq.
\end{thm}

A similar but slightly more elaborate analysis shows that a Darboux
transform of a \Bq\ is another \Bq.  For this we need a version of
the hexahedron lemma of \cite{HerHofPin99}:
\begin{lem}
\label{th:17}
Let $v,v_1,\vh,v_2$ be distinct concircular points in $\R^n$ with
Clifford algebra cross-ratio $C(v,v_1,\vh,v_2)=r_2/r_1$ and let
$v'\in\R^n$ distinct from $v,v_1,v_2$.  Then, for $r_3\in\R^\times$,
there are unique points $v'_1,\wh,v'_2$ such that
\begin{align*}
C(v,v_1,\vh,v_2)&=C(v',v'_1,\wh,v'_2)=r_2/r_1\\
C(v,v',v'_1,v_1)&=C(v_2,v'_2,\wh,\vh)=r_1/r_3\\
C(v,v',v'_2,v_2)&=C(v_1,v'_1,\wh,\vh)=r_2/r_3.
\end{align*}
Moreover, all 8 points lies on a single $2$-sphere or plane in
$\R^n$.
\end{lem}
\begin{proof}
The points $v,v_1,\vh,v_2,v'$ lie on a $2$-sphere or plane and so,
after a M\"obius transformation, may be taken to lie on a copy of
$\R^2$ where, as in the proof of Proposition~\ref{th:14}, all Clifford
algebra cross-ratios reduce to complex cross-ratios.  One now solves
\begin{align*}
C_\C(v,v',v'_1,v_1)&=r_1/r_3\\
C_\C(v,v',v'_2,v_2)&=r_2/r_3\\
C_\C(v',v'_1,\wh,v'_2)&=r_2/r_1
\end{align*}
to obtain, in turn, $v'_1,v'_2,\wh\in\C$ and then checks that the
remaining two equations hold: a task best left to a computer algebra
engine (c.f. \cite{HerHofPin99}).
\end{proof}

Now suppose that we start with a \Bq\ $(f,f_1,\fh,f_2)$ with
$f_1=\D_{r_1}f$, $f_2=\D_{r_2}f$ and take a third Darboux transform
$f'=\D_{r_3}f$ of $f$.  The permutability theorem yields isothermic
surfaces
\begin{align*}
f'_1&=\D_{r_3}f_1=\D_{r_1}f'\\
f'_2&=\D_{r_3}f_2=\D_{r_2}f'\\
\intertext{and, finally, thanks to Lemma~\ref{th:17},}
\fh'&=\D_{r_3}f'=\D_{r_1}f_2'=\D_{r_2}f'_1.
\end{align*}
Thus these 8 surfaces form the vertices of a cube all of whose faces
are \Bq s!  In particular, we have:
\begin{thm}
\label{th:18}
For suitably chosen initial conditions, the Darboux transform of a
\Bq\ is a \Bq.
\end{thm}

As a last application of these ideas, let us show that if the first
three surfaces in a \Bq\ are generalised $H$-surfaces with the same
$H\neq 0$ then so is the fourth.  We begin by examining a degenerate
case: so let $f$ be a generalised $H$-surface with $H\neq 0$ and
$f^c=Hf+N$.  We have seen that the parallel surface $f^c/H$ is a Darboux
transform of $f$: $f^c/H=\D_H f$.  Now take a second Darboux
transform $f_1=\D_{r_1}f$ which is also a generalised $H$-surface and
contemplate the \Bq\ $(f,f_1,\fh,f^c/H)$.
\begin{prop}
\label{th:19}
$\fh=f_1^c/H$.
\end{prop}
\begin{proof}
We must check that $C(f,f_1,f_1^c/H,f^c/H)=H/r_1$.  However, from
\eqref{eq:37}, we know that
\[
f_1^c=f^c+(r_1g_1)^{-1}=Hf_1-g_1Ng_1^{-1}
\]
whence
\begin{align*}
g_{12}&=f_1^c/H-f_1=-g_1Ng_1^{-1}\\
g_{21}&=f_1^c/H-f^c/H=(r_1g_1)^{-1}.
\end{align*}
Finally, $g_2=N/H$ so that
\[
C(f,f_1,f_1^c/H,f^c/H)=-g_1(g_1Ng_1^{-1})^{-1}(r_1g_1)^{-1}(N/H)^{-1}=H/r_1
\]
since $N^2=-1$.
\end{proof}

Thus a Darboux pair of generalised $H$-surfaces, together with their
parallel $H$-surfaces form a \Bq.

We are now in a position to prove:
\begin{thm}
Let $(f,f_1,\fh,f_2)$ be a \Bq\ with $f,f_1,f_2$ generalised
$H$-surfaces with the same $H\neq 0$.  Then $\fh$ is also a
generalised $H$ surface.
\end{thm}
\begin{proof}
Consider the configuration of 8 surfaces obtained from
Lemma~\ref{th:17} starting with $(f,f_1,\fh,f_2)$ and $f'=f^c/H$.
Proposition \ref{th:19} tells us that $f_1'=f_1^c$ and $f'_2=f_2^c$
while, from Theorem~\ref{th:16}, we have
\[
C(f^c,f_1^c,\fh^c,f_2^c)=r_2/r_1.
\]
Now, an obvious scaling symmetry of the cross-ratio gives
\[
C(f^c,f_1^c,\fh^c,f_2^c)=C(f^c/H,f_1^c/H,\fh^c/H,f_2^c/H)
\]
so that
\[
C(f^c/H,f_1^c/H,\fh^c/H,f_2^c/H)=r_2/r_1=C(f^c/H,f_1^c/H,\fh',f_2^c/H).
\]
We conclude that $\fh'=\fh^c/H$, that is, $\fh^c/H=\D_H\fh$ so that,
by Exercise~\ref{ex:8}, $\fh$ is a generalised $H$-surface also.
\end{proof}

\subsection{Isothermic surfaces via the vector Calapso equation }
\label{sec:isoth-surf-via}

Let us pause from our main development and digress\footnote{This
section may be omitted from a first reading.} to consider the
approach of Calapso \cite{Cal03,Cal15} to isothermic surfaces.  For
$n=3$, he reduced the problem to the study of a fourth-order
non-linear partial differential equation for a function that turns
out to be (the coefficient of) the conformal Hopf differential in CCL
coordinates.  This PDE is equivalent to the stationary version of the
second flow of the Davey--Stewartson II hierarchy \cite{Fer97} ---a
hierarchy of integrable PDE with mysterious\footnote{\textit{Note
added in December 2001}: these matters are now a little less
mysterious to me, see \cite{BurPedPin01}} (to this author)
connections to conformal geometry \cite{Kon00,KonLan00}.

In this section, we describe a simple generalisation of Calapso's
approach which treats isothermic surfaces in $\R^n$ and was also
arrived at independently by Schief \cite{Sch}.  For this we adapt
an argument of \cite{BurHerPed97} and so temporarily abandon our
Clifford algebra formalism to work with frames in $\Op$.

Fix a basis $e_0,\dots,e_{n+1}$ of $\lor$ with $e_1,\dots,e_n$
space-like orthogonal and $e_0,e_{n+1}\in\L^+$ with
$(e_0,e_{n+1})=-\half$.  A map $F:M\to\Op$ frames an immersion
$\sf:M\to\PL$ if $\pi(Fe_0)=\sf$, that is,
\[
Fe_0\in\sf.
\]
Let $\sf:M\to\PL$ be isothermic and fix $z=x+iy$ a CCL coordinate.
We are going to construct an essentially unique and M\"obius
invariant frame for $\sf$.  Firstly, choose $f:M\to\L^+$ to be the
(unique) lift of $\sf$ with
\[
(\d f,\d f)=\d x^2+\d y^2
\]
and set $X=f_x$, $Y=f_y$: these are orthonormal and space-like.  Now
contemplate the conformal Gauss map of $\sf$
(cf~page~\pageref{page:conf-gauss}):
\[
Z_{\sf}=\<f,f_x,f_y,f_{xx}+f_{yy}\>^\perp
\]
which is isomorphic to the normal bundle $\Nf$ and so a flat bundle
with respect to its induced connection.  Choose orthonormal parallel
sections $N_1,\dots,N_{n-2}$ of $Z_{\sf}$.  Finally, let
$\fh:M\to\L^+$ be (uniquely) determined by the demands that $\fh$ is
orthogonal to $X,Y,N_1,\dots,N_{n-2}$ and that $(f,\fh)=-\half$.

This data defines a frame $F:M\to\Op$ of $\sf$ such that
\begin{gather*}
Fe_0=f\\
Fe_1=X,\quad Fe_2=Y\\
Fe_i=N_{i-2}\quad\text{for $3\leq i\leq n$}\\
Fe_{n+1}=\fh
\end{gather*}
which is completely determined by $\sf$ and $z$ up to the right
action of $\O[n-2]$ permuting the choice of parallel framing of
$Z_{\sf}$.

Each $N_i$ is parallel so that $\d N_i\in\<f,f_x,f_y\>$.  Moreover,
$x,y$ are curvature line coordinates so there are functions
$\kappa^{(1)}_i,\kappa^{(2)}_i$ such that
\[
\d N_i=-\kappa^{(1)}_i f_x\d x-\kappa^{(2)}_i f_y\d y+\tau_i f
\]
for some $1$-form $\tau_i$.  Now $(N_i,f_{xx}+f_{yy})=0$ while
\begin{align*}
\kappa^{(1)}_i&=-(N_{i,x},f_x)=(N_i,f_{xx})\\
\kappa^{(2)}_i&=(N_i,f_{yy})
\end{align*}
so that
\[
\kappa^{(1)}_i+\kappa^{(2)}_i=0.
\]
We therefore set $\kappa_i=\kappa^{(1)}_i$ and conclude
\begin{equation}
\label{eq:45}
\d N_i=-\kappa_i X\d x+\kappa_i Y\d y+\tau_i f.
\end{equation}
The $\kappa_i$ are the components of the conformal Hopf differential
with respect to the frame $N_1,\dots,N_{n-2}$ of $Z_{\sf}$ and our
CCL coordinate $z=x+iy$:
\begin{ex}
Recall the definition of the conformal Hopf differential from
page~\pageref{page:conf-Hopf}.  Show that
\[
K_{\sf}(N_i+\sf)=\kappa_i
\]
\end{ex}
\begin{rem}
If, instead of the isometric lift, we take a Euclidean lift
$f':M\to E_{v_\infty}\subset\L^+$, we can use the Euclidean normal
bundle and parallel sections $N_1',\dots,N_{n-2}'$ to compute
$K_{\sf}$.  We then get
\[
K_{\sf}(N'_i+\sf)=\frac{e^{u}}{2}(\kappa_i'-\kappa_i'')
\]
where $(\d f',\d f')=e^{2u}(\d x^2+\d y^2)$ and the
$\kappa_i',\kappa_i''$ are the Euclidean principal curvatures for
$N'_i$.  This gives the formulation of Calapso \cite{Cal03} and
Schief \cite{Sch}.
\end{rem}

Returning to our frame, we note that
\[
\d X,\d Y\perp\<X,Y\>
\]
since $X,Y$ are an orthonormal coordinate frame for a flat metric on
$M$ and, taking this together with \eqref{eq:45}, we compute the \MC\
form of $F$:
\[
B=\mc{F}=
\begin{pmatrix}
&\chi_1&\chi_2&\tau&\\
\d x&&&-\kappa\d x&-\chi_1\\
\d y&&&\kappa\d y&-\chi_2\\
&\kappa^T\d x&-\kappa^T\d y&&-\tau^T\\
&-\d x&-\d y&&
\end{pmatrix}
\]
where $\kappa=(\kappa_1,\dots,\kappa_{n-2})$,
$\tau=(\tau_1,\dots,\tau_{n-2})$ and $\chi_1,\chi_2$ are two more
$1$-forms.

Now $B$ satisfies the \MC\ equations.  Conversely, any $B$ of the
above form that satisfies the \MC\ equations can be locally
integrated to give $F:M\to\Op$ with $B=\mc{F}$.  If we then define
$f=Fe_0$, $N_i=Fe_{i+2}$, $1\leq i\leq n-2$, we see that
\[
f_x=Fe_1\qquad f_y=Fe_2
\]
so that the $N_i$ are normal to $f$.  Moreover, we have
\[
\d N_i=-\kappa_i f_x\d x+\kappa_i f_y\d y+\tau_i f
\]
which shows that $x,y$ are CCL coordinates so that $\sf$ is
isothermic and, in addition, that the $N_i$ are a parallel frame for
the conformal Gauss map of $\sf$.

So let us examine the \MC\ equations of $B$: these amount to
\begin{subequations}
\begin{gather}
\chi_1\wedge\d x+\chi_2\wedge\d y=0\label{eq:46}\\
\chi_2\wedge\d x-\chi_1\wedge\d y+(\kappa,\kappa)\d y\wedge\d
x=0\label{eq:47}\\
\d(\kappa\d x)+\tau\wedge\d x=0\label{eq:48}\\
\d(\kappa\d y)-\tau\wedge\d y=0\label{eq:49}\\
\d\tau-\chi_1\wedge\kappa\d x+\chi_2\wedge\kappa\d
y=0\label{eq:50}\\
\d\chi_1+\tau\wedge\kappa\d x=0\label{eq:51}\\
\d\chi_2-\tau\wedge\kappa\d y\label{eq:52}
\end{gather}
\end{subequations}
where we have written $(\kappa,\kappa)$ for
$\sum_{i=1}^{n-2}\kappa_i^2$.

Write
\[
\chi_i=\chi_{i1}\d x+\chi_{i2}\d y.
\]
Then \eqref{eq:46} is equivalent to $\chi_{12}=\chi_{21}$ and we
denote this common value by $\psi$.

Similarly, \eqref{eq:47} is equivalent to
\begin{equation}
\label{eq:53}
\chi_{11}+\chi_{22}=-(\kappa,\kappa)
\end{equation}
so we write
\begin{equation}
\label{eq:54}
\chi_{11}=\half\bigl(u-(\kappa,\kappa)\bigr),\quad
\chi_{22}=\half\bigl(-u-(\kappa,\kappa)\bigr)
\end{equation}
for some function $u:M\to\R$.

The vector valued equations \eqref{eq:48} and \eqref{eq:49} amount to
\begin{equation}
\label{eq:55}
\tau=\kappa_x\d x-\kappa_y\d y
\end{equation}
while \eqref{eq:50} gives
\[
\d\tau=2\psi\kappa\d y\wedge\d x
\]
or, using \eqref{eq:55},
\[
\kappa_{xy}=\psi\kappa.
\]

Finally, \eqref{eq:51} and \eqref{eq:52} give
\begin{subequations}
\label{eq:56}
\begin{align}
\half u_y&=\psi_x+(\kappa,\kappa)_y\\
\half u_x&=-\psi_y-(\kappa,\kappa)_x.
\end{align}
\end{subequations}

Now $\d u=0$ which is the same as
\[
\Delta\psi+2(\kappa,\kappa)_{xy}=0.
\]
Thus the \MC\ equations for $B$ boil down to the \emph{vector Calapso
equation}:
\begin{subequations}
\label{eq:57}
\begin{gather}
\kappa_{xy}=\psi\kappa\\
\Delta\psi+2(\kappa,\kappa)_{xy}=0.
\end{gather}
\end{subequations}

\begin{rem}
When $n=3$, $\kappa$ is scalar and we can eliminate $\psi$ to obtain
Calapso's original equation\footnote{In fact, this equation first
appeared in the thesis of Rothe \cite{Rot97}.}:
\[
\Delta\biggl(\frac{\kappa_{xy}}{\kappa}\biggr)+2(\kappa^2)_{xy}=0.
\]
\end{rem}

Conversely, given a solution $\kappa,\psi$ of the vector Calapso
equation \eqref{eq:57}, we integrate \eqref{eq:56} to obtain $u$,
define $\tau$ by \eqref{eq:55} and finally $\chi_i$ by \eqref{eq:54}
together with $\chi_{12}=\chi_{21}=\psi$ to get a \MC\ solution and
so a frame of an isothermic surface, unique up to a M\"obius
transformation.

In fact, we get more from this analysis: there is a constant of
integration in the definition of $u$.  Replacing $u$ by $u+r$ gives
us a new \MC\ solution
\[
B_{r/2}=B+\frac{r}{2}
\begin{pmatrix}
&\d x&-\d y&&\\&&&&-\d x\\&&&&\d y\\&&&&\\&&&&
\end{pmatrix}
\]
and so a new isothermic surface $\sf_{r/2}$.

We have seen this before.  In our Clifford algebra formulation,
\[
B=
\begin{pmatrix}
*&e_1\d x+e_2\d y\\e_1\chi_1+e_2\chi_2&*
\end{pmatrix}
\]
and
\[
B_{r/2}=B+\frac{r}{2}
\begin{pmatrix}
0&0\\e_1\d x-e_2\d y&0
\end{pmatrix}.
\]
One easily checks that
\[
(e_1\d x+e_2\d y)^c=e_1\d x-e_2\d y
\]
so that, by Theorem~\ref{th:12}, $\sf_{r/2}$ is the $T$-transform
$\T_{r/2}\sf$ of $\sf$.

To summarise: \emph{each solution of the vector Calapso equation
\eqref{eq:57} gives rise to the $1$-parameter family of
$T$-transforms of an isothermic surface and conversely.}

\section{Curved flats}
\label{sec:curved-flats}

The rich transformation theory of isothermic surfaces strongly
suggests the presence of an underlying integrable system.  This is
indeed the case: the integrable system in question is that of
\emph{curved flats} discovered by Ferus--Pedit \cite{FerPed96} which
is very closely related to the ``$n$-dimensional system'' of Terng
\cite{Ter97}.

It is a main result of \cite{BurHerPed97} that Darboux pairs in
$\R^3$ amount to curved flats in a certain Grassmannian.  In this
section, we shall show that such a result holds in arbitrary
codimension and, in so doing, unify much of the transformation theory
of Section~\ref{sec:isoth-surf-class}.

\subsection{Curved flats in symmetric spaces}
\label{sec:curv-flats-symm}

Let $G/K$ be a symmetric space.  Thus $G$ is a Lie group (usually, for
us, semisimple) with an involution $\tau:G\to G$ and $K$ is a closed
subgroup open in the fixed set of $\tau$.  The derivative at $1$ of
$\tau$ is an involution, also called $\tau$, of the Lie algebra $\g$
of $G$.  We have a decomposition
\begin{equation}
\label{eq:58}
\g=\k\oplus\p
\end{equation}
into $\pm1$-eigenspaces of $\tau$.  The $+1$-eigenspace $\k$ is the
Lie algebra of $K$ and, since $\tau$ is an involution of $\g$, we
have:
\begin{equation}
\label{eq:59}
[\k,\k]\subset\k,\qquad[\k,\p]\subset\p,\qquad[\p,\p]\subset\k.
\end{equation}

The left action of $G$ on $G/K$ differentiates to give a surjection
$\g\to T_{gK}G/K$:
\[
\xi\mapsto\dt(\exp t\xi)gK
\]
with kernel $\Ad(g)\k$ which therefore restricts to give an
isomorphism $\Ad(g)\p\cong T_{gk}G/K$.  In this way, we view each
tangent space to $G/K$ as a subspace of $\g$.

\begin{defn}[\cite{FerPed96}]
An immersion $\phi:M\to G/K$ of a manifold $M$ is a \emph{curved
flat} if each $\d\phi(T_p M)$ is an abelian subalgebra of $\g$ (where
$T_{\phi(p)}G/K\subset\g$ as above).
\end{defn}
Under mild conditions on $G$, this amounts to the demand that the
curvature operator of the canonical connection of $G/K$ vanishes on
each $\bigwedge^2\d\phi(T_pM)$.

A \emph{frame} of $\phi$ is a map $F:M\to G$ which is mapped onto
$\phi$ by the coset projection $G\to G/K$:
\[
\phi=FK.
\]
Since the coset projection is locally trivial, frames exist locally
and if $F$ is one such, any other is of the form $Fk$ with $k:M\to
K$.

As we have already seen, a map $F:M\to G$ is determined by its \MC\
form $A=\mc{F}\in\Omega^1\otimes\g$ which satisfies the \MC\
equations:
\begin{equation}
\label{eq:60}
\d A+\half[A\wedge A]=0
\end{equation}
where
\[
[A\wedge B]_{X,Y}=[A_X,B_Y]-[A_Y,B_X].
\]

Conversely, if $A\in\Omega^1\otimes\g$ solves \eqref{eq:60} then we
can locally integrate to find $F:M\to G$, unique up to left
multiplication by constants, with $A=\mc{F}$.

For $F$ a frame of $\phi$ and $A=\mc{F}$, write
\[
A=A_\k+A_\p
\]
according to the decomposition \eqref{eq:58}.  Viewing $\d\phi$ as a
$\g$-valued $1$-form, we have
\[
\d\phi=\Ad(F)A_\p
\]
so that $\phi$ is a curved flat if and only if
\[
[A_\p\wedge A_\p]=0.
\]

Now the \MC\ equations \eqref{eq:60} decompose into their components
in $\k$ and $\p$ which, in view of \eqref{eq:59}, read
\begin{gather*}
\d A_k+\half[A_\k,A_\k]+\half[A_\p\wedge A_\p]=0\\
\d A_\p+[A_\k\wedge A_\p]=0
\end{gather*}
so that $\phi$ is a curved flat if and only if these equations
decouple further to give:
\begin{subequations}
\label{eq:61}
\begin{gather}
\label{eq:62}\d A_k+\half[A_\k,A_\k]=0\\
\d A_\p+[A_\k\wedge A_\p]=0\\
[A_\p\wedge A_\p]=0
\end{gather}
\end{subequations}
Now observe that \eqref{eq:61} are the coefficients of a
\emph{spectral parameter} $\lambda\in\R$ in the \MC\ equations for
the pencil of $1$-forms $A_\lambda\in\Omega^1\otimes\g$ given by
\[
A_\lambda=A_\k+\lambda A_\p.
\]
That is,
\begin{prop}
Let $F:M\to G$ with $\mc{F}=A_\k+A_\p$.  Then $F$ frames a curved
flat if and only if $A_\lambda=A_\k+\lambda A_\p$ satisfies
\[
\d A_\lambda+\half[A_\lambda\wedge A_\lambda]=0
\]
for all $\lambda\in\R$.
\end{prop}

We have therefore arrived at a \emph{zero curvature formulation} of
the curved flat condition.

As an immediate consequence, we see that curved flats come in
$1$-parameter families: for each $\lambda\in\R$, we can locally
integrate to find $F_\lambda: M\to G$ with $\mc{F_\lambda}=A_\lambda$
and, since each $(A_\lambda)_\p=\lambda A_\p$, we have
\[
[(A_\lambda)_\p\wedge (A_\lambda)_\p]=0
\]
so that, when $\lambda\neq 0$, $F_\lambda$ frames a new curved flat
$\phi_\lambda:M\to G/K$.  Moreover, this construction is
independent of our original choice of frame $F$:
\begin{ex}
If $F$ and $\hF=Fk$ are two frames of a curved flat $\phi$ then
$\hF_\lambda=F_\lambda k$.
\end{ex}
In fact, the only ambiguity in our construction comes from the
possibility of left multiplying each $F_\lambda$ by a constant
$c_\lambda\in G$.  Thus, the curved flats $\phi_\lambda$ are defined
up to the action of $G$ on $G/K$.

Note that since $A_1=A$, we may take $F_1=F$ and so $\phi_1=\phi$.
Similarly, since $A_0$ is $\k$-valued, $F_0$ may be taken to be
$K$-valued so that $\phi_0$ is constant.

To summarise: 
\begin{thm}\label{th:20}
Let $\phi:M\to G/K$ be a curved flat with $M$ simply connected.
Then, for each $\lambda\in\R$, there is a map $\phi_\lambda:M\to
G/K$, uniquely determined up to the action of $G$, such that
\begin{enumerate}
\item For $\lambda\in\R^\times$, $\phi_\lambda$ is a curved flat;
\item $\phi_1=\phi$;
\item $\phi_0$ is constant.
\end{enumerate}
\end{thm}
We say that the $\phi_\lambda$ comprise the \emph{associated family}
of $\phi$.

So far, our discussion requires no special choice of frame.  However,
special choices are available and useful:  if $F$ frames a curved
flat $\phi$ then \eqref{eq:62} says that $A_\k$ solves the \MC\
equations so that there is a map $k:M\to K$ with $\mc{k}=A_\k$.  We
now have a new frame $\hF=Fk^{-1}$ of $\phi$ with
\[
\mc{\hF}=\Ad k(A-\mc{k})=\Ad(k)A_\p\in\Omega^1\otimes\p.
\]
This prompts:
\begin{defn}
A \emph{flat frame} of a curved flat is a frame $F$ with
$\mc{F}\in\Omega^1\otimes\p$.
\end{defn}

Note that if $F$ is a flat frame then so is each of the $F_\lambda$,
$\lambda\neq 0$:
\[
\mc{F_\lambda}=\lambda\mc{F},
\]
while $F_0$ is constant.

So let $F$ be a flat frame of a curved flat with $\mc{F}=A_\p$.  The
\MC\ equations \eqref{eq:60} read
\begin{gather*}
\d A_\p=0\\
[A_\p\wedge A_\p]=0.
\end{gather*}
We can therefore integrate to get a function $\psi:M\to\p$ with
$\d\psi=A_\p$ and thus
\begin{equation}
\label{eq:63}
[\d\psi\wedge\d\psi]=0.
\end{equation}
\begin{defn}
An immersion $\psi:M\to\p$ is $\p$-flat if it satisfies~\eqref{eq:63}.
\end{defn}

Thus any flat frame gives rise to a $\p$-flat map and, conversely, a
$\p$-flat map $\psi:M\to\p$ gives rise to a $1$-parameter family of
flat frames $F_\lambda$ framing an associated family of curved flats
by solving
\[
\mc{F_\lambda}=\lambda\d\psi
\]
for $\lambda\in\R^\times$.

While we will mostly work with flat frames, we remark that there is
another canonical choice of frame for curved flats.  For this, we
must assume that all $\d\phi(T_pM)$ are conjugate to a fixed
semisimple abelian subalgebra $\a\subset\p$ (this is certainly the
case when each $\d\phi(T_pM)$ is \emph{maximal} abelian and $G/K$ is
a Riemannian symmetric space of semisimple type\footnote{Thus $G$ is
semisimple and $K$ is compact.}).  In this case, one can find a frame
for which each $A_\p(T_pM)=\a$ and then one can prove:
\begin{enumerate}
\item $\d A_\p=0$ so that, for any basis $H_1,\dots,H_l$ of $\a$,
there are coordinates $x_1,\dots,x_l$ on $M$ such that $A_\p=\sum_i
H_i\d x_i$;
\item There is a unique function $u:M\to[\a,\k]\subset\p$ such that
\[
A_\k=[A_\p,u].
\]
\end{enumerate}

The \MC\ equations for $A$ reduce to a differential equation for $u$
called the \emph{$l$-dimensional system associated to $G/K$}
\cite{Ter97}.  This frame is the basis of the approach to curved
flats adopted by Terng and her collaborators
\cite{BruDuPar00,Ter97,TerUhl98,TerUhl00}.

\newpage
\subsection{Curved flats and isothermic surfaces}
\label{sec:curv-flats-isoth}

\subsubsection{The symmetric space $S^n\times S^n\setminus\Delta$}
\label{sec:symm-spac-sntim}

Denote by $Z$ the space $S^n\times S^n\setminus\Delta$ of pairs of
distinct points of $S^n=\R^n\cup\{\infty\}$.  There is a diagonal
action of $\Op$ (and so $\SL$) on $Z$:
\[
g(x,y)=(g\cdot x,g\cdot y).
\]
\begin{ex}
Show that this action is transitive.
\end{ex}

Let $K\subset\SL$ be the stabiliser of $(0,\infty)\in Z$.  From
\eqref{eq:11} we see that $K$ is precisely the subgroup of diagonal
matrices in $\SL$:
\[
K=\left\{
\begin{pmatrix}
a&0\\0&a^{-t}
\end{pmatrix}:a\in\Gamma_n\right\}
\]
which is the fixed set of the automorphism $\tau$ of $\SL$ given by
conjugation by $\displaystyle
\begin{pmatrix}
1&0\\0&-1
\end{pmatrix}\in\Pin$.  The corresponding decomposition $\o=\k+\p$ is
the familiar decomposition into diagonal and off-diagonal matrices:
\[
\k=\left\{
\begin{pmatrix}
\xi&0\\0&-\xi^t
\end{pmatrix}:\xi\in[\R^n,\R^n]\oplus\R
\right\}
\qquad
\p=\left\{
\begin{pmatrix}
0&x\\y&0
\end{pmatrix}: x,y\in\R^n
\right\}.
\]
Finally, $gK\mapsto(g\cdot0,g\cdot\infty)$ is a diffeomorphism so
that $Z$ is identified with the symmetric space $\SL/K$.

\begin{rem}
There is another model for the symmetric space $Z$: it can be viewed
as the Grassmannian of oriented $(1,1)$-planes in $\lor$.  Indeed,
any pair of distinct points in $\PL$ span such a plane while any such
plane contains a unique pair of light-lines which are ordered via the
orientation.
\end{rem}

\subsubsection{Curved flats are Darboux pairs}
\label{sec:curved-flats-are}

A map $\phi:M\to Z=S^n\times S^n\setminus\Delta$ is the same as a
pair of maps $f,\fh:M\to S^n$ whose values never coincide.  Use the
identification of $Z$ with $\SL/K$ to view $\phi$ as a map into
$\SL/K$ and let $F:M\to\SL$ be a frame of $\phi$.  Then
\[
(f,\fh)=(F\cdot0,F\cdot\infty)
\]
so that $F$ frames the pair $(f,\fh)$ in the sense of
Section~\ref{sec:moving-frames}.  Now let
\[
A=\mc{F}=
\begin{pmatrix}
\alpha&\beta\\\gamma&\delta
\end{pmatrix}
\]
so that
\[
A_\p=
\begin{pmatrix}
0&\beta\\\gamma&0
\end{pmatrix}
\]
for $\beta,\gamma\in\Omega^1\otimes\R^n$.  The curved flat condition
$[A_\p\wedge A_\p]=0$ amounts to $\beta\wedge\gamma=0$ which, as long
as $f,\fh$ are immersions, is precisely the condition of
Theorem~\ref{th:10} that $(f,\fh)$ be a Darboux pair of isothermic
surfaces\footnote{Lemma~\ref{th:4} tells us that with
$\beta\wedge\gamma=0$, $\rank\beta=2$.  But $\rank\beta=\rank\d f$ so
$\dim M=2$.}.

Say that a map $(f,\fh):M\to Z$ is \emph{non-degenerate} if both $f$
and $\fh$ are immersions and conclude:
\begin{thm}\label{th:21}
A non-degenerate map $(f,\fh):M\to Z$ is a curved flat if and only if
$(f,\fh)$ is a Darboux pair of isothermic surfaces.
\end{thm}

\subsubsection{Spectral deformation is $T$-transform}
\label{sec:spectr-deform-t}

Given a Darboux pair $\phi=(f,\fh)$, Theorems~\ref{th:21} and \ref{th:20}
provide us with the $1$-parameter associated family
$\phi_\lambda=(f_{(\lambda)},\fh_{(\lambda)})$ of such with
$(f_{(1)},\fh_{(1)})=(f,\fh)$.  In fact, these new isothermic surfaces
are $T$-transforms of $f$ and $\fh$.  To see this, fix a polarisation
$Q$ and thus a \Ch\ transform $f^c$ of $f$ so that $\fh=\D_r f$ for
some $r\in\R^\times$.  As usual, take
\[
F=
\begin{pmatrix}
\fh\gi&f\\\gi&1
\end{pmatrix}
\]
so that
\[
A_\p=
\begin{pmatrix}
0&\d f\\-r\d f^c
\end{pmatrix}.
\]
Then $(f_{(\lambda)},\fh_{(\lambda)})$ is framed by
$F_\lambda:M\to\SL$ with
\[
\mc{F_\lambda}=A_\k+\lambda A_\p.
\]
Now replace $F_\lambda$ with the frame $\displaystyle F_\lambda
\begin{pmatrix}
\sqrt{\lambda}&0\\0&1/\sqrt{\lambda}
\end{pmatrix}$ which has \MC\ form
\[
A_\k+
\begin{pmatrix}
0&\d f\\-\lambda^2r\d f^c&0\end{pmatrix}
=A_\k+A_\p+(1-\lambda^2)r
\begin{pmatrix}
0&0\\\d f^c&0
\end{pmatrix}.
\]
Thus, by Theorem~\ref{th:12}, $f_{(\lambda)}=\T_{(1-\lambda^2)r}f$.
\begin{ex}
Contemplate the frame
\[
F_\lambda
\begin{pmatrix}
0&-1/\sqrt{\lambda}\\\sqrt{\lambda}&0
\end{pmatrix}
\]
of $(\fh_{(\lambda)},f_{(\lambda)})$ to conclude that
$\fh_{(\lambda)}=\T_{(1-\lambda^2)r}\fh$.
\end{ex}
To summarise:
\begin{thm}
The associated family of a Darboux pair $(f,\fh)$ consists
of $T$-transforms of the pair: 
\[
f_{(\lambda)}=\T_{(1-\lambda^2)r}f,\qquad
\fh_{(\lambda)}=\T_{(1-\lambda^2)r}\fh
\]
\end{thm}

\begin{rem}
The extraction of roots in our gauge transformations means we must
take $\lambda>0$.  However, since
\[
\tau A_\lambda=A_\k-\lambda A_\p=A_{-\lambda},
\]
$\tau F_\lambda$ and $F_{-\lambda}$ differ by a constant so that the
pairs $(f_{(\lambda)},\fh_{(\lambda)})$ and
$(f_{(-\lambda)},\fh_{(-\lambda)})$ differ by a M\"obius
transformation.  We shall have more to say about this symmetry below.
\end{rem}

\subsubsection{$\p$-flat maps are \Ch\ pairs}
\label{sec:p-flat-maps}

The alert reader will have noticed by now that there is a second way
to construct a pair of isothermic surfaces from a curved flat: the
\MC\ form of a flat frame of a curved flat is the derivative of a
$\p$-flat map $\psi:M\to\p$:
\[
A_\p=\d\psi.
\]
In our case, write
\[
\psi=
\begin{pmatrix}
0&f_0\\f_0^c&0
\end{pmatrix}
\]
for $f_0,f^c_0:M\to\R^n$.  Then $[\d\psi\wedge\d\psi]=0$ amounts to
\[
\d f_0\wedge\d f_0^c=0
\]
and its transpose so that a $\p$-map is precisely a dual pair of
isothermic surfaces!

It is important to emphasise that this pair is \emph{not} the Darboux
pair comprising the curved flat.  Rather, the two pairs are
$T$-transforms of each other: indeed, if the flat frame $F$ frames
the Darboux pair $(f,\fh)$ we have
\[
\mc{F}=
\begin{pmatrix}
0&\d f_0\\\d f^c_0&0
\end{pmatrix}
\]
so that $f=\T_1 f_0$ and, by Theorem~\ref{th:13}, $\fh=\T_1 f_0^c$.

Conversely, given a \Ch\ pair $(f_0,f_0^c)$, we integrate to obtain
the associated family of flat frames $F_\lambda$ with
\[
\mc{F_\lambda}=\lambda\begin{pmatrix}
0&\d f_0\\\d f^c_0&0
\end{pmatrix}.
\]
The $F_\lambda$ frame Darboux pairs $(f_{(\lambda)},\fh_{(\lambda)})$
and we argue as in Section~\ref{sec:spectr-deform-t} to prove:
\begin{ex}
$f_{(\lambda)}=\T_{\lambda^2}f_0$, $\fh_{(\lambda)}=\T_{\lambda^2}f_0^c$.
\end{ex}

As we shall see in Section~\ref{sec:extended-flat-frames}, if the
constants of integration are chosen correctly, we can recover
$(f_0,f_0^c)$ up to a translation from the frames $F_\lambda$ via the
Sym formula \cite{CieGolSym95}:
\[
\begin{pmatrix}
0&f_0\\f_0^c&0
\end{pmatrix}=
\dl{F_\lambda}
\]
so that a \Ch\ transform is a limit of Darboux transforms.

In conclusion, we have seen that an associated family of curved flats
in $Z$ amounts to the family of $T$-transforms (for $r>0$) of a \Ch\
pair of isothermic surfaces, each $T$-transform being, as we know
from Theorem~\ref{th:11}, a Darboux pair of isothermic surfaces.
However, the curved flat formulation gives us more: curved flats
admit a zero curvature formulation which means that we can apply the
powerful loop group approach to integrable systems and, in doing so,
find a completely different view-point on the topics we have been
studying.   It is to this that we now turn.

\section{Loop groups and B\"acklund transformations}
\label{sec:loop-groups-backlund}

We are going to show that associated families of curved flats (or
rather their flat frames) are the same as certain maps into an
infinite dimensional group $\Gp$ of holomorphic maps from $\C$ into a
complex Lie group.   Completely general principles, first enunciated
by Zakharov and his collaborators \cite{Zakilo78,ZakSha78}, then
allow us to construct a local action of a second infinite-dimensional
group $\Gm$ on these families.  In general, computation of this
action amounts to solving a Riemann--Hilbert problem but, as has been
made clear in a series of papers by Terng and Uhlenbeck
\cite{TerUhl98,TerUhl00,Uhl89,Uhl92}, the action of certain elements
of $\Gm$, the \emph{simple factors}, can be computed explicitly.  In
several geometric problems, the action of these simple factors amount
to known B\"acklund transformations.

We shall show that this is the case for isothermic surfaces: the
action of simple factors will turn out to be precisely by Darboux
transforms of the underlying \Ch\ pair.  This places our theory in a
well-understood context in integrable systems theory and, in
particular, general arguments of Terng--Uhlenbeck \cite{TerUhl00} can
be exploited to establish Bianchi permutability of Darboux
transforms.  In this way, we find a second approach to the results of
Section~\ref{sec:bianchi-perm-cliff}.

\subsection{Extended flat frames}
\label{sec:extended-flat-frames}

Henceforth $M$ will be simply connected with a fixed base-point $o\in
M$.

Let $G/K$ be a symmetric space.  Further let $G^\C$ be the
complexification of $G$ and denote by $g\mapsto \bar{g}$ the
conjugation across the real form $G$.  Thus $g\mapsto \bar{g}$ is the
anti-holomorphic involution on $G^\C$ with fixed set $G$.

Let $\psi:M\to\p$ be a $\p$-flat map and set $A_\p=\d\psi$.  We have
already seen how $\psi$ gives rise to a family of flat frames
$F_\lambda$ with
\[
\mc{F_\lambda}=\lambda A_\p,
\]
for $\lambda\in\R$.  We now extend this construction to
$\lambda\in\C$ and fix the constants of integration: for
$\lambda\in\C$, let $F_\lambda:M\to G^\C$ be the \emph{unique} map
with
\begin{gather*}
\mc{F_\lambda}=\lambda A_\p\\
F_\lambda(o)=1.
\end{gather*}
The existence of each $F_\lambda$ is guaranteed since $\lambda A_\p$
solves the \MC\ equations and $M$ is simply connected.

We note:
\begin{enumerate}
\item $F_0=1$ since $\mc{F_0}=0$ and $F_0(o)=1$.
\item For each $p\in M$, $\lambda\mapsto F_\lambda(p):\C\to G^\C$ is
holomorphic since $\lambda\mapsto\lambda A_\p$ is certainly
holomorphic as is $\lambda\mapsto F_\lambda(o)$.
\item For all $\lambda\in\C$,
\[
\overline{F_\lambda}=F_{\bl}
\]
or, equivalently, $F_\lambda:M\to G$ when $\lambda\in\R$.  This holds
since
\[
\overline{\lambda A_\p}=\bl A_p
\]
so that $\overline{F_\lambda}$ and $F_{\bl}$ have the same \MC\ form
and the same value at $o$ and so must coincide.
\item Similarly, since $\tau(\lambda A_\p)=-\lambda A_\p$, we
conclude that, for all $\lambda\in\C$,
\[
\tau F_\lambda=F_{-\lambda}.
\]
\end{enumerate}

We now change our point of view and assemble the $F_\lambda$ into a
single map $\Phi:M\to\Map(\C,G^\C)$ by setting
\[
\Phi(p)(\lambda)=F_\lambda(p).
\]
Observe that $\Phi$ takes values in the group $\Gp$ of
\emph{holomorphic} maps $g:\C\to G^\C$ satisfying
\begin{subequations}
\label{eq:64}
\begin{gather}
\label{eq:65}
g(0)=1,\\
\label{eq:66}
\tau g(\lambda)=g(-\lambda),\\
\label{eq:67}
\overline{g(\lambda)}=g(\bl),
\end{gather}
\end{subequations}
for all $\lambda\in\C$.  It is easy to see that $\Gp$ is a group
under point-wise multiplication.

\begin{defn}
A map $\Phi:M\to\Gp$ is an \emph{extended flat frame} if and only if
\begin{equation}
\label{eq:68}
\mc{\Phi}(\lambda)=\lambda A_\p
\end{equation}
with $A_\p\in\Omega^1\otimes\p$ independent of $\lambda$.

$\Phi$ is additionally said to be \emph{based} if $\Phi(o)=1$.
\end{defn}

The property of being an extended flat frame is characterised
entirely by the behaviour at $\lambda=\infty$ of $\mc{\Phi}$:
\begin{lem}\label{th:22}
$\Phi:M\to\Gp$ is an extended flat frame if and only if, for each
$p\in M$, $\mc{\Phi}_{|p}$ has a simple pole at $\lambda=\infty$.
\end{lem}
\begin{proof}
Let $\Phi:M\to\Gp$ and contemplate the power series expansion of
$\mc{\Phi}$:
\[
\mc{\Phi}=\sum_{n\geq0}\lambda^n A_n
\]
with $A_n\in\Omega^1\otimes\g^\C$.  The twisting and reality
conditions \eqref{eq:66} and \eqref{eq:67} force
\[
\tau\sum_{n\geq0}\lambda^n A_n=\sum_{n\geq0}(-\lambda)^n A_n,\qquad
\overline{\sum_{n\geq0}\lambda^n A_n}=\sum_{n\geq0}\bl^n A_n
\]
whence
\[
A_{2n}\in\Omega^1\otimes\k,\qquad A_{2n-1}\in\Omega^1\otimes\p.
\]
Moreover, $\Phi(0)\equiv1$ so that $A_0=0$.

Thus $\mc{\Phi}$ has a simple pole at $\lambda=\infty$ if and only if
all $A_n=0$ for $n>1$ which is the case precisely when
$\mc{\Phi}=\lambda A_1$ for some $A_1\in\Omega^1\otimes\p$.
\end{proof}

We can recover the generating $\p$-flat map up to translation from
$\Phi$ by a popular device known as the Sym formula:
\begin{prop}\label{th:23}
Let $\Phi$ be an extended flat frame with $\mc{\Phi}=\lambda A_\p$
and define $\psi_0:M\to\g$ by
\begin{equation}
\label{eq:69}
\psi_0=\dl{\Phi}
\end{equation}
Then
\begin{enumerate}
\item $\psi_0:M\to\p$;
\item $\d\psi_0=A_\p$.
\end{enumerate}
\end{prop}
\begin{proof}
We have $\tau\Phi(\lambda)=\Phi(-\lambda)$ and differentiating with
respect to $\lambda$ gives
\[
\tau\dl{\Phi}=-\dl{\Phi},
\]
that is, $\psi_0:M\to\p$.

Now view $\Phi$ as a map $M\times\C\to G^\C$ with \MC\ form
$\alpha$.  Then, for $p\in M$ and $X\in T_pM$, we have
\begin{align*}
\alpha_{(p,0)}(\del/\del\lambda)&=\psi_0(p);\\
\alpha_{(p,\lambda)}(X)&=\lambda A_\p(X).
\end{align*}
The \MC\ equations for $\alpha$ give
\[
\d \alpha_{(p,0)}(\del/\del\lambda,X)+
[\alpha_{(p,0)}(\del/\del\lambda),\alpha_{(p,0)}(X)]=0.
\]
However, $\alpha_{(p,0)}(X)=0$ since $\Phi(p)(0)=1$ for all $p\in M$
so we are left with
\[
\dl{\alpha(X)}-\d_X\alpha(\del/\del\lambda)=0,
\]
that is, $A_\p(X)=\d_X\psi_0$.
\end{proof}
Thus, if $\psi:M\to\p$ is a $\p$-flat map and $\Phi$ is the
corresponding \emph{based} extended flat frame, then 
\[
\dl{}\Phi(o)=0
\]
so that
\begin{equation}
\label{eq:70}
\psi=\dl{\Phi}+\psi(o).
\end{equation}
This gives us a bijective correspondence:
\begin{align*}
\{\text{$\p$-flat maps}\}&\to\{\text{based extended flat
frames}\}\times\p\\
\psi&\mapsto(\Phi,\psi(o))\\
\intertext{with inverse}
(\Phi,\xi)&\mapsto\dl{\Phi}+\xi.
\end{align*}

The Sym formula has geometric content: for $\lambda\in\R$, let
$\phi_\lambda:M\to G/K$ be the curved flat framed by
$\Phi(\lambda)$.  In particular $\phi_0\equiv eK$, the identity
coset.  With the usual identification $T_{eK}G/K\cong\p$, one sees
that
\[
\psi_0=\dl{\phi_\lambda}.
\]
In particular, in the isothermic surface case, we have $\p\cong
T_0S^n\oplus T_\infty S^n$ and an associated family of Darboux pairs
$(f_{(\lambda)},\fh_{(\lambda)})$ with
\[
f_{(0)}\equiv0,\qquad \fh_{(0)}\equiv\infty.
\]
The generating Christoffel pair $(f,f^c)$ are recovered by ``blowing
up'' their $T$-transforms as $\lambda\to 0$:
\begin{align*}
f&=\dl{f_{(\lambda)}}:M\to T_0S^n;\\
\fh&=\dl{\fh_{(\lambda)}}:M\to T_\infty S^n.
\end{align*}

\subsection{The dressing action}
\label{sec:dressing-action}

We are going to define a local action of a group of rational maps on the
set of extended flat frames and so, eventually, on the set of
$\p$-flat maps.  Our action will be by point-wise application of a
local action on $\Gp$ which we now describe.

Let $\G$ denote the group of holomorphic maps
$g:\dom(g)\subset\P^1\to G^\C$ of affine subsets of the Riemann sphere
which are twisted and real in the sense that
\begin{subequations}
\label{eq:71}
\begin{align}
\tau g(\lambda)=g(-\lambda),\\
\overline{g(\lambda)}=g(\bl),
\end{align}
for all $\lambda\in\dom(g)$.
\end{subequations}
Clearly $\Gp$ is a subgroup of $\G$.  We define a second subgroup
$\Gm$ by
\[
\Gm=\{g\in\G:\text{$g$ is rational on $\P^1$ and holomorphic near $\infty$}\}.
\]
Thus $\Gp$ consists of those elements of $\G$ which are holomorphic
on $\C$ while $\Gm$ consists of those which are rational\footnote{The
restriction to rational maps is not really necessary: one could work
with the group of germs at $\infty$ of maps to $G^\C$ with
\eqref{eq:71}.  While not appropriate here, such generality is
necessary in some contexts, see \cite{BurPed95} for a discussion in a
related situation.} and holomorphic near $\infty$.

\begin{lem}\label{th:24}
$\Gp\cap\Gm=\{1\}$.
\end{lem}
\begin{proof}
If $g\in\Gp\cap\Gm$ then $g$ is holomorphic on $\P^1$ and so is
constant.  Moreover $g(0)=1$ whence $g=1$.
\end{proof}

The basis of our action is the Birkhoff-Grothendieck decomposition
theorem in a formulation due to Pressley--Segal \cite{PreSeg86}:
\begin{thm}
\label{th:25}
Set $\mathcal{U}=\Gp\Gm$.  Then $\mathcal{U}$ is a dense
open\footnote{The reader may object that I have not topologised $\G$:
in fact, the compact open topology will do (or any stronger one).}
subset of $\G$.
\end{thm}

Thus $g\in\mathcal{U}$ if and only if we can write
\begin{equation}
\label{eq:72}
g=g_+g_-
\end{equation}
with $g_\pm\in\G^\pm$.
\begin{ex}\label{ex:9}
Use Lemma~\ref{th:24} to show that the decomposition \eqref{eq:72} is
unique when it exists.
\end{ex}

For $g_-\in\Gm$, set
$\mathcal{U}_{g_-}=\gi_-\mathcal{U}g^{\vphantom{1}}_-\cap\Gp$:
this is an open neighbourhood of $1$ in $\Gp$.
\begin{lem}
$g_+\in\mathcal{U}_{g_-}$ if and only if there are unique
$\gh_\pm\in\G^\pm$ such that
\begin{equation}
\label{eq:73}
g_-g_+=\gh_+\gh_-
\end{equation}
on $\C\cap\dom(g_-)$.
\end{lem}
\begin{proof}
If \eqref{eq:73} holds then
\[
g_+=\gi_-\gh^{\vphantom{1}}_+\gh^{\vphantom{1}}_-
=\gi_-\gh^{\vphantom{1}}_+\gh^{\vphantom{1}}_-\gi_-g^{\vphantom{1}}_-
\in\gi_-\Gp\Gm g^{\vphantom{1}}_-\cap\Gp=
\mathcal{U}_{g_-}.
\]
Conversely, if $g_+\in\mathcal{U}_{g_-}$ then
$g^{\vphantom{1}}_-g^{\vphantom{1}}_+\gi_-\in\mathcal{U}$ so we can write
\[
g^{\vphantom{1}}_-g^{\vphantom{1}}_+\gi_-=h_+h_-
\]
with $h_\pm\in\G^\pm$.  Now put $\gh_+=h_+$ and $\gh_-=h_-g_-$.  

The uniqueness assertion is proved as in Exercise~\ref{ex:9}.
\end{proof}

\begin{notation}
Write $g_-\#g_+$ for $\gh_+$ in \eqref{eq:73}.
\end{notation}

Thus $g_-\#g_+=g_-g_+\gh_-^{-1}$.
\begin{ex}
\label{ex:10}
Show:
\begin{enumerate}
\item $\mathcal{U}_1=\Gm$ and $1\#g_+=g_+$ for all $g_+\in\Gp$.
\item For all $g_-\in\Gm$, $g_-\#1=1$.
\end{enumerate}
\end{ex}

Now let $g_1,g_2\in\Gm$, $g_+\in\Ug[1]$ and suppose $g_1\#g_+\in\Ug[2]$
so that $g_2\#(g_1\# g_+)$ is defined.  This means we have
$\gh_1,\gh_2\in\Gm$ such that
\[
g^{\vphantom{1}}_2(g^{\vphantom{1}}_1g^{\vphantom{1}}_+\gh_1^{-1})\gh_2^{-1}=g_2\#(g_1\# g_+)\in\Gp
\]
whence
\[
(g_2g_1)g_+=(g_2\#(g_1\# g_+))\gh_2\gh_1.
\]
Since $\gh_2\gh_1\in\Gm$, we conclude that
$g_+\in\mathcal{U}_{g_2g_1}$ and that
$(g_2g_1)\#g_+=g_2\#(g_1\#g_+)$.
Taking this together with Exercise~\ref{ex:10}, we conclude:
\begin{thm}
$g_-\#g_+$ defines a local action of $\Gm$ on $\Gp$.
\end{thm}

Now let $\Phi:M\to\Gp$ be a map and $g_-\in\Gm$.  Define
$g_-\#\Phi:\Phi^{-1}(\Ug)\subset M\to\Gp$ by
\[
(g_-\#\Phi)(p)=g_-\#\bigl(\Phi(p)\bigr).
\]
The whole point of this is contained in the following theorem:
\begin{thm}\label{th:26}
If $\Phi:M\to\Gp$ is a (based) extended flat frame then so is
$g_-\#\Phi$.
\end{thm}
\begin{proof}
Set $\hat{\Phi}=g_-\#\Phi$ and let $\hat{A}$ be its \MC\ form.  By
Lemma~\ref{th:22}, we must show that $\hat{A}$ has a simple pole at
$\lambda=\infty$.  However, in a punctured neighbourhood of $\infty$,
we have
\[
\hat{\Phi}=g^{\vphantom{1}}_-\Phi\gh_-^{-1}
\]
with $\gh_-^{-1}:\Phi^{-1}(\Ug)\to\Gm$ so that
\[
\hat{A}=\Ad\gh^{\vphantom{1}}_-(\mc{\Phi})-\d\gh^{\vphantom{1}}_-\gh_-^{-1}.
\]
Since $\gh_-$ is holomorphic at $\lambda=\infty$, we immediately
conclude that $\hat{A}$ has the same pole at $\infty$ as $\mc{\Phi}$.

Finally, $g_-\#1=1$ so that $\hat{\Phi}$ is based if and only if
$\Phi$ is.
\end{proof}

\begin{rem}
To get this far, we have used very little of the specifics of the
situation.  To get a local action of $\Gm$ on $\Gp$ we only used that
$\Gp\Gm$ is open in $\G$ with $\Gp\cap\Gm=\{1\}$.  Moreover the
argument of Theorem~\ref{th:26} is also very general: the only
ingredient is that membership of the class of extended frames is
determined by the pole behaviour of the \MC\ form.  For then, the
pointwise action of the group of maps holomorphic near these poles will
preserve that class.  Thus one can use exactly the same techniques to
produce actions of such groups in a variety of geometric problems.
See \cite{BurGue97,BurPed95,Uhl89} among others for the case of
harmonic maps and the work of Terng--Uhlenbeck
\cite{TerUhl98,TerUhl00} for many other examples.
\end{rem}

We would like our action of $\Gm$ on based extended flat frames to
induce an action on $\p$-flat maps.  However, since the frame only
determines the $\p$-flat map up to translation, we must work with a
slightly smaller group which has an action on $\p$ also.

For this, define $\Gmb\subset\Gm$ by
\[
\Gmb=\{g\in\Gm:g(0)=1\}.
\]
For $g_-\in\Gmb$ and $\psi:M\to\p$ a $\p$-flat map, let $\Phi$ be the
based extended flat frame with $\mc{\Phi}=\lambda\d\psi$ and define
$g_-\#\psi:\Phi^{-1}(\Ug)\to\p$ by
\begin{equation}
\label{eq:74}
g_-\#\psi=\psi(o)+\dl{}(g_-\#\Phi)-\dl{g_-}.
\end{equation}
Note that $\del g_-/\del\lambda|_{\lambda=0}\in\p$ so that the right
hand side is indeed $\p$-valued.

$g_-\#\psi$ differs from $\del/\del\lambda|_{\lambda=0}(g_-\#\Phi)$
by constants so that, by Proposition~\ref{th:23}, $g_-\#\psi$ is again a
$\p$-flat map.
\begin{ex}
\label{ex:11}
Show that, for $\psi$ a $\p$-flat map and $g_1,g_2\in\Gmb$,
\begin{gather*}
1\#\psi=\psi\\
g_1\#(g_2\#\psi)=(g_1g_2)\#\psi
\end{gather*}
whenever the left hand side is defined.
\end{ex}

Thus we conclude:
\begin{thm}
There is a local action of $\Gmb$ on $\p$-flat maps given by \eqref{eq:74}.
\end{thm}

We obtain a more efficient formula for this action as follows: write
\[
g_-\Phi=(g_-\#\Phi)\gh_-
\]
with $\gh_-:\Phi^{-1}(\Ug)\to\Gmb$ so that
\[
g_-\#\Phi=g^{\vphantom{1}}_-\Phi\gh_-^{-1}.
\]
Thus
\[
\dl{}(g_-\#\Phi)=\dl{g_-}+\dl{\Phi}+\dl{\gh_-^{-1}}.
\]
Now \eqref{eq:70} gives
\[
\dl{\Phi}=\psi-\psi(o)
\]
and feeding all this into \eqref{eq:74} gives:
\begin{equation}
\label{eq:75}
g_-\#\psi=\psi+\dl{\gh_-^{-1}}.
\end{equation}

\subsection{Simple factors}
\label{sec:simple-factors}

Given $g_-\in\Gm$ and $\Phi$ an extended flat frame, a basic problem
is to compute $g_-\#\Phi$.  This amounts to performing the
factorisation
\[
g_-\Phi(p)=\gh_+(p)\gh_-(p)
\]
for each $p\in M$ and, in general, this is a Riemann--Hilbert
problem.  Part of the philosophy of Terng--Uhlenbeck is that there
are special elements of $\Gm$, the \emph{simple factors}, for which
one can explicitly perform the factorisation by algebra alone and,
moreover, that the action of these factors amount to B\"acklund-type
transformations of the underlying geometric problem.  Of course, the
Art in this approach is to put one's hands on these simple factors!

Let us look for some hints.  We are given $g_\pm\in\G^\pm$ and seek
$\gh_\pm\in\G^\pm$ so that
\begin{equation}
\label{eq:76}
g_-g_+=\gh_+\gh_-.
\end{equation}
First observe that any $g_-\neq 1\in\Gmb$ must have some
singularities in $\C^\times$, that is, $\lambda\in\C^\times$ where
$g$ either fails to be defined or fails to be invertible.  In view of
the twisting and reality conditions \eqref{eq:71}, if $\alpha$ is
such a singularity, so is $-\alpha$ and $\bar{\alpha}$.

Secondly, rearrange \eqref{eq:76} to get
\[
\gh_-=\gh_+^{-1}g^{\vphantom{1}}_-g^{\vphantom{1}}_+
\]
with both $g_+,\gh_+$ holomorphic on $\C$.  Thus $\gh_-$ has the same
singularities as $g_-$.

The idea now is to work with $g_-$ having the minimum number of
singularities.  In our case, this number is two and we are
contemplating $g_-$ with poles at $\pm\alpha$ and demand that either
$\bar{\alpha}=\alpha$ so that $\alpha\in\R$ or $\bar{\alpha}=-\alpha$
so that $\alpha\in\I\R$.

To get further, we begin by considering the case where $G$ is
compact.  Here the situation reduces to one which is completely
understood.  When $G$ is compact, we can have no singularities on
$\R$ and thus $\alpha\in\I\R$.  Now use a linear fractional
transformation to move the singularities at $\pm\alpha$ to $0$ and
$\infty$: define $t_\alpha:\P^1\to\P^1$ by
\[
t_\alpha(\lambda)=\frac{\alpha-\lambda}{\alpha+\lambda}
\]
so that
\begin{gather*}
t_\alpha(\alpha)=0,\quad t_\alpha(-\alpha)=\infty,\quad
t_\alpha(0)=1\\
t_\alpha(\R)\subset S^1=\{\lambda:\abs{\lambda}=1\}
\end{gather*}
We write
\[
g_-=h_-\circ t_\alpha
\]
for $h_-:\C^\times\to G^\C$.  Since $g_-$ is rational, $h_-$ is a
Laurent polynomial and, moreover, we have $h_-(1)=1$ and
$h(S^1)\subset G$.  Otherwise said, $h_-$ lies in the based algebraic
loop group $\Um$ of Laurent polynomial maps $h:\C^\times\to G^\C$
satisfying
\begin{gather*}
h(1)=1\\
\overline{h(\lambda)}=h(1/\bl).
\end{gather*}
This group features in a factorisation problem which can always be
solved: let $\Up$ denote the group of maps $h_+$ to $\G^\C$ which are
defined and holomorphic near $0$ and $\infty$ and have the reality
condition
\[
\overline{h_+(\lambda)}=h_+(1/\bl),
\]
for $\lambda\in\dom(h_+)$.
It follows from the results of Pressley-Segal \cite{PreSeg86} that,
for $h_+\in\Up$, $h_-\in\Um$, there is always a unique decomposition
\[
h_+h_-=\hh_-\hh_+
\]
with $\hh_-\in\Um$ and $\hh_+\in\Up$.  Just as before, we set
\[
h_+\cdot h_-=\hh_-
\]
to get a (now global) action of $\Up$ on $\Um$ which is well
understood.

The relevance of all this to our own factorisation problem is that,
after taking inverses and moving the poles with $t_\alpha$, our
decomposition problem becomes that of Pressley--Segal:
\begin{ex}
For $g_\pm\in\G^\pm$ with $g_-$ having only singularities at
$\pm\alpha$, write
\[
g_\pm=h_\pm\circ t_\alpha
\]
so that $h_-\in\Um$ and $h_+\in\Up$ (since $g_+$ is holomorphic at
$\pm\alpha$).  Then
\[
g_-\# g_+=(h_+^{-1}\cdot h_-^{-1})^{-1}\circ t_\alpha^{-1}.
\]
\end{ex}
In particular, we deduce
\begin{prop}
If $G$ is compact and $g_-$ has only two poles then $\Ug=\G^+$.
\end{prop}

There is more: in this setting, the action of such a $g_-$ is, in principle,
computable algebraically:
\begin{fact}
The orbits of $\Up$ on $\Um$ are finite-dimensional: they form the
\emph{Bruhat decomposition} of $\Um$ \cite{PreSeg86}.  In fact,
$h_+\cdot h_-$ depends only on a finite jet of $h_+$ at $\lambda=0$.
\end{fact}

As a consequence, $g_-\#g_+$ can be computed from $g_-$ and a finite
jet at $\alpha$ of $g_+$.  The maximally desirable situation is when
only the $0$-jet $g_+(\alpha)$ is involved: again, this amounts to a
feature of the Bruhat decomposition.
\begin{fact}{\cite{PreSeg86,BurGue97}}
$h_+\cdot h_-$ depends only on $h_+(0)$ if and only if $h_-:\C^\times\to
G^C$ is a \emph{homomorphism} such that $\Ad
h_-:\C^\times\to\Ad(G^\C)$ has simple poles only.  In this case,
$\hh_-$ is another homomorphism in the same (real) conjugacy class. 
\end{fact}

We therefore conclude that if we wish to be able to compute
$g_-\#\Phi(p)$ from just $g_-$ and the value $\Phi(p)(\alpha)$ then
we are compelled to take
\begin{equation}
\label{eq:77}
g_-(\lambda)=\gamma\bigl(\frac{\alpha-\lambda}{\alpha+\lambda}\bigr)
\end{equation}
where $\Ad\gamma:\C^\times\to\Ad(G^\C)$ is a homomorphism with only
simple poles.  This last is a very restrictive condition: it implies
that the real conjugacy class of $\gamma$ is a Hermitian symmetric
space and so excludes the exceptional Lie groups $\mathrm{G}_2$,
$\mathrm{F}_4$ and $\mathrm{E}_8$.  in fact, there is more: we need
$\tau g_-(\lambda)=g_-(-\lambda)$ which amounts to demanding
\[
\tau\gamma(\mu)=\gamma(1/\mu),
\]
for $\mu\in\C^\times$ and this excludes many symmetric
spaces\footnote{For example, if $G/K$ is a projective space $\R\P^n$,
$\C\P^n$ or $\H\P^n$, then such $\gamma$ exist only when $n=1$.}.

Be that as it may, for $G$ compact, we have shown that any $g_-$ with
the properties we want must be of the form \eqref{eq:77}.  For $G$
non-compact, we make this our ansatz:
\begin{defn}
$g_-\in\Gmb$ is a \emph{simple factor} if it is the form
\[
g_-=\gamma\circ t_\alpha
\]
with $\alpha^2\in\R$ and $\gamma:\C^\times\to G^{\C}$ a homomorphism for
which $\Ad\gamma$ has simple poles.
\end{defn}

It turns out that simple factors retain their desirable property of
having algebraically computable action even for non-compact $G$.
However, to develop the theory any further in this general setting
will take us too far afield so we now turn to the case of relevance
to isothermic surfaces.

\subsection{Simple factors for $S^n\times S^n\setminus\Delta$}
\label{sec:simple-fact-sntim}

We are going to classify the simple factors for $G=\Op$ and so begin by
determining the homomorphisms $\gamma:\C^\times\to\Op^\C=\O[n+2,\C]$ for which
$\Ad\gamma$ has simple poles.

Let $\gamma:\C^\times\to\O[n+2,\C]$ be a homomorphism.  There is a
decomposition of $\C^{n+2}$ into common eigenspaces of the $\gamma(\lambda)$:
\[
\C^{n+2}=\oplus_{i=-k}^k V_i
\]
so that, with $\pi_i$ the projection onto $V_i$ along $\oplus_{i\neq j}V_j$, we
have
\[
\gamma(\lambda)=\sum_{i=-k}^k \lambda^i\pi_i
\]
(we must allow the possibility that some $V_i=\{0\}$).  Since
$\gamma(\lambda)\in\O[n+2,\C]$, we have $V_i\perp V_j$ for $i+j\neq 0$ so that
each $V_i$ is isotropic for $i\neq 0$, $\dim V_i=\dim V_{-i}$ and
$V_0^\perp=\oplus_{i\neq 0}V_i$.  As $\O[n+2,\C]$-modules,
$\o[n+2,\C]\cong\bigwedge^2\C^{n+2}$ via
\[
(u\wedge v)w=(u,w)v-(v,w)u
\]
and using this identification we immediately see that $\Ad\gamma(\lambda)$ has
eigenvalues $\lambda^{2i}$ on $\bigwedge^2 V_i$ and $\lambda^{i+j}$ on
$V_i\otimes V_j$, $i\neq j$.  Thus $\Ad\gamma$ has simple poles exactly when
$k=1$ and $\dim V_1=1$ (to ensure $\bigwedge^2 V_1=\{0\}$).  We are therefore
working with $\gamma$ of the form
\[
\gamma(\lambda)=\lambda\pi_+ +\pi_0+\lambda^{-1}\pi_-
\]
corresponding to a decomposition
\[
\C^{n+2}=L_+\oplus L_0\oplus L_-
\]
with $L_\pm$ $1$-dimensional isotropic subspaces and $L_0=(L_+\oplus
L_-)^\perp$.

The key to computing the dressing action of the corresponding simple factor is
the following lemma:
\begin{lem}
\label{th:27}
Let $\gamma(\lambda)=\lambda\pi_+ +\pi_0+\lambda^{-1}\pi_-$ and
$\hg=\lambda\hpi_+ +\hpi_0+\lambda^{-1}\hpi_-$ be homomorphisms as
above with $\Ad\gamma$, $\Ad\hg$ having simple poles and let
\[
\C^{n+2}=L_+\oplus L_0\oplus L_-=\hL_+\oplus \hL_0\oplus \hL_-
\]
be the corresponding eigenspace decompositions. 

Let $E$ be the germ at $0$ of a map into $\O[n+2,\C]$.  Then $\gamma E\hg^{-1}$
is holomorphic and invertible at $0$ if and only if
\[
\hL_+=E(0)^{-1}L_+.
\]
\end{lem}
\begin{proof}
Write $E$ as a power series:
\[
E(\lambda)=\sum_{k\geq0}\lambda^k E_k.
\]
Comparing coefficients of $\lambda$, we see that $\gamma E\hg^{-1}$
is holomorphic at zero if and only if 
\begin{enumerate}
\item $\pi_-E_0\hpi_+=0$ (this is the coefficient of $\lambda^{-2}$);
\item $\pi_0E_0\hpi_+=\pi_-E_0\hpi_0=\pi_-E_1\hpi_+=0$ (these are the
components of the coefficient of $\lambda^{-1}$).
\end{enumerate}
Now observe that
\[
\pi_-E_0\hpi_+=\pi_0E_0\hpi_+=0
\]
if and only if $E_0\hL_+=L_+$ and then, since $E_0\in\O[n+2,\C]$,
\[
L_+\oplus L_0=L_+^\perp=E_0(\hL_+^\perp)=E_0(\hL_+\oplus\hL_0)
\]
whence $\pi_-E_0\hpi_0$ vanishes automatically.

This leaves the term involving $E_1$.  However, when $E_0\hL_+=L_+$,
we have $E_1\hL_+=E_1E_0^{-1}L_+$ and
\[
E_1E_0^{-1}=\dl{E}\in\o[n+2,\C]
\]
so that $E_1E_0^{-1}$ is skew-symmetric.  Thus, since $L$ is
$1$-dimensional\footnote{It is at this point of the argument that we
are really using the hypothesis that $\Ad\gamma$ has only simple
poles.}, we have
\[
(E_1E_0^{-1}L_+,L_+)=0
\]
giving $\pi_-E_1\hpi_+=0$.

Thus $\gamma E\hg^{-1}$ is holomorphic at zero if and only if
$\hL_+=E(0)^{-1}L_+$.  The invertibility now follows by applying this
result to $\hg E^{-1}\gamma^{-1}$.
\end{proof}

Fix such a $\gamma=\lambda\pi_++\pi_0+\lambda^{-1}\pi_-$ and set
$g_-=\gamma\circ t_\alpha$:
\[
g_-(\lambda)=\gamma\bigl(\frac{\alpha-\lambda}{\alpha+\lambda}\bigr),
\]
with $\alpha^2\in\R^\times$.  Thus
$g_-:\P^1\setminus\{\pm\alpha\}\to\O[n+2,\C]$ and $g_-(0)=1$.  We
want $g_-\in\Gmb$ which means imposing two further conditions:
firstly, we must have
\[
\tau g_-(\lambda)=g_-(-\lambda)
\]
or, equivalently,
\[
\tau\gamma(\lambda)=\gamma(1/\lambda).
\]
In our setting, $\tau$ is conjugation by the reflection
$\rho:\C^{n+2}\to\C^{n+2}$ in $\R^n=(\R^{1,1})^\perp$ so that this
condition reads
\[
\rho L_+=L_-.
\]
In particular, this forces $\rho L_+\neq L_+$ and shows that $\gamma$
is completely determined by $L_+$ since $L_-=\rho L_+$ and
$L_0=(L_+\oplus\rho L_+)^\perp$.

Secondly, we must impose the reality condition
\[
\overline{g_-(\lambda)}=g_-(\bl)
\]
which amounts to
\[
\overline{\gamma(\lambda)}=
\begin{cases}
\gamma(\bl)&\text{if $\alpha\in\R^\times$;}\\
\gamma(1/\bl)&\text{if $\alpha\in\I\R^\times$;}
\end{cases}
\]
or, equivalently,
\[
\overline{L_+}=
\begin{cases}
L_+&\text{if $\alpha^2>0$;}\\
L_-&\text{if $\alpha^2<0$.}
\end{cases}
\]
Now we can put all this together: for $L\subset\C^{n+2}$ a
$1$-dimensional isotropic subspace with $\rho L\neq L$, let $\gamma_L$ be the
homomorphism $\C^\times\to\O[n+2,\C]$ given by
\[
\gamma_L(\lambda)=\lambda\pi_++\pi_0+\lambda^{-1}\pi_-
\]
with $\im\pi_+=L$, $\im\pi_-=\rho L$ and $\im\pi_0=(L\oplus\rho
L)^\perp$.  Further, for $\alpha\in\C^\times$, set
$p_{\alpha,L}=\gamma_L\circ t_\alpha$ so that
\[
p_{\alpha,L}(\lambda) =\frac{\alpha-\lambda}{\alpha+\lambda}\pi_+
+\pi_0+\frac{\alpha+\lambda}{\alpha-\lambda}\pi_-.
\]
We have shown that the simple factors in $\Gmb$ are precisely the
$p_{\alpha,L}$ with either
\begin{enumerate}
\item $\alpha^2>0$ and $L=\ell^\C$, the complexification of
$\ell\in\PL$ with $\rho\ell\neq\ell$, or,
\item $\alpha^2<0$ and $L$ is the complexification of a light-line
$\ell$ in $\R^n\oplus\I\R^{1,1}$ with $\rho\ell\neq\ell$.
\end{enumerate}

With all this in hand, we can now compute the dressing action of our
simple factors.  With an eye to proving Bianchi permutability, we
formulate a slightly more general result:
\begin{prop}
\label{th:28}
Let $p_{\alpha,L}\in\Gmb$ and let $E$ be a germ at $\alpha$ of a
holomorphic map into $\O[n+2,\C]$ such that
\[
\overline{E(\lambda)}=E(\bl),\qquad\tau E(\lambda)=E(-\lambda).
\]
Suppose further that $\rho(E(\alpha)^{-1}L)\neq E(\alpha)^{-1}L$.
Then
\begin{enumerate}
\item $p_{\alpha,E(\alpha)^{-1}L}\in\Gmb$;
\item $p^{\vphantom{1}}_{\alpha,L}Ep^{-1}_{\alpha,E(\alpha)^{-1}L}$
is holomorphic and invertible at $\alpha$.
\end{enumerate}
\end{prop}
\begin{proof}
For the first assertion we must establish the reality condition for
$p_{\alpha,E(\alpha)^{-1}L}$ and there are two cases.  First, if
$\alpha\in\R$, we must show that
$\overline{E(\alpha)^{-1}L}=E(\alpha)^{-1}L$.  However, in this case,
$\overline{L}=L$ and $\overline{E(\alpha)}=E(\alpha)$ so this follows
immediately.

When $\alpha\in\I\R$, we must show that $\overline{E(\alpha)}=\rho
E(\alpha)$ and, in this case, we have $\overline{L}=\rho L$ while
\[
\overline{E(\alpha)}=E(\overline{\alpha})=E(-\alpha)=\tau E(\alpha)=
\rho\circ E(\alpha)\circ\rho^{-1}.
\]
Thus
\[
\overline{E(\alpha)^{-1}L}=\overline{E(\alpha)^{-1}}\rho L=\rho
E(\alpha)^{-1}L
\]
as required.

The second assertion follows at once from Lemma~\ref{th:27}:
\[
p^{\vphantom{1}}_{\alpha,L}Ep^{-1}_{\alpha,E(\alpha)^{-1}L}=
(\gamma_L\circ t_\alpha)E(\gamma_{E(\alpha)^{-1}L}\circ
t_\alpha)^{-1}
\]
which is holomorphic at $\alpha$ if and only if
\[
\gamma^{\vphantom{1}}_L (E\circ t_\alpha^{-1})\gamma_{E(\alpha)^{-1}L}^{-1}
\]
is holomorphic at $0$.  However, $E\circ t^{-1}_\alpha$ is
holomorphic at $0$ with value $E(\alpha)$ there so Lemma~\ref{th:27}
applies.
\end{proof}

As a corollary we have:
\begin{thm}
\label{th:29}
$\mathcal{U}_{p_{\alpha,L}}=
\{g_+\in\Gp:g_+(\alpha)^{-1}L\neq\<v_0\>^\C,\<v_\infty\>^\C\}$ and,
for $g_+\in\mathcal{U}_{p_{\alpha,L}}$,
\begin{equation}
\label{eq:78}
p^{\vphantom{1}}_{\alpha,L}\# g_+=
p^{\vphantom{1}}_{\alpha,L}g_+ p^{-1}_{\alpha,g_+^{-1}(\alpha)L}.
\end{equation}
\end{thm}
\begin{proof}
Let $g_+\in\Gp$ be such that
$g_+(\alpha)^{-1}L\neq\<v_0\>^\C,\<v_\infty\>^\C$.  The first part of
Proposition~\ref{th:28} assures us that
$p_{\alpha,g_+^{-1}(\alpha)L}\in\Gmb$ so all we need do is see that
$p_{\alpha,L}\# g_+$ given by \eqref{eq:78} defines an element of
$\Gp$.  It is clear that $p_{\alpha,L}\# g_+$ has the reality and
twisting conditions as it is a product of maps with these conditions
so the only issue is that of holomorphicity and invertibility at
$\pm\alpha$.  However, holomorphicity at $\alpha$ follows at once
from Proposition~\ref{th:28} and then we get holomorphicity at
$-\alpha$ from the twisting condition:
\[
\tau(p_{\alpha,L}\# g_+)(\lambda)=p_{\alpha,L}\# g_+(-\lambda).
\]
\end{proof}

\begin{ex}
Complete the proof of Theorem~\ref{th:29} by showing that if
$g_+\in\Gp$ has $g_+^{-1}(\alpha)L=\<v_0\>^\C$ or $\<v_\infty\>^\C$
then $g_+\not\in\mathcal{U}_{p_{\alpha,L}}$.
\end{ex}

\subsection{The action of simple factors on Christoffel pairs}
\label{sec:acti-simple-fact}

We are finally in a position to compute the dressing action of simple
factors on \Ch\ pairs of isothermic surfaces.  Let us begin by
recalling all the ingredients: a $\p$-flat map $\psi:M\to\p$ is the
same as a \Ch\ pair $(f,f^c)$:
\[
\psi=
\begin{pmatrix}
0&f\\f^c&0.
\end{pmatrix}
\]
$g_-\in\Gmb$ acts on $\psi$ by \eqref{eq:75}:
\[
g_-\#\psi=\psi+\dl{}\gh^{-1}_-
\]
where $\gh_-:M\to\Gmb$ comes from the factorisation
\[
g_-\Phi=(g_-\#\Phi)\gh_-
\]
and $\Phi:M\to\Gp$ solves
\begin{gather*}
\mc{\Phi}=\lambda\d \psi=\lambda
\begin{pmatrix}
0&\d f\\\d f^c&0
\end{pmatrix},\\
\Phi(o)=1.
\end{gather*}
Now take $g_-=p_{\alpha,L}$.  Then Theorem~\ref{th:29} gives
\[
g_-\#\Phi=p^{\vphantom{1}}_{\alpha,L}\Phi p^{-1}_{\alpha,\Phi(\alpha)^{-1}L}
\]
so that $\gh_-=p_{\alpha,\Phi(\alpha)^{-1}L}$ and we have
\begin{equation}
\label{eq:79}
p_{\alpha,L}\#
\begin{pmatrix}
0&f\\f^c&0
\end{pmatrix}
=\begin{pmatrix}
0&f\\f^c&0
\end{pmatrix}+\dl{}p^{-1}_{\alpha,\Phi(\alpha)^{-1}L}.
\end{equation}
All that remains to do is to compute the second summand in
\eqref{eq:79}.  For this, write
\[
\gamma_{\Phi^{-1}(\alpha)L}(\lambda)=\lambda\hpi_++\hpi_0+\lambda^{-1}\hpi_-
\]
so that $\im\hpi_+=\Phi^{-1}(\alpha)L$.  Then
\[
p_{\alpha,\Phi(\alpha)^{-1}L}(\lambda)=
\frac{\alpha-\lambda}{\alpha+\lambda}\hpi_+
+\hpi_0+\frac{\alpha+\lambda}{\alpha-\lambda}\hpi_-
\]
so that
\[
\dl{}p^{-1}_{\alpha,\Phi(\alpha)^{-1}L}=\frac{2}{\alpha}(\hpi_+-\hpi_-).
\]
\begin{lem}
Fix $\omega_o\in L^\times$ and set
$\omega=\Phi^{-1}(\alpha)\omega_o:M\to\C^{n+2}$.  Then:
\begin{enumerate}
\item $\omega$ is the unique solution of
\begin{subequations}
\begin{gather}
\d\omega+\mc{\Phi_\alpha}\omega=0\label{eq:80}\\
\omega(o)=\omega_o.\label{eq:81}
\end{gather}
\end{subequations}
\item Viewing $\o$ as $[\lor,\lor]\subset\Cl_{n+1,1}$,
\[
\hpi_+-\hpi_-=\half\frac{[\omega,\rho\omega]}{\{\omega,\rho\omega\}}.
\]
\end{enumerate}
\end{lem}
\begin{proof}
We have $\omega_o=\Phi(\alpha)\omega$ and differentiating gives
\[
0=\d\Phi(\alpha)\omega+\Phi(\alpha)\d\omega
\]
whence \eqref{eq:80}.  Further, $\Phi(o)(\alpha)=1$ whence \eqref{eq:81}.

For the second part, recall that under the isomorphism
$[\lor,\lor]\cong\o$, $\xi\in[\lor,\lor]$ acts on $\lor$ by $v\mapsto
[\xi,v]$.  We must therefore show that, with
$\xi=\half[\omega,\rho\omega]/\{\omega,\rho\omega\}$, we have
\[
[\xi,\omega]=\omega,\quad[\xi,\rho\omega]=-\rho\omega,\quad[\xi,v]=0,
\]
for $v\perp\<\omega,\rho\omega\>$.  For
$v\perp\<\omega,\rho\omega\>$, $v$ anti-commutes with both $\omega$
and $\rho\omega$ and so commutes with $[\omega,\rho\omega]$.  Again,
using $\omega^2=0$, we have
\begin{align*}
[[\omega,\rho\omega],\omega]&=
(\omega\rho\omega-\rho\omega\omega)\omega-
\omega(\omega\rho\omega-\rho\omega\omega)\\
&=2\omega\rho\omega\omega=2\{\omega,\rho\omega\}\omega.
\end{align*}
Similarly, we have
\[
[[\omega,\rho\omega],\rho\omega]=-2\{\omega,\rho\omega\}\rho\omega.
\]
\end{proof}

Write
\[
\omega=
\begin{pmatrix}
v&s\\t&-v
\end{pmatrix}
\]
so that $v\in\C^n$ and $s,t\in\C$ with $v^2+st=0$.
\begin{ex}
\[
\hpi_+-\hpi_-=\half
\begin{pmatrix}
0&v/t\\t/v&0
\end{pmatrix}.
\]
\end{ex}
Thus, setting $h=v/t$, we have
\[
p_{\alpha,L}\#
\begin{pmatrix}
0&f\\f^c&0
\end{pmatrix}=
\begin{pmatrix}
0&f+h/\alpha\\f^c+h^{-1}/\alpha&0
\end{pmatrix}
\]
while \eqref{eq:80} reads
\[
\d
\begin{pmatrix}
v&s\\t&-v
\end{pmatrix}
+\left[
\begin{pmatrix}
0&\alpha\d f\\\alpha\d f^c&0
\end{pmatrix},
\begin{pmatrix}
v&s\\t&-v
\end{pmatrix}
\right]=0.
\]
We now argue as on page~\pageref{page:darboux-system} to conclude
that
\begin{equation}
\label{eq:82}
\d h=\alpha h \d f^c h-\alpha\d f.
\end{equation}
Finally, set $g=h/\alpha=v/t\alpha$.  Since $g$ is homogeneous in the
entries of $\omega$, without loss of generality, we may take $v$ to
be $\R^n$-valued and $t\in\R$ or $\I\R$ according to whether
$\alpha\in\R$ or $\I\R$.  Either way, $t\alpha\in\R$ so that
$g:M\to\R^n$ and \eqref{eq:82} becomes the familiar Riccati equation
\[
\d g=\alpha^2 g\d f^c g-\d f
\]
while
\[
p_{\alpha,L}\#
\begin{pmatrix}
0&f\\ f^c&0
\end{pmatrix}=
\begin{pmatrix}
0&f+g\\f^c+(\alpha^2g)^{-1}&0.
\end{pmatrix}
\]
Thus
\[
p_{\alpha,L}\#\begin{pmatrix}
0&f\\ f^c&0
\end{pmatrix}=
\begin{pmatrix}
0&\D_{\alpha^2}f\\\D_{\alpha^2}f^c&0
\end{pmatrix}
\]
and we have proved:
\begin{thm}
\label{th:30}
The dressing action of the simple factor $p_{\alpha,L}$ on a \Ch\
pair $(f,f^c)$ is by the Darboux transform $\D^v_{\alpha^2}$ where
$L$ is the complexification of the null-line corresponding to
$\alpha(v-f(o))$.
\end{thm}
In particular, Darboux transforms $\D_r$ correspond to the two types
of simple factor according to the sign of $r$.

\begin{rem}
Our action on $\p$-flat maps is only local: $p_{\alpha,L}\#\Phi$
fails to be defined at points $p\in M$ where
$\Phi(p)\not\in\mathcal{U}_{p_{\alpha,L}}$, that is, when
$\Phi^{-1}(p)(\alpha)L=\<v_0\>$ or $\<v_\infty\>$.  The geometric
meaning of this restriction is now clear: these are the points where
$g(p)=0$ or $g(p)=\infty$ and so are exactly the singularities of our
Riccati equation.  In the first case, we have $f(p)=\fh(p)$ and, in
the second, $f^c(p)=\fh^c(p)$.  In either case, we have genuine
singularities of the corresponding curved flats $(f,\fh)$ or
$(f^c,\fh^c)$.
\end{rem}

\subsection{Applications}
\label{sec:applications}

This new viewpoint on Darboux transformations allows several standard
arguments from the loop group formalism to be applied.  We conclude
our study by considering some of these.

\subsubsection{Explicit solutions}
\label{sec:explicit-solutions}

In general, computation of a Darboux transform involves solving a
differential equation: either the Riccati equation for $g$ or, what
is essentially the same thing, the \MC\ equations for the based
extended frame $\Phi$ at $\lambda=\alpha$.  However, the loop group
approach has the following advantage: if one based extended frame is
known then the based extended frame of any Darboux transform can be
found \emph{algebraically} via:
\[
p^{\vphantom{1}}_{\alpha,L}\#\Phi=p^{\vphantom{1}}_{\alpha,L}\Phi
p^{-1}_{\alpha,\Phi(\alpha)^{-1}L}.
\]
In this way, one can iteratively construct infinitely many explicit
examples given one known based extended frame---this is the procedure
of ``dressing the vacuum''.  The issue is, of course, to find a
suitable ``seed'' \Ch\ pair with known extended frame.

Experience with other problems (see, for example, \cite{BurPed95})
suggests that a good starting point is to look for surfaces framed by
a $2$-dimensional abelian subgroup of $G$ for then the \MC\ equations
are solved by exponentiation and extended flat frames are readily
computed.

For example: let $e_1,e_2,e_3$ denote the standard basis of $\R^3$
and let $f:\R^2\to\R^3$ be given by
\[
f(x,y)=xe_1+ye_2.
\]
The plane parametrised by $f$ is trivially isothermic with \Ch\
transform
\[
f^c(x,y)=xe_1-ye_2
\]
and the corresponding $\p$-flat map has
\[
\d\psi=E_1\d x+E_2\d y
\]
where
\[
E_1=
\begin{pmatrix}
0&e_1\\e_1&0
\end{pmatrix},\quad
\begin{pmatrix}
0&e_2\\-e_2&0.
\end{pmatrix}
\]
\begin{ex}
Show that $[E_1,E_2]=0$, $E_1^2=-1$ and $E_2^2=1$.
\end{ex}
Thus the extended flat frame $\Phi$ based at $0\in\R^2$ with
$\mc{\Phi}=\lambda\d\psi$ is given by
\begin{align*}
\Phi(x,y)(\lambda)&=(\exp \lambda xE_1)(\exp \lambda yE_2)\\
&=\bigl(\cos\lambda x+(\sin \lambda x)E_1\bigr)
\bigl(\cosh\lambda y+(\sinh \lambda y)E_1\bigr).
\end{align*}
Having got our hands on $\Phi$, we can compute the $T$-transforms of
$f$:
\[
T_r f=\Phi_{\sqrt{r}}\cdot0.
\]
\begin{ex}
Show that
\[
(\T_rf)(x,y)=\frac{(\sin 2\sqrt{r}x)e_1+(\sinh2\sqrt{r}y)e_2}%
{2(\cos^2\sqrt{r}x+\sinh^2\sqrt{r}y)}
\]
\end{ex}
Thus the $T$-transforms of $f$ are different parametrisations of the
same plane as is to be expected as all these isothermic surfaces
share the same solution $\kappa\equiv0$ of Calapso's equation.
\begin{ex}
Compute the Darboux transforms of $f$.
\end{ex}
A discussion of this example and its Darboux transforms can be found
in \cite{Cie97A}.

A somewhat less trivial example arises as follows: set
\begin{align*}
E_3&=
\begin{pmatrix}
0&e_3\\e_3&0
\end{pmatrix}\\
E=E_\k+E_\p=
\begin{pmatrix}
e_1e_2&e_1\\-e_1&e_1e_2
\end{pmatrix}
\end{align*}
and observe that $[E,E_3]=0$.  Taking $\k$ and $\p$ components gives
(since $E_3\in\p$)
\[
[E_\k,E_3]=[E_\p,E_3]=0
\]
so that
\[
(E_\k+\lambda E_\p)\d x+E_3\d y
\]
solves the \MC\ equations for all $\lambda$ and so integrates to give
a frame $\hat{\Phi}$ of an associated family of curved flats.
Indeed, since
\[
E_3^2=-1,\quad (E_\k+\lambda E_\p)^2=\lambda^2-1,
\]
we readily compute that
\begin{align*}
\hat{\Phi}(x,y)(\lambda)&=\exp\bigl(x(E_\k+\lambda E_\p)+\lambda y
E_3\bigr)\\&=
\bigl(\cosh x\sqrt{\lambda^2-1}+
\dfrac{\sinh x\sqrt{\lambda^2-1}}{\sqrt{\lambda^2-1}}(E_\k+\lambda
E_\p)\bigr)
\bigl(\cos \lambda y+(\sin\lambda y)E_3\bigr).
\end{align*}
Now $\hat{\Phi}$ is not an extended flat frame since
$\mc{\hat{\Phi}}$ has non-zero $\k$-component but the analysis of
Section~\ref{sec:curv-flats-symm} assures us that gauging by
$\hat{\Phi}_{|\lambda=0}$ gives such a frame.  Thus we define $\Phi$
by $\Phi=\hat{\Phi}\hat{\Phi}^{-1}_{|\lambda=0}$ to get a based
extended flat frame with $\Phi(x,y)(\lambda)$ given by
\begin{multline*}
\bigl(\cosh x\sqrt{\lambda^2-1}+
\dfrac{\sinh x\sqrt{\lambda^2-1}}{\sqrt{\lambda^2-1}}(E_\k+\lambda
E_\p)\bigr)\\
\bigl(\cos \lambda y+(\sin\lambda y)E_3\bigr)
\bigl(\cos x-(\sin x)E_\k\bigr).
\end{multline*}

\begin{ex}
\begin{enumerate}
\item Show that
\[
\dl{\Phi}=
\begin{pmatrix}
0&f\\f^c&0
\end{pmatrix}
\]
where
\begin{align*}
f(x,y)&=\half(\sin 2x)e_1+\half(1-\cos 2x)e_2+ye_3\\
f^c(x,y)&=-\half(\sin 2x)e_1-\half(1-\cos 2x)e_2+ye_3
\end{align*}
so that the \Ch\ pair associated to $\Phi$ is a right cylinder of
radius $\half$ (and so $H\equiv 1$) together with (up to a
translation) the parallel (that is, identical) cylinder parametrised
by $f+N$.
\item Compute the $T$-transforms of the cylinder.
\item Compute the Darboux transforms of the cylinder.
\item Persuade a computer to draw pictures of the surfaces you have
found.
\end{enumerate}
\end{ex}
A detailed analysis of this example and its Darboux transforms, using
somewhat different methods, has been carried out by Bernstein \cite{Ber99}.

\subsubsection{Bianchi permutability}
\label{sec:bianchi-perm}

Recall the assertion of Theorem~\ref{th:15}: given an isothermic
surface $f$ and Darboux transforms $f_i=\D_{r_i}f$, $i=1,2$, there is
a fourth isothermic surface $\fh$ such that
\[
\fh=\D_{r_1}f_2=\D_{r_2}f_1.
\]
Moreover, Theorem~\ref{th:16} says that the \Ch\ transform of such a
\Bq\ is another such so that
\[
\fh^c=\D_{r_1}f_2^c=\D_{r_2}f_1^c.
\]
In view of Theorem~\ref{th:30}, both these results can be formulated
in terms of simple factors: given a $\p$-flat map $\psi$ and Darboux
transforms $\psi_1=p_{\alpha_1,L_1}\#\psi$,
$\psi_2=p_{\alpha_2,L_2}\#\psi$, there is a $\p$-flat map $\hat{\psi}$
and light-lines $L_1'$, $L_2'$ such that
\[
\hat{\psi}=p_{\alpha_1,L_1'}\#(p_{\alpha_2,L_2\vphantom{'}}\#\psi)=
p_{\alpha_2,L_2'}\#(p_{\alpha_1,L_1\vphantom{'}}\#\psi),
\]
that is,
\[
(p_{\alpha_1,L_1'}p_{\alpha_2,L_2\vphantom{'}})\#\psi=
(p_{\alpha_2,L_2'}p_{\alpha_1,L_1\vphantom{'}})\#\psi.
\]
We shall therefore have found an alternative (and simultaneous!)
proof of both the Bianchi Permutability Theorem~\ref{th:15} and its
\Ch\ transform Theorem~\ref{th:16} as soon as we establish:
\begin{prop}
\label{th:31}
Let $p_{\alpha_i,L_i}\in\Gmb$, $i=1,2$, with
$\alpha_1^2\neq\alpha_2^2$.

Set
\begin{align*}
L'_1&=p_{\alpha_2,L_2}(\alpha_1)L_1\\
L'_2&=p_{\alpha_1,L_1}(\alpha_2)L_2
\end{align*}
and assume that $L'_i\neq\<v_0\>^\C,\<v_\infty\>^\C$, $i=1,2$.

Then $p_{\alpha_1,L_i'}\in\Gmb$, $i=1,2$ and
\begin{equation}
\label{eq:83}
p_{\alpha_1,L_1'}p_{\alpha_2,L_2\vphantom{'}}=
p_{\alpha_2,L_2'}p_{\alpha_1,L_1\vphantom{'}}.
\end{equation}
\end{prop}
\begin{proof}
Since $\alpha_1\neq\pm\alpha_2$, we have that $p_{\alpha_2,L_2}^{-1}$
is holomorphic near $\alpha_2$ and so we may apply
Proposition~\ref{th:28} with $E=p_{\alpha_2,L_2}^{-1}$ to conclude
that $p_{\alpha_1,L_1'}\in\Gmb$ and, further, that
\[
p^{\vphantom{1}}_{\alpha_1,L_1} p_{\alpha_2,L_2}^{-1}p_{\alpha_1,L_1'}^{-1}
\]
is holomorphic and invertible at $\pm\alpha_1$.

Similarly $p_{\alpha_2,L_2'}\in\Gmb$ and
\[
p^{\vphantom{1}}_{\alpha_2,L_2} p_{\alpha_1,L_1}^{-1}p_{\alpha_2,L_2'}^{-1}
\]
is holomorphic and invertible at $\pm\alpha_2$.

Now contemplate
\[
p^{\vphantom{1}}_{\alpha_1,L_1'}
(p^{\vphantom{1}}_{\alpha_2,L_2} p_{\alpha_1,L_1}^{-1}p_{\alpha_2,L_2'}^{-1})=
(p^{\vphantom{1}}_{\alpha_1,L_1} p_{\alpha_2,L_2}^{-1}p_{\alpha_1,L_1'}^{-1})^{-1}
p_{\alpha_2,L_2'}^{-1}.
\]
Looking at the left hand side, we see that this expression is
holomorphic at $\pm\alpha_2$ and, from the right hand side, we see
that is is holomorphic at $\pm\alpha_1$.  Thus it is holomorphic on
$\P^1$ and so constant.  Evaluating at $\lambda=0$ now gives
\[
p^{\vphantom{1}}_{\alpha_1,L_1'}
(p^{\vphantom{1}}_{\alpha_2,L_2} p_{\alpha_1,L_1}^{-1}p_{\alpha_2,L_2'}^{-1})=1
\]
that is
\[
p_{\alpha_1,L_1'}p_{\alpha_2,L_2\vphantom{'}}=
p_{\alpha_2,L_2'}p_{\alpha_1,L_1\vphantom{'}}.
\]
\end{proof}

There is another way to think about this result which shows what a
general phenomenon it is that we are dealing with here: \eqref{eq:83}
amounts to a factorisation
\[
p^{\vphantom{1}}_{\alpha_1,L_1}p_{\alpha_2,L_2}^{-1}=
p_{\alpha_2,L_2'}^{-1}p_{\alpha_1,L_1'}^{\vphantom{-1}}
\]
corresponding to the subgroups $\G_{\alpha_i}$ of $\Gmb$ consisting
of those $g\in\Gmb$ that are holomorphic on
$\P^1\setminus\{\pm\alpha_i\}$.  Just as before, we get from such a
factorisation a local action of $\G_{\alpha_1}$ on $\G_{\alpha_2}$
which we denote by $*_{\alpha_1}$ and then
\[
p^{\vphantom{1}}_{\alpha_2,L'_2}=
(p^{\vphantom{1}}_{\alpha_1,L_1}*_{\alpha_1}p^{-1}_{\alpha_2,L_2})^{-1}.
\]
More generally, for $g_i\in\G_{\alpha_i}$, we can find
$g'_i\in\G_{\alpha_i}$ with
\[
g'_1g_2^{\vphantom{1}}=g'_2g_1^{\vphantom{1}}
\]
by setting
\[
g'_2=(g^{\vphantom{1}}_1*_{\alpha_1}g_2^{-1})^{-1},\quad
g'_1=(g^{\vphantom{1}}_2*_{\alpha_2}g_1^{-1})^{-1}.
\]
This shows that Bianchi permutability is not a consequence of the
fact that our simple factors have simple poles but rather that these
factors have only two poles.

In auspicious circumstances (for example $\alpha\in\I\R$ and $G$
compact) one can argue as in Section~\ref{sec:simple-factors} and
precompose everything with $t_{\alpha_2}^{-1}$ to reduce $*_{\alpha_1}$
to the globally defined Pressley--Segal action.  For example, with
$G=\mathrm{SU}(2)$, this accounts for the classical B\"acklund
transform of pseudo-spherical surfaces and their Bianchi
permutability.

As a final advertisement for this technology, let us give another
proof of Theorem~\ref{th:18} which asserts that the Darboux transform
of a \Bq\ is another \Bq\ thus giving a configuration of 8 isothermic
surfaces forming the vertices of a cube all of whose faces are \Bq
s.  For this, choose $\alpha_1,\alpha_2,\alpha_3\in\C^\times$ with
all $\alpha_i^2$ real and distinct and let $q_i\in\G_{\alpha_i}$ be
three simple factors with poles at $\pm\alpha_i$.
Proposition~\ref{th:31} now gives simple factors
$q_i^j\in\G_{\alpha_i}$ with
\begin{subequations}
\begin{align}
\label{eq:84}q_3^1q_1^{\vphantom{1}}&=q^3_1q_3^{\vphantom{1}}\\
\label{eq:85}q^2_1q_2^{\vphantom{1}}&=q^1_2q_1^{\vphantom{1}}\\
\label{eq:86}q^3_2q_3^{\vphantom{1}}&=q^2_3q_2^{\vphantom{1}}
\end{align}
\end{subequations}
and then simple factors $q^{jj}_i\in\G_{\alpha_i}$ with
\begin{subequations}
\begin{align}
\label{eq:87}q_1^{33}q_2^3&=q_2^{33}q_1^3\\
\label{eq:88}q_1^{22}q_3^2&=q_3^{22}q_1^2\\
\label{eq:89}q_2^{11}q_3^1&=q_3^{11}q_2^1.
\end{align}
\end{subequations}
(This notation becomes a little easier to stomach when one sees that
the subscripts locate the poles of the simple factor.
Figure~\ref{fig:Bq} on page \pageref{fig:Bq} may also help.)

The key to our result is the following lemma that asserts that the
$q^{jj}_i$ are determined solely by their poles:
\begin{lem}\label{th:32}
$q_1^{22}=q_1^{33}$, $q_2^{11}=q_2^{33}$, $q_3^{11}=q_3^{22}$.
\end{lem}
\begin{proof}
Multiply \eqref{eq:87} by $q_3$ to get
\[
q_1^{33}q_2^3q_3^{\vphantom{1}}=q_2^{33}q_1^3q_3^{\vphantom{1}}
\]
and use \eqref{eq:84} and \eqref{eq:86} to get
\[
q_1^{33}q_3^2q_2^{\vphantom{1}}=q^{33}_2q_3^1q_1^{\vphantom{1}}.
\]
Rearranging this and using \eqref{eq:85} yields
\[
q_1^{33}q_3^2=q^{33}_2q_3^1q_1^{\vphantom{1}}q_2^{-1}=
q^{33}_2q_3^1(q_2^1)^{-1}q^2_1
\]
whence
\begin{equation}
\label{eq:90}
q_1^{33}q_3^2(q^2_1)^{-1}=q^{33}_2q_3^1(q_2^1)^{-1}.
\end{equation}
Temporarily denote by $q$ the common value in \eqref{eq:90}.  From
the left hand side, we see that $q$ is holomorphic except possibly at
$\pm\alpha_1,\pm\alpha_3$ while the right hand side tells us that $q$
is holomorphic except possibly at $\pm\alpha_2,\pm\alpha_3$.  We
therefore conclude that $q$ has poles at $\pm\alpha_3$ only, that is,
$q\in\G_{\alpha_3}$ so that we have factorisations
\begin{align*}
q_1^{33}q_3^2&=q q^2_1\\
q^{33}_2q_3^1&=q q_2^1.
\end{align*}
However, for $i\neq j$, $\G_{\alpha_i}\cap\G_{\alpha_j}=\{1\}$ so
factorisations of this kind are unique (recall Exercise~\ref{ex:9}!)
and, comparing with \eqref{eq:88}, \eqref{eq:89}, we get
\begin{align*}
q_1^{33}&=q_1^{22}& q&=q_3^{22}\\
q_2^{33}&=q_2^{11}& q&=q_3^{11}.
\end{align*}
\end{proof}

With this in hand, start with a $\p$-flat map $\psi$ and set
\[
\psi_1=q_1\#\psi,\quad \psi_2=q_2\#\psi,\quad \psi'=q_3\#\psi.
\]
We then obtain \Bq s $(\psi,\psi_1,\hat{\psi},\psi_2)$,
$(\psi,\psi',\psi'_1,\psi_1^{\vphantom{1}})$,
$(\psi,\psi',\psi'_2,\psi^{\vphantom{1}}_2)$ with
\begin{align*}
\hat{\psi}&=(q^2_1q_2^{\vphantom{1}})\#\psi=(q^1_2q_1^{\vphantom{1}})\#\psi\\
\psi_1'&=(q^3_1q_3^{\vphantom{1}})\#\psi=(q^1_3q_1^{\vphantom{1}})\#\psi\\
\psi_2'&=(q^3_2q_3^{\vphantom{1}})\#\psi=(q^2_3q_2^{\vphantom{1}})\#\psi
\end{align*}
and then a \Bq\ $(\psi',\psi'_1,\hat{\psi}',\psi_2')$ with
\[
\hat{\psi}'=(q_1^{33}q^3_2)\#\psi'=(q^{33}_2q^3_1)\#\psi'.
\]
The situation is summarised in Figure~\ref{fig:Bq}.
\begin{figure}[ht]
\begin{center}
\includegraphics{feb-fig-3.mps}
\caption{Darboux transform of a Bianchi quadrilateral}
\label{fig:Bq}
\end{center}
\end{figure}
The claim is that the remaining two faces
$(\psi_2^{\vphantom{1}},\hat{\psi},\hat{\psi}',\psi_2')$ and
$(\psi_1^{\vphantom{1}},\psi_1'\hat{\psi}',\hat{\psi})$ are also \Bq
s. That is,
\[
\hat{\psi}'=(q_1^{22}q^2_3)\#\psi_2^{\vphantom{1}}=
(q_2^{11}q^1_3)\#\psi_1^{\vphantom{1}}.
\]
But Lemma~\ref{th:32} with \eqref{eq:86} gives
\begin{align*}
(q_1^{22}q^2_3)\#\psi_2&=(q_1^{33}q^2_3q_2^{\vphantom{1}})\#\psi\\
&=(q_1^{33}q^3_2q_3^{\vphantom{1}})\#\psi=(q_1^{33}q^3_2)\#\psi'=\hat{\psi}'.
\end{align*}
A similar argument establishes the second equation.

While this argument requires some book-keeping it seems less involved
than our Clifford algebra cross-ratio argument of
Section~\ref{sec:bianchi-perm-cliff} and has a certain universal
character which applies to all other B\"acklund transforms which are
given by the dressing action of simple factors.  For example, working
with $G=\mathrm{SU}(2)$ and the extended frames of pseudo-spherical
surfaces, we immediately read off a result which was doubtless known
to Bianchi:
\begin{thm}
The B\"acklund transform of a \Bq\ of pseudo-spherical surfaces is
another such.
\end{thm}

\section{Coda}
\label{sec:coda}

We have developed a fairly complete theory of isothermic surfaces in
$\R^n$ but there is more to be said and more to be understood.  I
draw this (already over-long) work to a close by indicating some
recent developments in the area and some open problems.

\subsection{Recent developments }
\label{sec:recent-developments}

\subsubsection{Symmetric $R$-spaces}
\label{sec:symmetric-r-spaces}

The conformal geometry of $S^n$ is an example of a parabolic geometry
of a kind possessed by any \emph{symmetric $R$-space}
\cite{KobNag64,KobNag65,Tak65}.  According to Nagano \cite{Nag65},
these can be characterised as those Riemannian symmetric spaces of
compact type which admit a Lie groups of diffeomorphisms strictly
larger than the isometry group.  Thus examples include:
\begin{enumerate}
\item $S^n$ with its group $\Mob$ of conformal diffeomorphisms and
more generally the conformal compactification $S^p\times S^q$ of
$\R^{p,q}$ with the corresponding group of conformal diffeomorphisms;
\item Any Grassmannian $G_k(\R^n)$ of $k$-planes in $\R^n$ with the
action of $\mathrm{PSL}(n,\R)$.  In particular, taking $k=1$, we find
the setting of projective differential geometry.
\item Any Hermitian symmetric space of compact type with its group of
biholomorphisms.
\end{enumerate}
All symmetric $R$-spaces have a common algebraic
structure\footnote{The stabilisers of points in the ``big'' group are
parabolic subgroups with \emph{abelian} nilradical.} which accounts
for all the structure we have exploited in this work: one has
analogues of stereographic projection, the pseudo-Riemannian
symmetric space $Z$ of point pairs and, most importantly, an
invariant formulation of the notion of an isothermic submanifold.
\Ch, Darboux and $T$-transformations are all available in this
general context and the delicate inter-relations between them remain
true as does the loop group interpretation described in
section~\ref{sec:loop-groups-backlund}.

In particular, these ideas provide a manifestly conformally invariant
definition of an isothermic surface in $S^n$, the lack of which may
be viewed as a weakness of the present work.

These ideas will be described in \cite{BurPedPin}.

\subsubsection{Meromorphic functions as isothermic surfaces}
\label{sec:merom-funct-as}

One can specialise our existing theory to the case $n=2$: this amounts
to studying meromorphic functions on a Riemann surface $M$.  In this
case, the isothermic surface condition is vacuous---any meromorphic
function is isothermic---so one must change one's point of view and
emphasis the role of the holomorphic quadratic differential $Q$.
Thus, on a polarised Riemann surface $(M,Q)$, the \Ch\ transform
$f^c$ of a meromorphic function $f$ is given by specialising
\eqref{eq:28} to this setting and demanding
\[
\del f^c=Q/\del f.
\]
If $f$ is viewed as the Gauss map of a minimal surface with Hopf
differential $Q$ via the Weierstrass--Enneper formula, this
transformation gives rise to an intriguing transformation of minimal
surfaces that has been studied by McCune \cite{McC}.

\subsubsection{Willmore surfaces in $S^4$}
\label{sec:willmore-surfaces-s4}

In low dimensions, Clifford algebras are most conveniently studied as
transformations of spinors.  For $n=4$, this amounts to viewing $S^4$
as the quaternionic projective line $\H P^1$.  Here a central topic
is the study of Willmore surfaces---extremals of a conformally
invariant functional that are characterised by the harmonicity of
their conformal Gauss map.  There are strong formal analogies between
such conformal Gauss maps and the Euclidean Gauss maps of CMC
surfaces.  One can exploit this analogy along with the methods of
Section~\ref{sec:darb-transf-gener} to obtain a large family of
``Darboux'' transformations of Willmore surfaces.

Similarly, another class of transformations can be obtained by
adapting the methods of McCune \cite{McC} to this context.

A detailed exposition of these ideas may be found in \cite{BurFerLes00}.

\subsection{Open problems}
\label{sec:open-problems}

I list some problems to which I would like to know the answers!
\begin{enumerate}
\item Is there any interesting theory of isothermic submanifolds of
$\R^n$ of dimension greater than two?  The problem here is to find a
suitable definition that is not too restrictive: certainly our
formulation only works in $2$ dimensions and the same is true of the
symmetric $R$-space approach.  One way forward might be to study
submanifolds admitting a conformal Ribaucour sphere congruence.  The
work of Dajczer--Tojeiro \cite{DajToj} may be relevant here.
\item Motivated by considerations concerning surfaces isometric to
quadrics that this writer does not understand, Darboux \cite{Dar99A}
distinguished the class of \emph{special isothermic surfaces} in
$\R^3$ and these were studied intensively by Bianchi
\cite{Bia05,Bia05A} and Calapso \cite{Cal15}.  Characterised by a
differential equation on the mean curvature, this class includes CMC
surfaces as a degenerate case and is stable under all the
transformations of the theory.
\begin{prob}
Find a simple geometric characterisation of special isothermic
surfaces in $\R^3$.

Is there an interesting extension of the notion to surfaces in $\R^n$?
\end{prob}
\item The theory of constant mean curvature (CMC) surfaces in $\R^3$
lies at the intersection of two integrable geometries: via their
Gauss maps, they are the same as harmonic maps into $S^2$, a
well-studied integrable system.  In particular, they admit a spectral
deformation, the ``associated family'', through CMC surfaces $f_\mu$
for $\mu\in S^1$.  On the other hand, viewed as isothermic surfaces,
they have the spectral deformation $\T_r f$ through isothermic
surfaces for $r\in\R$ which amounts to the Guichard--Lawson
deformation through CMC surfaces in other space forms (see
\cite{HerMusNic} for a recent account).  The relation between these
deformations is not well-understood although there is some evidence
to suggest that they should be viewed as the angular and radial parts
of a single complex deformation\footnote{\textit{Note added in
December 2001}: this issue has now been clarified in \cite{BurPedPin01}}.

Again, the Darboux transforms of CMC surfaces described herein amount
to the (iterated) B\"acklund transforms of the harmonic map theory
\cite{HerPed97} despite the fact that the underlying symmetry groups
seem quite different.  Thus we formulate:
\begin{prob}
Find a theory of CMC surfaces that unifies the harmonic map and
isothermic surface theories. 
\end{prob}
\end{enumerate}

\providecommand{\bysame}{\leavevmode\hbox to3em{\hrulefill}\thinspace}

\end{document}